\documentclass[a4paper,fleqn]{cas-sc}

\usepackage[authoryear,longnamesfirst]{natbib}

\usepackage{amscd}
\usepackage{tikz-cd}
\usetikzlibrary{cd}
\usepackage[all,cmtip]{xy}

\def\tsc#1{\csdef{#1}{\textsc{\lowercase{#1}}\xspace}}
\tsc{WGM}
\tsc{QE}

\newtheorem{assum}{Assumption}
\newtheorem{theorem}{Theorem}[section]
\newtheorem{lemma}[theorem]{Lemma}

\newtheorem{example}[theorem]{Example}

\newtheorem{remark}[theorem]{Remark}
\newproof{proof}{Proof}
\newtheorem{alg}{Algorithm}

\newcommand{\dd}{\,{\rm d}}

\newcommand{\bs}{\boldsymbol}

\newcommand{\dt}{\Delta t}
\newcommand{\red}{\color{red}}

\newcommand{\curl}{{\rm curl\,}}
\renewcommand{\div}{\operatorname{div}}
\newcommand{\grad}{{\rm grad\,}}

\begin{document}
\let\WriteBookmarks\relax
\def\floatpagepagefraction{1}
\def\textpagefraction{.001}

\shorttitle{MFEM for Rosensweig FHD model}    

\shortauthors{W. Wu and X. Xie}  

\title [mode = title]{Energy-stable mixed finite element methods for the Rosensweig ferrofluid flow model}  

\tnotemark[1] 

\tnotetext[1]{This work was supported by the National Natural Science Foundation of China (11971094, 12171340).} 

%

\author[1]{Yongke Wu}


\fnmark[1]

\ead{wuyongke1982@sina.com}



\affiliation[1]{organization={School of Mathematical Sciences, University of Electronic Science and Technology of China},
            city={Chengdu},
            postcode={611731}, 
            state={Sichuan},
            country={China}}

\author[2]{Xiaoping Xie}
\cormark[1]

\fnmark[2]

\ead{xpxie@scu.edu.cn}



\affiliation[2]{organization={School of Mathematical, Sichuan University},
            city={Chengdu},
            postcode={610065}, 
            state={Sichuan},
            country={China}}

\cortext[1]{Corresponding author}

\fntext[1]{ }


\begin{abstract}
In this paper, we consider  mixed finite element semi-/full discretizations  of the Rosensweig ferrofluid flow model.We first establish some regularity results for the model under  several basic assumptions. 
 Then we show that the energy stability of the weak solutions   is preserved exactly for both the semi- discrete and fully discrete finite element solutions. Moreover, we prove the existence and uniqueness of the discrete solutions. We also    derive optimal error estimates for the discrete schemes. Finally, we provide numerical experiments to  verify the theoretical results.
\end{abstract}




\maketitle

\section{Introduction}\label{sec:intr}

Ferrofluids, also known as ferromagnetic fluids  or magnetic liquids, are stable colloidal suspensions of nanoparticles  within aqueous or nonaqueous liquids that possess strong magnetic features. Such fluids allow to be manipulated and controlled when exposed to   external magnetic fields and have a wide range of applications~(cf. \citet{Pankhurst2003,Zahn2001}) such as magnetic resonance imaging, magneto-optic sensors, magnetic separation, vibration damping and acoustics, and vacuum technology.

There are two widely accepted ferrohydrodynamics (FHD) models,  i.e. the Shliomis FHD model (cf. \citet{Shliomis1972Effective,Shliomis2002Ferrofluids}) and  the Rosensweig FHD model (cf. \citet{Rosensweig1985,Rosensweig1987}). They both treat ferrofluids as homogeneous monophase fluids, but differ in dealing with internal rotation of magnetic nanoparticles. The Shliomis FHD model treats the rotation as a magnetic torque and is a coupled nonlinear system of the Navier-Stokes equations and the Maxwell equations, while the  Rosensweig   model considers the internal rotation by introducing an   angular momentum variable and is a coupled nonlinear system of the Navier-Stokes equations, the Maxwell equations and an angular momentum equation. 
We refer to \citet{Amirat2008,Amirat2009Strong,Amirat2010UniqueR,Amirat2008GlobalR,Nochetto2019} for several works on the
    existence of global-in-time weak solutions or local-in-time strong solutions for these two FHD models.
    
    For the the Shliomis FHD  model and related problems, there are considerable research efforts on the   finite element analysis.  
    \citet{Nochetto2016-two} 
    proposed a first-order energy-stable fully discrete finite element  scheme   for a two-phase ferrofluid flow model which couples  the Shliomis   model and the Cahn-Hilliard equations. 
    \citet{Zhang.G;2021} developed an unconditionally energy-stable decoupled fully finite element   scheme for a two-phase ferrofluid flow model. 
     \citet{Wu-Xie2023} presented  and analyzed mixed finite element methods for  a steady   Shliomis FHD model   with magnetization  paralleled to the magnetic field.   \citet{WuXie2023}  also  proposed   energy-stable mixed finite element methods for the unsteady Shliomis   model    and derived   error estimates for both the semi- and fully discrete schemes under some regularity assumptions.  Recently \citet{Mao2024} 
      investigated a first order energy-stable fully discrete finite element method for a simplified Shliomis   model with the effective magnetizing field  being treated as a known function. They  proved some regularity results for the  model   and  obtained   error estimates for the  numerical solutions.

As for the  finite element  analysis of the Rosensweig FHD model, there are very limited research works   in the literature.    \citet{Nochetto2015The}   presented a first-order energy-stable fully discrete finite element  scheme for the  model and  analyzed the convergence of the discrete solution, by using compactness argument, for a simplified Rosensweig FHD model with a negligible demagnetizing field.

In this paper, we 
first establish an energy stability estimate and  several regularity results for the weak solutions of the Rosensweig FHD  model under some assumptions on the external magnetic field, the initial values, and the domain $\Omega$. Then we develop   energy-stable semi-discrete and  fully discrete  mixed finite element methods for the  model.   We show the existence and, under certain conditions, the uniqueness of the   semi-discrete  solution and the fully discrete   solution, and  derive optimal error estimates.  Finally, we construct a decoupled quasi-Newton  scheme to solve the nonlinear fully discrete scheme and give some numerical examples.

The rest of this paper is organized as follows. In section 2 we introduce notations as well as  the governing equations and the weak form of the Rosensweig FHD model. In Section 3 we derive the energy estimate and some regularity results.  We devote Sections 4 and 5   to the analyses of the semi-discrete and fully discrete mixed finite element schemes, respectively, including   the energy stability, the existence/uniqueness   and  error estimates of  the discrete solutions.   In section 6 we apply the quasi-Newton iterative algorithm to perform several numerical tests.

\section{Preliminaries} \label{sec:pre}

\subsection{Notations}
Let $\Omega \subset \mathbb R^3$ be a bounded simply connected convex domain with Lipschitz boundary $\partial\Omega$, and $T>0$ be the final time. We set $\Omega_T := \Omega \times (0,T]$ and $\Gamma_T := \partial\Omega \times (0,T]$. Let $\bs n$ be the unit outward norm vector along $\partial\Omega$.

For any nonnegative integer $l$, we denote by $H^l(\Omega)$ the usual $l$-th order Sobolev space with norm $\|\cdot\|_l$ and semi-norm $|\cdot|_l$. In particular, $H^0(\Omega) = L^2(\Omega)$ denotes the space of all square integrable functions on $\Omega$, with the inner product $(\cdot,\cdot)$ and the norm $\|\cdot\|$. For vector spaces $(H^l(\Omega))^3$ and $(L^2(\Omega))^3$, we use the same notations of norm, semi-norm and inner product as those for the scalar cases. We also introduce the spaces
\begin{align*}
H(\curl) = \{\bs v = & (v_1,v_2,v_3)^\intercal\in (L^2(\Omega))^3:\\
&\curl\bs v:= (\partial_y v_3- \partial_z v_2,\partial_z v_1 - \partial_x v_3,\partial_x v_2 - \partial_y v_1)^\intercal \in (L^2(\Omega))^3\}
\end{align*}
and
$$
H(\div) = \{\bs v =  (v_1,v_2,v_3)^\intercal\in (L^2(\Omega))^3:~\div\bs v := \partial_x v_1 + \partial_y v_2 + \partial_z v_3 \in L^2(\Omega)\},
$$
with the norms
$$
\|\bs v\|_c:=(\|\bs v\|^2 + \|\curl\bs v\|^2)^{1/2},\quad\text{and}\quad \|\bs v\|_d: = (\|\bs v\|^2 + \|\div\bs v\|^2)^{1/2}.
$$
For simplicity we set
\begin{align*}
& \bs S:=(H_0^1(\Omega))^3 = \{\bs v \in (H^1(\Omega))^3: ~\bs v = 0\text{ on }\partial\Omega\}, \\
& \bs U: = H_0(\curl) =\{ \bs v \in H(\curl):~\bs v \times \bs n = 0 \text{ on }\partial\Omega\}, \\
& \bs V := H_0(\div) = \{\bs v \in H(\div):~\bs v\cdot\bs n = 0\text{ on }\partial\Omega\}, \\
& W := L_0^2(\Omega) = \{ v \in L^2(\Omega): \int_\Omega v\dd \bs x = 0\}, \\
& H(\div 0) := \{\bs v \in H(\div):~\div \bs v = 0\},\\
& H_0(\div 0): = H_0(\div)\cap H(\div 0), \\
& \bs L_n^2(\Omega) := \{\bs v \in (L^2(\Omega))^3:~\bs v \cdot\bs n = 0\text{ on }\partial\Omega\}.
\end{align*}
For any $\bs v \in H_0(\curl)\cap H(\div)$ or $H(\curl)\cap H_0(\div)$, we define
$$
\|\bs v\|_n: = (\|\curl\bs v\|^2 + \|\div\bs v\|^2)^{1/2}.
$$

For any scalar- or vector-valued space $X$, defined on $\Omega$, with norm $\|\cdot\|_X$, we set
$$
W^{l,p}([0,T];X):=\{v:~[0,T] \rightarrow X;~\|v\|_{W^{l,p}(X)} <\infty\},
$$
with 
$$
\|v\|_{W^{l,p}(X)} := \left\{ 
\begin{array}{ll}
\left(\int_0^T\sum\limits_{k = 0}^l\|\partial_t^kv(\cdot,t)\|_X^p \dd t\right)^{1/p} & \text{if } 1\leq p <\infty, \\
\max\limits_{0\leq k\leq l} \text{ess}\sup\limits_{0\leq t \leq T} \|\partial_t^kv(\cdot,t)\|_X & \text{if }p = \infty. 
\end{array}
\right.
$$
For simplicity we set $W^{l,p}(X) := W^{l,p}([0,T];X)$. For any integer $r \geq 0$, the spaces $H^r(X):=H^r([0,T];X)$ and $\mathcal C^r(X) := \mathcal C^r([0,T];X)$ can be defined similarly.

Let $\mathcal T_h$ be a quasi-uniform shape regular tetrahedron triangulation of $\Omega$ with mesh size $h:= \max_{K\in \mathcal T_h} h_K$, with $h_K$ the diameter of the element $K\in \mathcal T_h$. For any integer $l\geq 0$, let $\mathbb P_l(K)$ be the set of polynomials, defined on $K$, of degree no more than $l$.

Let  $N$ be a positive integer, and let
 $$
 \mathcal J_{\dt} = \{ t_n:~t_n = n\dt,~0\leq n\leq N\}
 $$
 be a uniform partition of the time interval $[0,T]$ with the time step
size $\dt = T/ N$.

Throughout this paper,  we denote by $C$ (or $C$ with a subscript) a generic    positive constant, independent of the  mesh sizes $h$ and $\dt$, and $C$ may be different at its each occurrence. Sometimes for simplicity we may use $a\lesssim (or \gtrsim) b$ to denote $a\leq (or \geq) Cb$.

\subsection{Governing equations of the Rosensweig  model}

Consider the flow of an incompressible, viscous Newtonian ferrofluid, filling $\Omega$, under the action of a known external magnetic field $\bs H_e$ satisfying 
$$
\bs H_e\cdot\bs n = 0\qquad\text{on }\Gamma_T.
$$
The magnetic field $\bs H_e$ induces a demagnetizing field $\bs H$ and a magnetic induction $\bs B$ satisfying the law
$$
\bs B =\mu_0( \bs H+\bs m),
$$
with $\bs m$ the magnetization inside $\Omega$. 

In the Rosensweig model  (cf. \citet{Rosensweig1985,Rosensweig1987})  the   velocity $\bs u$, the pressure $p$, the angular momentum $\bs\omega$, and the magnetization $\bs m$ are governed by the conservation of linear momentum, the incompressible condition, the conservation of angular momentum, the transport of magnetization, and the magneto static equation for $\bs H$, and the governing equations read as follows:  
\begin{align}
\label{eq:FHD-1}
& \rho(\partial_t\bs u + (\bs u\cdot\nabla)\bs u) - (\eta+\zeta)\Delta \bs u +\nabla p = \mu_0(\bs m\cdot\nabla)\bs H + 2\zeta \curl\bs\omega &\text{in }\Omega_T,\\
\label{eq:FHD-2}& \div\bs u = 0 & \text{in }\Omega_T,\\
\label{eq:FHD-3}&\rho\kappa (\partial_t\bs\omega + (\bs u\cdot\nabla)\bs\omega) -\eta^\prime\Delta \bs\omega - (\eta^\prime + \lambda^\prime)\nabla\div\bs\omega = \mu_0\bs m\times \bs H 
+ 2\zeta(\curl\bs u - 2\bs \omega) & \text{in }\Omega_T,\\
\label{eq:FHD-4} & \partial_t\bs m + (\bs u\cdot\nabla)\bs m - \sigma\Delta\bs m = \bs\omega \times \bs m - \frac{1}{\tau}(\bs m - \chi_0\bs H) & \text{in }\Omega_T,\\
\label{eq:FHD-5} &\curl \bs H = 0 & \text{in }\Omega_T,\\
\label{eq:FHD-6} & \mu_0\div(\bs H+\bs m) = -\div\bs H_e& \text{in }\Omega_T,
\end{align}
subject to the boundary conditions
\begin{equation}\label{eq:bd}
\bs u = 0,\quad\bs\omega = 0,\quad \bs m\cdot\bs n = 0,\quad\curl\bs m\times\bs n = 0,\quad\text{and}\quad \bs H\cdot\bs n = 0\qquad\text{on }\Gamma_T
\end{equation}
and the initial conditions
\begin{equation}\label{eq:initial}
\bs u(\bs x,0) = \bs u_0(\bs x),\quad\bs \omega(\bs x,0) = \bs \omega_0(\bs x),\quad\bs m(\bs x,0) = \bs m_0(\bs x)\qquad\forall~\bs x \in \Omega.
\end{equation}
Here $\bs u_0,\ \bs \omega_0,\ \bs m_0$ are given functions, and the parameters $\rho,\ \eta,\ \zeta,\ \mu_0,\ \kappa,\ \eta^\prime,\ \lambda^\prime,\ \sigma,\ \tau$ and $\chi_0$ are positive constants and their physical means can be found in, for example, \citet{Shliomis1972Effective,Rosensweig1985,Rosensweig1987,Shliomis2002Ferrofluids}.

Introduce new variables
\begin{equation}
\label{eq:kcurlm}
\bs z :=\bs u \times\bs m,\quad \bs k := \curl\bs m, \quad \tilde p := p - \frac{\mu_0}{2}\bs m\cdot\bs H
\end{equation}
and let $\phi:~[0,T]\rightarrow H_0^1(\Omega)$ be such that $\bs H = \nabla \phi$ due to the fact that $\curl\bs H = 0$.  Then we apply  the identities
\begin{equation*}
\left\{
\begin{array}{ll}
(\bs u\cdot\nabla)\bs m & = \frac{1}{2}\bs u\div\bs m - \frac{1}{2}\bs m\div\bs u + \frac{1}{2}\nabla (\bs u\cdot\bs m) - \frac{1}{2}\bs m \times \curl\bs u 
- \frac{1}{2}\bs u \times \curl\bs m -\frac{1}{2}\curl(\bs u \times \bs m),\\ \\
(\bs m \cdot\nabla)\bs H,\bs v) & = -\frac{1}{2}(\bs H\cdot\bs v,\div\bs m) - \frac{1}{2}(\bs m\cdot\bs H,\div\bs v) + \frac{1}{2}(\bs v\cdot\bs m,\div\bs H) 
+\frac{1}{2}(\bs v \times \curl\bs m,\bs H) + \frac{1}{2}(\bs m\times \curl\bs v,\bs H)\\
& \qquad\qquad\qquad\qquad\qquad\qquad\qquad\qquad\qquad\qquad\qquad\qquad\qquad\qquad\qquad\qquad\quad \forall~\bs v \in \bs H_0^1(\Omega)
\end{array}
\right.
\end{equation*}
and get the following weak problem: Find $\bs u:[0,T]\rightarrow \bs S$, $\tilde p:[0,T] \rightarrow W$, $\bs\omega:[0,T]\rightarrow \bs S$, $\bs m:[0,T]\rightarrow \bs V$, $\bs z:[0,T]\rightarrow \bs U$, $\bs k:[0,T]\rightarrow \bs U$, $\bs H:[0,T]\rightarrow \bs V$  and $\phi:[0,T]\rightarrow W$ such that
\begin{align}
\label{eq:FHD-weak-1}
\begin{split}
 \rho(\partial_t\bs u,\bs v) + \rho b(\bs u,\bs u,\bs v) + \eta(\nabla\bs u,\nabla\bs v) + \zeta (\curl\bs u,\curl\bs v)
  - (\tilde p,\div\bs v)
  - \mu_0 c(\bs v,\bs m,\bs H) & \\
  - \frac{\mu_0}{2}(\bs v \times \bs k,\bs H) 
  - \frac{\mu_0}{2}(\bs m\times \curl\bs v,\bs H) 
 - 2\zeta(\bs\omega,\curl\bs v) &= 0
\end{split} & \forall ~\bs v \in \bs S,\\
\label{eq:FHD-weak-2}
 (\div\bs u,q) & = 0 &\forall~q\in W, \\
\label{eq:FHD-weak-3}
\begin{split}
\rho\kappa(\partial_t\bs\omega,\bs s) + \eta^\prime(\nabla\bs\omega,\nabla\bs s) + (\eta^\prime + \lambda^\prime)(\div\bs\omega,\div\bs s) 
+ \rho\kappa b(\bs u,\bs\omega,\bs s) -\mu_0(\bs m\times \bs H,\bs s) &\\
- 2\zeta(\curl\bs u,\bs s) 
+ 4\zeta(\bs \omega,\bs s)& = 0
\end{split} & \forall~\bs s \in \bs S,\\
\label{eq:FHD-weak-4}
\begin{split}
(\partial_t\bs m,\bs F) - c(\bs u,\bs m,\bs F) +\sigma(\div\bs m,\div\bs F) + \sigma(\curl\bs k,\bs F) 
-\frac{1}{2}(\bs m\times\curl\bs u,\bs F) &\\
- \frac{1}{2}(\bs u \times \bs k,\bs F) 
- \frac{1}{2}(\curl\bs z,\bs F)
 - (\bs\omega\times \bs m,\bs F) + \frac{1}{\tau}(\bs m,\bs F) 
-\frac{\chi_0}{\tau}(\bs H,\bs F) & = 0
\end{split}& \forall~\bs F\in \bs V,\\
\label{eq:FHD-weak-5} (\bs z,\bs\Lambda) - (\bs u\times \bs m,\bs\Lambda) & = 0 & \forall~\bs \Lambda \in \bs U, \\
\label{eq:FHD-weak-6} (\bs k,\bs \Theta) - (\bs m,\curl\bs\Theta) & = 0 & \forall~\bs \Theta \in \bs U, \\
\label{eq:FHD-weak-7} (\bs H,\bs G) + (\phi,\div\bs G) & = 0 &\forall~\bs G \in \bs V, \\
\label{eq:FHD-weak-8}\mu_0(\div\bs H+\div\bs m,r) + (\div\bs H_e,r) & = 0 &\forall~ r \in W,
\end{align}
with the initial data \eqref{eq:initial}, where 
\begin{align*}
& b(\bs u,\bs v,\bs w) := \frac{1}{2}\left[(\bs u\cdot \nabla)\bs v,\bs w) - ((\bs u\cdot\nabla)\bs w,\bs v)\right]\quad\text{and}\quad c(\bs v,\bs m,\bs F) := \frac{1}{2}\left[(\bs m\cdot\bs v,\div\bs F) - (\bs F\cdot\bs v,\div\bs m) \right].
\end{align*}

\subsection{Basic assumptions}
For the sake of  stability and regularity analysis,  we make several assumptions on the prescribed data and the solution of the model problem \eqref{eq:FHD-1} - \eqref{eq:initial}.
\begin{assum}[Regularity assumption on  external magnetic field]\label{ass:He}
The external magnetic field $\bs H_e$ satisfies
$$
\bs H_e \in  H^1( (H^2(\Omega))^3 )\cap  H^2((L^2(\Omega))^3 )\cap L^\infty(( H^2(\Omega))^3).
$$
\end{assum}
\begin{assum}[Regularity assumptions on initial data]\label{ass:initial-data}
The initial data $\bs u_0 , \bs\omega_0 $ and $\bs m_0 $ satisfy
{$$\bs u_0 \in \bs S\cap H_0(\div 0)\cap (H^2(\Omega))^3, \quad \bs\omega_0\in \bs S \cap (H^2(\Omega))^3 , \quad \bs m_0\in \bs V \cap (H^2(\Omega))^3.$$}
\end{assum}

\begin{assum}[Boundedness assumptions on solution]\label{ass:solu} 
The solution components  $ \bs u$ and $\bs m$ of the problem \eqref{eq:FHD-1} - \eqref{eq:initial} satisfies 
$$
\int_0^T(\|\bs u\|_1^4 + \|\bs m\|_n^4)\dd t \leq C.
$$ 
\end{assum}

To derive higher regularity of the weak solution we need the following smoothness assumption on the domain $\Omega$:
\begin{assum}[Smoothness assumption on   domain]\label{ass:omega}
The boundary of $\Omega$ is smooth so that 
the unique solution $(\bs v,q)$ of the steady Stokes problem
$$
-\Delta\bs v + \nabla q  = \bs f,\qquad\div\bs v = 0\quad\text{in }\Omega,\qquad \bs v = 0\quad\text{on }\partial\Omega,
$$
the unique solution $\bs\psi \in \bs H_0^1(\Omega)$ of the equation
$$
-\Delta\bs\psi -(1+\frac{\lambda'}{\eta'})\grad\div\bs\psi = \bs g\quad\text{in }\Omega,\qquad\bs\psi = 0\quad\text{on }\partial\Omega,
$$
and the unique solution $\bs \phi \in \bs H^1(\Omega)$ of the equation
$$
-\Delta\bs\phi = \bs j\quad \text{in }\Omega,\qquad\curl\bs\phi\times\bs n = 0,\qquad\bs\phi\cdot\bs n = 0\quad\text{on }\partial\Omega,
$$
for $\bs f,~\bs g,~\bs j \in \bs L^2(\Omega)$ satisfy the following regularity estimates:
\begin{align*}
\|\bs v\|_2 + \|q\|_1 & \lesssim \|\bs f\|, \\
\|\bs\psi\|_2 & \lesssim \|\bs g\|, \\
\|\bs\phi\|_2 & \lesssim \|\bs j\|.
\end{align*}
\end{assum}

\section{Regularity estimates}\label{sec:reg}

\subsection{Energy stability}
Let $\bs u \in \bs S$, $\tilde p\in W$, $\bs\omega \in \bs S$, $\bs m\in  \bs V$, $\bs z\in \bs U$, $\bs k \in \bs U$, $\bs H \in \bs V$ and $\phi \in W$ solve the weak problem \eqref{eq:FHD-weak-1} - \eqref{eq:FHD-weak-8}. We define for any $t \in [0,T]$   the energy
 $$
 \mathcal E(t): = \rho\|\bs u(\cdot,t)\|^2 + \rho\kappa\|\bs\omega(\cdot,t)\|^2 + \|\bs m(\cdot,t)\|^2 + \mu_0\|\bs H(\cdot,t)\|^2
 $$
 and the dissipated energy
 \begin{align*}
 \mathcal F(t) := & \eta\|\nabla\bs u(\cdot,t)\|^2 + \eta^\prime \|\nabla\bs\omega(\cdot,t)\|^2 + (\eta^\prime+\lambda^\prime)\|\div\bs\omega(\cdot,t)\|^2 + \frac{1}{\tau}\|\bs m(\cdot,t)\|^2 + \frac{\chi_0+(1+\mu_0)\chi_0}{\tau}\|\bs H(\cdot,t)\|^2
 \\
 & + \sigma\|\div\bs m(\cdot,t)\|^2  + \sigma\|\bs k(\cdot,t)\|^2 + \mu_0\sigma\|\div\bs H(\cdot,t)\|^2  +\zeta\|\curl\bs u(\cdot,t) - 2\bs\omega(\cdot,t)\|^2.
 \end{align*}
 We cite a useful inequality for the forthcoming energy stability analysis.
 \begin{lemma}(\citep[Lemma 1]{Kirby2015} ))  \label{lem:Gr}
 Let  $x\geq 0$ satisfy 
 $$
 x^2 \leq \gamma^2 + \beta x
 $$
 for $\beta$, $\gamma\geq 0$. Then  
 $
 x \leq \beta + \gamma.
 $
 \end{lemma}
 We have the following energy stability result.
 \begin{theorem}
 \label{the:ener}
Under   Assumption \ref{ass:He}   the energy inequality
 $$
 \mathcal E(t) +C_1\int_0^t\mathcal F(s) \dd s \leq \mathcal E(0)  + C_2\int_0^t (\|\bs H_e(\cdot,s)\|_{H(\div)}^2  + \|\partial_t\bs H_e(\cdot,s)\|^2)\dd s
 $$
 holds for all $t \in [0,T]$, where $C_1$ and $C_2$ are positive constants depending only on $\chi_0$, $\mu_0$ and $\tau$.
 \end{theorem}
 \begin{proof}
 Taking $\bs v = \bs u$ in \eqref{eq:FHD-weak-1} and using \eqref{eq:FHD-weak-2} with $q =\tilde p$, we get
 \begin{equation}\label{eq:ener-pr-1}
 \begin{split}
 \frac{\rho}{2}\frac{\dd}{\dd t}\|\bs u\|^2 + \eta\|\nabla\bs u\|^2 + \zeta\|\curl\bs u\|^2-2\zeta(\bs\omega,\curl\bs u) 
 =  \mu_0 c(\bs u,\bs m,\bs H) 
  + \frac{\mu_0}{2}(\bs u\times \bs k,\bs H)  + \frac{\mu_0}{2}(\bs m \times \curl\bs u,\bs H).
 \end{split}
\end{equation}
Taking $\bs s = \bs\omega$ in \eqref{eq:FHD-weak-3}, we have
\begin{equation}\label{eq:ener-pr-2}
\begin{split}
\frac{\rho\kappa}{2}\frac{\dd}{\dd t}\|\bs\omega\|^2 +\eta^\prime\|\nabla\bs\omega\|^2 + (\eta^\prime+\lambda^\prime)\|\div\bs\omega\|^2 - 2\zeta(\curl\bs u,\bs \omega) + 4\zeta\|\bs\omega\|^2  
= \mu_0(\bs m\times\bs H,\bs \omega).
\end{split}
\end{equation}
Taking $\bs F = \bs m$ in \eqref{eq:FHD-weak-4} and using the fact that $(\bs z,\bs k) = (\bs u \times \bs m,\bs k) = (\bs m,\curl\bs z)$ and $\|\bs k\|^2 = (\curl\bs m,\bs k)$, we obtain
\begin{equation}\label{eq:ener-pr-3}
\frac{1}{2}\frac{\dd}{\dd t}\|\bs m\|^2 +\sigma\|\div\bs m\|^2 + \sigma\|\bs k\|^2
+ \frac{1}{\tau}\|\bs m\|^2 - \frac{\chi_0}{\tau}(\bs H,\bs m) = 0.
\end{equation}
Taking $\bs G = \bs H$ and $\bs G = \bs m$ in \eqref{eq:FHD-weak-7}, respectively,  and taking $r = \phi$ in \eqref{eq:FHD-weak-8}, we arrive at
\begin{align*}
& \|\bs H\|^2 =- (\phi,\div\bs H),\\
& (\bs H,\bs m) = -(\phi,\div\bs m),\\
& (\div\bs H,\phi) + (\div\bs m,\phi) = -\frac{1}{\mu_0}(\div\bs H_e,\phi) = \frac{1}{\mu_0}(\bs H,\bs H_e),
\end{align*}
which imply
$$
(\bs H,\bs m) = -(\phi,\div\bs m) = -\|\bs H\|^2 - \frac{1}{\mu_0}(\bs H,\bs H_e).
$$
Thus, \eqref{eq:ener-pr-3} becomes
\begin{equation*}\label{eq:ener-pr-4}
\frac{1}{2}\frac{\dd}{\dd t}\|\bs m\|^2 + \sigma\|\div\bs m\|^2 + \sigma\|\bs k\|^2
+ \frac{1}{\tau}\|\bs m\|^2 + \frac{\chi_0}{\tau}\|\bs H\|^2  = -\frac{\chi_0}{\tau\mu_0}(\bs H,\bs H_e).
\end{equation*}
This relation, together with \eqref{eq:ener-pr-1}-\eqref{eq:ener-pr-2}, yields%
%
%
\begin{equation}\label{eq:ener-pr-5}
\begin{split}
& \frac{1}{2}\frac{\dd}{\dd t}\left(\rho\|\bs u\|^2 + \rho\kappa\|\bs\omega\|^2 + \|\bs m\|^2 \right) + \eta\|\nabla\bs u\|^2 +\eta^\prime\|\nabla\bs\omega\|^2+\sigma\|\bs k\|^2  + \frac{1}{\tau}\|\bs m\|^2 + \sigma\|\div\bs m\|^2 \\
& \quad + (\eta^\prime +\lambda^\prime)\|\div\bs\omega\|^2 
  +\zeta\|\curl\bs u - 2\bs\omega\|^2 + \frac{\chi_0}{\tau}\|\bs H\|^2\\
= & \mu_0c(\bs u,\bs m,\bs H) + \frac{\mu_0}{2}(\bs u\times \bs k,\bs H) + \frac{\mu_0}{2}(\bs m\times \curl\bs u,\bs H)  + \mu_0(\bs m\times \bs H,\bs\omega) - \frac{\chi_0}{\tau\mu_0}(\bs H,\bs H_e).
\end{split}
\end{equation}

Differentiating \eqref{eq:FHD-weak-8} with respect to $t$ and taking $r = \phi$ in the resulting equation, we get
\begin{align}\label{eq:FHD-weak-008}
(\div\partial_t\bs H,\phi) + (\div\partial_t\bs m,\phi)  = -\frac{1}{\mu_0}(\div\partial_t\bs H_e,\phi).
\end{align}
Taking $\bs G = \partial_t\bs H$, $\bs G = \partial_t\bs m$ and $\bs G = \partial_t\bs H_e$ in \eqref{eq:FHD-weak-7}, respectively,  we have
\begin{align}\label{eq:FHD-weak-007}
\frac{1}{2}\frac{\dd}{\dd t}\|\bs H\|^2 = -(\phi,\div\partial_t\bs H),\quad (\bs H,\partial_t\bs m) = -(\phi,\div\partial_t\bs m),\quad (\bs H,\partial_t\bs H_e) = -(\phi,\div\partial_t\bs H_e).
\end{align}
Taking $\bs F = \bs H$ in \eqref{eq:FHD-weak-4} and using the fact that $(\curl\bs k,\bs H) = (\phi,\div\curl\bs k) = 0$, we obtain
\begin{equation*}\label{eq:ener-pr-6}
\begin{split}
(\partial_t\bs m,\bs H) - c(\bs u,\bs m,\bs H) + \sigma(\div\bs m,\div\bs H)  - \frac{1}{2}(\bs m\times \curl\bs u,\bs H) 
- \frac{1}{2}(\bs u \times \bs k,\bs H) 
- \frac{1}{2} (\curl\bs z,\bs H)&\\- (\bs \omega \times \bs m,\bs H) + \frac{1}{\tau}(\bs m,\bs H) 
- \frac{\chi_0}{\tau}\|\bs H\|^2 & = 0,
\end{split}
\end{equation*}
which plus \eqref{eq:FHD-weak-007} and \eqref{eq:FHD-weak-008} gives
\begin{equation}
\label{eq:ener-pr-7}
\begin{split}
\frac{1}{2}\frac{\dd}{\dd t}&\|\bs H\|^2  = (\div\partial_t\bs m,\phi) + \frac{1}{\mu_0}(\div\partial_t\bs H_e,\phi)  = -(\partial_t\bs m,\bs H) - \frac{1}{\mu_0}(\partial_t\bs H_e,\bs H)\\
& = -\frac{1}{\mu_0}(\partial_t\bs H_e,\bs H) - c(\bs u,\bs m,\bs H) + \sigma(\div\bs m,\div\bs H)    - \frac{1}{2}(\bs m\times \curl\bs u,\bs H)  - \frac{1}{2}(\bs u\times \bs k,\bs H) \\
&\quad- \frac{1}{2}(\curl\bs z,\bs H) - (\bs \omega\times \bs m,\bs H)   - \frac{1+\chi_0}{\tau}\|\bs H\|^2 -\frac{1}{\tau\mu_0}(\bs H_e,\bs H).
\end{split}
\end{equation}
Multiplying \eqref{eq:ener-pr-7} with $\mu_0$, adding the resulting equation to \eqref{eq:ener-pr-5},  and using the fact that $\div\bs H= -\frac{1}{\mu_0}\div\bs H_e- \div\bs m$, we get
\begin{equation}
\label{eq:ener-pr-8}
\begin{split}
\frac{1}{2}\frac{\dd}{\dd t}\mathcal E(t) + \mathcal F(t) = & -\frac{\chi_0+\mu_0}{\tau\mu_0}(\bs H_e,\bs H) - (\partial_t\bs H_e,\bs H) - \sigma(\div\bs H_e,\div\bs m),
\end{split}
\end{equation}
where we have used the fact   $(\curl\bs z,\bs H) = 0$, 
which is obtained by taking $\bs G = \curl\bs z$ in \eqref{eq:FHD-weak-7} and using   $\div\curl = 0$. For any $t \in [0,T]$, integrating \eqref{eq:ener-pr-8} with respect to $t$ in the interval $(0,t)$, we finally obtain the desired conclusion. 
 \end{proof}
 \begin{remark}
 The energy estimate in Theorem \ref{the:ener} implies that under   Assumption \ref{ass:He}  the   solution of problem  \eqref{eq:FHD-1} - \eqref{eq:initial}  satisfies
 \begin{align*}
 & \bs u \in L^\infty((L^2(\Omega))^3)\cap L^2((H^1(\Omega))^3), & \bs \omega \in L^\infty((L^2(\Omega))^3)\cap L^2((H^1(\Omega))^3), \\
 &\bs m \in L^\infty((L^2(\Omega))^3) \cap L^2((H^1(\Omega))^3), & \bs H \in L^\infty((L^2(\Omega))^3)\cap L^2((H^1(\Omega))^3).
 \end{align*}
 \end{remark}
 
 \subsection{Regularity estimates}
 
To discuss the regularity of the  solution of  the   problem  \eqref{eq:FHD-1} - \eqref{eq:initial}, we need to introduce the Stokes operator $\mathscr A_1$ and the Maxwell type operator $\mathscr A_2$ with
 \begin{eqnarray*}
  \mathscr A_1:&~\mathcal D(\mathscr A_1) :=  (H^2(\Omega))^3 \cap \bs S \cap  H(\div 0) \rightarrow  H_0(\div 0),& \mathscr A_1 = -\mathscr P_1(-\Delta),\\
  \mathscr A_2:&~\mathcal D(\mathscr A_2) :=  (H^2(\Omega))^3\cap H_0(\div) \rightarrow \bs L_n^2(\Omega),& \mathscr A_2 = \mathscr P_2(\curl\curl-\grad\div),
 \end{eqnarray*}
where $\mathscr P_1$ and  $\mathscr P_2$ are  the  $L^2 $-projection operators  from $( L^2(\Omega))^3$ to $H_0(\div 0)$ and  from $( L^2(\Omega))^3$ to $\bs L_n^2(\Omega)$, respectively. There hold the following inequalities (cf. \citet{Heywood1982,Gerbeau2006}):
\begin{align}
\label{eq:A-i-1} & \|\bs v\|_{L^\infty} + \|\bs v\|_{W^{1,3}} \leq C\|\bs v\|_1^{1/2}\|\mathscr A_i\bs v\|^{1/2},& \forall ~\bs v \in \mathcal D(\mathscr A_i), \\
\label{eq:A-i-2} & \|\bs v\|_2 \leq C\|\mathscr A_i\bs v\|,&\forall ~\bs v \in \mathcal D(\mathscr A_i),
\end{align}
for $i=1,2$. 

\begin{theorem}\label{the:reg-2}
Under Assumptions \ref{ass:He}
-\ref{ass:solu} there holds
\begin{equation}
\label{eq:reg-2}
\begin{split}
\sup\limits_{0\leq t\leq T}(\|\nabla\bs u(t)\|^2 + \|\nabla\bs\omega(t)\|^2 + \|\bs m(t)\|_n^2 + \|\bs H(t)\|^2_n )
+ \int_0^T(\|\bs u\|_2^2 + \|\bs \omega\|_2^2 + \|\bs m\|_2^2 + \|\bs H\|_2^2)\dd t & \leq C. 
\end{split} 
\end{equation}
\end{theorem}
\begin{proof}
Testing equation \eqref{eq:FHD-1} with $\mathscr A_1(\bs u)$, \eqref{eq:FHD-3} with $-\Delta \bs\omega$, and \eqref{eq:FHD-4} with $\mathscr A_2(\bs m)$ and $\mathscr A_2(\bs H)$, respectively, we have
\begin{equation}\label{eq:reg-pr-1}
\begin{split}
\frac{\rho}{2}\frac{\dd}{\dd t}\|\nabla\bs u\|^2 + \eta\|\mathscr A_1(\bs u)\|^2 + \zeta \|\curl\curl\bs u\|^2 = \mu_0((\bs m\cdot\nabla)\bs H,\mathscr A_1(\bs u)) 
+ 2\zeta (\curl\bs \omega,\curl\curl\bs u) - \rho ((\bs u\cdot\nabla)\bs u,\mathscr A_1(\bs u)),
\end{split}
\end{equation}
\begin{equation}\label{eq:reg-pr-2}
\begin{split}
\frac{\rho\kappa}{2}\frac{\dd }{\dd t}\|\nabla\bs\omega\|^2 +\eta^\prime \|\Delta\bs\omega\|^2 + (\eta^\prime + \lambda^\prime)\|\nabla\div\bs\omega\|^2 = \mu_0(\bs m\times \bs H,-\Delta\bs\omega) 
+ 2\zeta (\curl\curl\bs u,\curl\bs \omega) \\- 4\zeta (\|\curl\bs\omega\|^2 + \|\div\bs\omega\|^2),
\end{split}
\end{equation}
\begin{equation}\label{eq:reg-pr-3}
\begin{split}
\frac{1}{2}\frac{\dd}{\dd t}\|\bs m\|_{n}^2 +\sigma\|\mathscr A_2(\bs m)\|^2 + \frac{1}{\tau}\|\bs m\|_n^2 = (\bs\omega \times \bs m,\mathscr A_2(\bs m)) 
+\frac{\chi_0}{\tau}(\bs H,\mathscr A_2(\bs m)) -((\bs u\cdot\nabla)\bs m,\mathscr A_2(\bs m)),
\end{split}
\end{equation}
and
\begin{equation}\label{eq:reg-pr-4}
\begin{split}
(\partial_t\bs m,\mathscr A_2(\bs H)) + ((\bs u\cdot\nabla)\bs m,\mathscr A_2(\bs H)) - \sigma(\Delta\bs m,\mathscr A_2(\bs H)) =(\bs\omega\times \bs m,\mathscr A_2(\bs H)) 
-\frac{1}{\tau}(\bs m,\mathscr A_2(\bs H))+ \frac{\chi_0}{\tau}\|\div\bs H\|^2 .
\end{split}
\end{equation}
Equation \eqref{eq:FHD-6} implies 
$$
\mu_0( \div\bs H_t +\div\bs m_t ) = -\div\partial_t\bs H_e,
$$
testing which with $\div\bs H$ gives
\begin{equation}\label{eq:reg-pr-5}
\frac{\mu_0}{2}\frac{\dd}{\dd t}\|\div\bs H\|^2 = \mu_0(\partial_t\bs m,\mathscr A_2(\bs H)) + (\partial_t\bs H_e,\mathscr A_2(\bs H)).
\end{equation}
Testing   \eqref{eq:FHD-6} with $\div\bs m$, we get
\begin{equation}
\label{eq:reg-pr-6}
\mu_0(\mathscr A_2(\bs H),\bs m) + \mu_0\|\div\bs m\|^2 = - (\div\bs H_e,\div\bs m).
\end{equation}
By   \eqref{eq:FHD-6}, we also have
\begin{equation}\label{eq:reg-pr-7}
\begin{split}
\sigma\mu_0(\Delta\bs m,\mathscr A_2(\bs H)) & = \sigma\mu_0(\grad\div\bs m,\mathscr A_2(\bs H)) 
& = -\mu_0\sigma\|\mathscr A_2(\bs H)\|^2 + \sigma(\grad\div\bs H_e,\mathscr A_2(\bs H))
\end{split}
\end{equation}
and
\begin{equation}
\begin{split}\label{eq:reg-pr-8}
(\bs H,\mathscr A_2(\bs m)) & = (\bs H,-\grad\div\bs m) = (\div\bs H,\div\bs m)  = -\|\div\bs m\|^2 -\frac{1}{\mu_0}(\div\bs H_e,\div\bs m).
\end{split}
\end{equation}

Adding equations \eqref{eq:reg-pr-1} - \eqref{eq:reg-pr-3} and \eqref{eq:reg-pr-5} together and using \eqref{eq:reg-pr-4} and \eqref{eq:reg-pr-6}-\eqref{eq:reg-pr-8}, we get
\begin{equation}\label{eq:reg-pr-9}
\begin{split}
\frac{1}{2}\frac{\dd}{\dd t}\left(\rho\|\nabla\bs u\|^2 + \rho\kappa\|\nabla\bs\omega\|^2 + \|\bs m\|_n^2 + \mu_0\|\bs H\|_n^2 \right) + \eta\|\mathscr A_1(\bs u)\|^2 +\eta^\prime \|\Delta\bs\omega\|^2 
+ (\eta^\prime + \lambda^\prime)\|\nabla\div\bs\omega\|^2 + \sigma\|\mathscr A_2(\bs m)\|^2\\ 
+\sigma\mu_0\|\mathscr A_2(\bs H)\|^2+ \frac{1}{\tau}\|\bs m\|_n^2 
+ \frac{1+\chi_0}{\tau}\|\div\bs m\|^2+ \frac{\mu_0\chi_0}{\tau}\|\div\bs H\|^2 
+4\zeta\|\div\bs\omega\|^2 
+\zeta\|\curl\curl\bs u - \curl\bs\omega\|^2\\
 = \mu_0((\bs m\cdot\nabla)\bs H,\mathscr A_1(\bs u)) 
- \rho((\bs u\cdot\nabla)\bs u,\mathscr A_1(\bs u)) + \mu_0(\bs m\times \bs H,-\Delta\bs\omega) 
+ (\bs\omega\times \bs m,\mathscr A_2(\bs m)) 
-((\bs u\cdot\nabla)\bs m,\mathscr A_2(\bs m)) \\
+ \mu_0(\bs\omega\times \bs m,\mathscr A_2(\bs H)) 
- \mu_0((\bs u\cdot\nabla)\bs m,\mathscr A_2(\bs H)) 
+(\partial_t\bs H_e,\mathscr A_2(\bs H)) \\
-\frac{\mu_0+\chi_0}{\tau\mu_0}(\div\bs H_e,\div\bs m) + \sigma(\grad\div\bs H_e,\mathscr A_2(\bs H)).
\end{split}
\end{equation}
Using \eqref{eq:A-i-1}, \eqref{eq:A-i-2}, the H\"older inequality, and the Young inequality, we have
\begin{align*}
|\mu_0((\bs m\cdot\nabla)\bs H,\mathscr A_1(\bs u))| & \leq\frac{\eta}{4} \|\mathscr A_1(\bs u)\|^2 + \frac{\sigma\mu_0}{10} \|\mathscr A_2(\bs H)\|^2 + C\|\bs m\|_n^2 \|\div\bs H\|^2, \\
|\rho((\bs u\cdot\nabla)\bs u,\mathscr A_1(\bs u))| & \leq \frac{\eta}{4}\|\mathscr A_1(\bs u)\|^2 + C\|\nabla\bs u\|^6,\\
|\mu_0(\bs m\times \bs H,-\Delta\bs\omega) |& \leq \frac{\eta^\prime}{2}\|\Delta \bs\omega\|^2 + C\|\bs m\|_n^2 \|\div\bs H\|^2,\\
|(\bs\omega\times \bs m,\mathscr A_2(\bs m))| & \leq \frac{\sigma}{4} \|\mathscr A_2(\bs m)\|^2 + C\|\nabla\bs\omega\|^2\|\bs m\|\|\bs m\|_n,\\
|(\bs u\cdot\nabla)\bs m,\mathscr A_2(\bs m))| & \leq \frac{\sigma}{4} \|\mathscr A_2(\bs m)\|^2 + C\|\nabla\bs u\|^4\|\bs m\|_n^2,\\
|\mu_0(\bs\omega\times \bs m,\mathscr A_2(\bs H))|& \leq \frac{\sigma\mu_0}{10}\|\mathscr A_2(\bs H)\|^2 + C\|\nabla\bs\omega\|^2 \|\bs m\|\|\bs m\|_n,\\
|\mu_0((\bs u\cdot\nabla)\bs m,\mathscr A_2(\bs H))|& \leq \frac{\sigma\mu_0}{10} \|\mathscr A_2(\bs H)\|^2 + C\|\nabla\bs u\|^4\|\bs m\|_n^2,\\
|(\partial_t\bs H_e,\mathscr A_2(\bs H))| & \leq \frac{\sigma\mu_0}{10} \|\mathscr A_2(\bs H)\|^2 + C\|\partial_t\bs H_e\|^2,\\
\frac{\mu_0+\chi_0}{\tau\mu_0}|(\div\bs H_e,\div\bs m)| & \leq \frac{1+\chi_0}{2\tau}\|\div\bs m\|^2 + C\|\div\bs H_e\|^2 ,\\
|\sigma (\grad\div\bs H_e,\mathscr A_2(\bs H))| & \leq \frac{\sigma\mu_0}{10} \|\mathscr A_2(\bs H)\|^2 + C\|\grad\div\bs H_e\|^2.
\end{align*}
Combining the above inequalities with \eqref{eq:reg-pr-9}, we obtain
\begin{equation}\label{eq:reg-pr-10}
\begin{split}
\frac{\dd}{\dd t}&\left(\rho\|\nabla\bs u\|^2 + \rho\kappa\|\nabla\bs\omega\|^2 + \|\bs m\|_n^2 + \mu_0\|\bs H\|_n^2 \right) + \eta\|\mathscr A_1(\bs u)\|^2 +\eta^\prime \|\Delta\bs\omega\|^2 
+ 2(\eta^\prime + \lambda^\prime)\|\nabla\div\bs\omega\|^2  
+ \sigma\|\mathscr A_2(\bs m)\|^2 \\
&+\sigma\mu_0\|\mathscr A_2(\bs H)\|^2+ \frac{2}{\tau}\|\bs m\|_n^2 
+ \frac{1+\chi_0}{\tau}\|\div\bs m\|^2+ \frac{2\mu_0\chi_0}{\tau}\|\div\bs H\|^2 
+8\zeta\|\div\bs\omega\|^2 
+2\zeta\|\curl\curl\bs u - \curl\bs\omega\|^2\\
\leq C &\left(\|\partial_t\bs H_e\|^2 + \|\div\bs H_e\|^2 +\|\grad\div\bs H_e\|^2 \right) + C\|\bs m\|_n^2 \|\bs H\|_n^2 
+ C\|\nabla\bs u\|^6 + C\|\nabla\bs\omega\|^2\|\bs m\|_n^2 +C\|\nabla\bs u\|^4\|\bs m\|_n^2.
\end{split}
\end{equation}
Integrating the inequality \eqref{eq:reg-pr-10} with respect to   $t$ from $0$ to $t_F$ for any $t_F\in (0,T]$, we have
\begin{equation}\label{eq:reg-pr-1010}
\begin{split}
&\rho\|\nabla\bs u(t_F)\|^2 +\rho\kappa \|\nabla\bs\omega(t_F)\|^2 + \|\bs m(t_F)\|_n^2 +\mu_0\|\bs H(t_F)\|_n^2 + \int_0^{t_F}(\frac{2}{\tau}\|\bs m\|_n^2 + \frac{2\mu_0\chi_0}{\tau}\|\bs H\|_n^2 )\dd t\\
&\qquad +\int_0^{t_F}(\eta\|\mathscr A_1(\bs u)\|^2 + \eta^\prime \|\Delta \bs\omega\|^2 +\sigma\|\mathscr A_2(\bs m)\|^2 + \sigma\mu_0\|\mathscr A_2(\bs H)\|^2)\dd t  \\
\leq & \rho\|\nabla\bs u_0\|^2 +\rho\kappa\|\nabla\bs\omega_0\|^2 + \|\bs m_0\|_n^2 + \mu_0\|\bs H(0)\|_n^2   +C\int_0^{T}\left(\|\partial_t\bs H_e\|^2 + \|\div\bs H_e\|^2 +\|\grad\div\bs H_e\|^2 \right)\dd t\\
&\qquad + C\int_0^{t_F} (\|\nabla\bs u\|^4 + \|\bs m\|_n^2 )(\|\nabla\bs u\|^2 +\|\nabla\bs\omega\|^2+ \|\bs m\|_n^2 + \|\bs H\|^2_n)\dd t.
\end{split}
\end{equation}
Theorem \ref{the:ener} and Assumptions \ref{ass:He} - \ref{ass:solu} imply 
\begin{align*}
& \rho\|\nabla\bs u_0\|^2 +\rho\kappa\|\nabla\bs\omega_0\|^2 + \|\bs m_0\|_n^2 + \mu_0\|\bs H(0)\|_n^2  \leq C, \\
& \int_0^{T}\left(\|\partial_t\bs H_e\|^2 + \|\div\bs H_e\|^2 +\|\grad\div\bs H_e\|^2 \right)\dd t \leq C,\\
& \int_0^T(\|\nabla\bs u\|^4 + \|\bs m\|_n^2)\dd t \leq C.
\end{align*}
Finally, applying the Gr\"onwall inequality to \eqref{eq:reg-pr-1010} and using the inequality \eqref{eq:A-i-2}, we obtain the desired result \eqref{eq:reg-2}.
\end{proof}
\begin{remark}
Theorem \ref{the:reg-2} implies that under Assumptions \ref{ass:He}
-\ref{ass:solu} the solution of the FHD model  \eqref{eq:FHD-1} - \eqref{eq:initial}  satisfies 
$$
\bs u,~\bs\omega,~\bs m,~\bs H \in  L^\infty((H^1(\Omega))^3)\cap  L^2((H^2(\Omega))^3).
$$
\end{remark}

We shall also derive some estimates on the time derivative of the solution, which will be used later in the error analysis. 

Differentiating both sides of \eqref{eq:FHD-1}, \eqref{eq:FHD-3} and \eqref{eq:FHD-4} with respect to $t$, we get
\begin{equation}\label{eq:FHD-t-1}
\begin{split}
\rho(\bs u_{tt} + (\bs u_t\cdot\nabla)\bs u + (\bs u\cdot\nabla)\bs u_t) - (\eta+\zeta)\Delta\bs u_t +\nabla p_t = \mu_0(\bs m_t\cdot\nabla)\bs H 
+ \mu_0(\bs m\cdot\nabla)\bs H_t + 2\zeta\curl\bs \omega_t,
\end{split}
\end{equation}
\begin{equation}\label{eq:FHD-t-3}
\begin{split}
\rho\kappa(\bs\omega_{tt} + (\bs u_t\cdot\nabla)\bs \omega + (\bs u\cdot\nabla)\bs \omega_t) - \eta^\prime\Delta\bs\omega_t - (\eta^\prime + \lambda^\prime)\nabla\div\bs\omega_t = \mu_0\bs m_t\times \bs H 
+ \mu_0\bs m\times \bs H_t + 2\zeta(\curl\bs u_t -2\bs\omega_t),
\end{split}
\end{equation}
and
\begin{equation}\label{eq:FHD-t-4}
\begin{split}
\bs m_{tt} + (\bs u_t\cdot\nabla)\bs m + (\bs u\cdot\nabla)\bs m_t -\sigma\Delta\bs m_t = \bs\omega_t \times \bs m + \bs \omega\times\bs m_t
-\frac{1}{\tau}(\bs m_t -\chi_0\bs H_t).
\end{split}
\end{equation}
\begin{theorem}\label{the:reg-t}
Under  Assumptions \ref{ass:He} - \ref{ass:solu}, the  solution components $ \bs u,~\bs\omega,~\bs m,$ and $\bs H$ of problem \eqref{eq:FHD-1} - \eqref{eq:initial} satisfy
\begin{equation}
\begin{split}
\sup\limits_{0\leq t \leq T} (\|\bs u_t(t)\|^2 + \|\bs\omega_t(t)\|^2 + \|\bs m_t(t)\|^2 + \|\bs H_t(t)\|^2) 
+ \int_0^{T}(\|\nabla\bs u_t\|^2 + \|\nabla\bs\omega_t\|^2 + \|\bs m_t\|_n^2 + \|\bs H_t\|_n^2)\dd t&  \leq C.
\end{split}
\end{equation}

\end{theorem}
\begin{proof}
Testing \eqref{eq:FHD-t-1} - \eqref{eq:FHD-t-4} with $\bs u_t$, $\bs \omega_t$ and $\bs m_t$, respectively, we have
\begin{equation}\label{eq:reg-t-pr-1}
\begin{split}
\frac{\rho}{2}\frac{\dd}{\dd t}\|\bs u_t\|^2 +\eta\|\nabla\bs u_t\|^2 +\zeta\|\curl\bs u_t\|^2 & = \mu_0((\bs m_t\cdot\nabla)\bs H,\bs u_t) 
+ \mu_0((\bs m\cdot\nabla)\bs H_t,\bs u_t) & +2\zeta(\curl\bs u_t,\bs \omega_t) - \rho((\bs u_t\cdot\nabla)\bs u,\bs u_t),
\end{split}
\end{equation}
\begin{equation}
\label{eq:reg-t-pr-2}
\begin{split}
\frac{\rho\kappa}{2}\frac{\dd }{\dd t} \|\bs\omega_t\|^2 +\eta^\prime \|\nabla\bs\omega_t\|^2 + (\eta^\prime + \lambda^\prime)\|\div\bs\omega_t\|^2  = \mu_0(\bs m_t\times \bs H,\bs\omega_t) 
+ \mu_0(\bs m\times \bs H_t,\bs \omega_t) + 2\zeta(\curl\bs u_t,\bs\omega_t)  \\
- 4\zeta\|\bs\omega_t\|^2 -\rho\kappa ((\bs u_t\cdot\nabla)\bs\omega,\bs\omega_t),
\end{split}
\end{equation}
and
\begin{equation}
\label{eq:reg-t-pr-3}
\begin{split}
\frac{1}{2}\frac{\dd}{\dd t}\|\bs m_t\|^2 + \sigma\|\bs m_t\|_n^2  + \frac{1}{\tau}\|\bs m_t\|^2 = (\bs\omega_t\times \bs m,\bs m_t) + \frac{\chi_0}{\tau}(\bs H_t,\bs m_t) 
- ((\bs u_t\cdot\nabla)\bs m,\bs m_t) .
\end{split}
\end{equation}
Testing \eqref{eq:FHD-t-4} with $\bs H_t$, we get
\begin{equation}\label{eq:reg-t-pr-4}
\begin{split}
 (\bs m_{tt},\bs H_t) + ((\bs u_t\cdot\nabla)\bs m,\bs H_t) + ((\bs u\cdot\nabla)\bs m_t,\bs H_t) + \sigma(\div\bs m_t,\div\bs H_t) 
=  (\bs\omega_t\times \bs m,\bs H_t) + (\bs\omega\times \bs m_t,\bs H_t) \\-\frac{1}{\tau}(\bs m_t,\bs H_t) + \frac{\chi_0}{\tau}\|\bs H_t\|^2.
\end{split}
\end{equation}
Differentiating \eqref{eq:FHD-6} with $t$, we obtain
$$
\mu_0(\div\bs H_t + \div\bs m_t) = -\div \partial_t\bs H_e
$$
and
$$
\mu_0(\div\bs H_{tt} + \div\bs m_{tt}) = -\div \partial_{tt}\bs H_e.
$$
Equation \eqref{eq:FHD-5} implies there exists a unique $\varphi \in H_0^1(\Omega)$ such that $\bs H= \nabla\varphi$. Testing the above two equations with $\varphi_t$, we get
\begin{equation}\label{eq:reg-t-pr-5}
\mu_0(\|\bs H_t\|^2 + (\bs m_t,\bs H_t)) = -(\partial_t\bs H_e,\bs H_t)
\end{equation}
and
\begin{equation}
\label{eq:reg-t-pr-6}
\frac{\mu_0}{2}\frac{\dd }{\dd t}\|\bs H_t\|^2 + \mu_0(\bs m_{tt},\bs H_t) = -(\partial_{tt}\bs H_e,\bs H_t)
\end{equation}
Adding equations \eqref{eq:reg-t-pr-1} - \eqref{eq:reg-t-pr-4} together, and using the equations \eqref{eq:reg-t-pr-5} - \eqref{eq:reg-t-pr-6}, we obtain
\begin{equation}\label{eq:reg-t-pr-7}
\begin{split}
& \frac{1}{2}\frac{\dd}{\dd t}(\rho\|\bs u_t\|^2 + \rho\kappa\|\bs\omega_t\|^2 + \|\bs m_t\|^2 + \mu_0\|\bs H_t\|^2) +\eta\|\nabla\bs u_t\|^2 + \eta^\prime\|\nabla\bs\omega_t\|^2  
+ (\eta^\prime + \lambda^\prime)\|\div\bs\omega_t\|^2 +\sigma\|\bs m_t\|_n^2 \\
&\qquad \qquad + \mu_0\sigma\|\bs H_t\|_n^2 + \frac{1}{\tau}\|\bs m_t\|^2 
+ \frac{2\mu_0\chi_0 + \mu_0^2}{\tau}\|\bs H_t\|^2 + \zeta\|\curl\bs u_t-2\bs\omega_t\|^2 \\
= & \mu_0((\bs m_t\cdot\nabla)\bs H,\bs u_t) + \mu_0((\bs m\cdot\nabla)\bs H_t,\bs u_t) - \rho((\bs u_t\cdot\nabla)\bs u,\bs u_t) 
+ \mu_0(\bs m_t\times \bs H,\bs\omega_t) -\rho\kappa ((\bs u_t\cdot\nabla)\bs\omega,\bs\omega_t)  +(\bs\omega_t\times \bs m,\bs m_t) \\
&
- ((\bs u_t\cdot\nabla)\bs m,\bs m_t) + \mu_0((\bs u_t\cdot\nabla)\bs m,\bs H_t) + \mu_0((\bs u\cdot\nabla)\bs m_t,\bs H_t)  -\mu_0(\bs\omega\times \bs m_t,\bs H_t) -(\partial_{tt}\bs H_e,\bs H_t) - \sigma(\div\partial_t\bs H_e,\div\bs H_t),
\end{split}
\end{equation}
which, together with the inequalities
\begin{align*}
\mu_0|((\bs m_t\cdot\nabla)\bs H,\bs u_t)| & \leq \ \frac{\sigma}{10}\|\bs m_t\|_n^2 + C\|\bs H\|_2^2\|\bs u_t\|^2,\\
\mu_0|((\bs m\cdot\nabla)\bs H_t,\bs m_t)| & \leq \frac{\mu_0\sigma}{6}\|\bs H_t\|_n^2 +\frac{\sigma}{10}\|\bs m_t\|_n^2 + C\|\bs m\|_n^4\|\bs m_t\|^2,\\
\rho|((\bs u_t\cdot\nabla)\bs u,\bs u_t)| & \leq \frac{\eta}{8}\|\nabla\bs u_t\|^2 + C\|\bs u\|_2^2 \|\bs u_t\|^2,\\
\mu_0|(\bs m_t\times \bs H,\bs\omega_t)| & \leq \frac{\sigma}{10}\|\bs m_t\|_n^2 + C\|\bs H\|_n^2 \|\bs\omega_t\|^2,\\
\rho\kappa|((\bs u_t\cdot\nabla)\bs \omega,\bs\omega_t)| & \leq \frac{\eta}{8}\|\nabla\bs u_t\|^2 + C\|\bs\omega\|_2^2\|\bs \omega_t\|^2,\\
|(\bs\omega_t \times \bs m,\bs m_t)| & \leq \frac{\eta^\prime}{2}\|\nabla\bs\omega_t\|^2 + C\|\bs m\|_n^2 \|\bs m_t\|^2,\\
|((\bs u_t\cdot\nabla)\bs m,\bs m_t)| & \leq \frac{\eta}{8}\|\nabla\bs u_t\|^2 + C\|\bs m\|_2^2\|\bs m_t\|^2,\\
\mu_0|((\bs u_t\cdot\nabla)\bs m,\bs H_t)| & \leq \frac{\eta}{8}\|\nabla\bs u_t\|^2 + C\|\bs m\|_2^2 \|\bs H_t\|^2,\\
\mu_0|((\bs u\cdot\nabla)\bs m_t,\bs H_t)| & \leq \frac{\sigma}{10}\|\bs m_t\|_n^2 + \frac{\mu_0\sigma}{6}\|\bs H_t\|_n^2 + C\|\nabla\bs u\|^4\|\bs H_t\|^2,\\
\mu_0|(\bs\omega\times \bs m_t,\bs H_t)| & \leq \frac{\sigma}{10}\|\bs m_t\|_n^2 + C\|\nabla\bs\omega\|^2\|\bs H_t\|^2,\\
|(\partial_{tt}\bs H_e,\bs H_t)| & \leq \frac{1}{2}(\|\partial_{tt}\bs H_e\|^2 + \|\bs H_t\|^2), \\
\sigma |(\div\partial_t\bs H_e,\div\bs H_t) | &\leq \frac{\mu_0\sigma}{6}\|\bs H_t\|_n^2 + C \|\div\partial_t\bs H_e\|^2,
\end{align*}
further gives
\begin{equation}\label{eq:reg-t-pr-8}
\begin{split}
\frac{\dd}{\dd t}(\rho\|\bs u_t\|^2 + \rho\kappa\|\bs\omega_t\|^2 + \|\bs m_t\|^2 + \mu_0\|\bs H_t\|^2) +\eta\|\nabla\bs u_t\|^2 + \eta^\prime\|\nabla\bs\omega_t\|^2  
+ 2 (\eta^\prime + \lambda^\prime)\|\div\bs\omega_t\|^2 +\sigma\|\bs m_t\|_n^2 \\
+ \mu_0\sigma\|\bs H_t\|_n^2 + \frac{2}{\tau}\|\bs m_t\|^2 
+ \frac{4\mu_0\chi_0 + 2\mu_0^2}{\tau}\|\bs H_t\|^2 + 2\zeta\|\curl\bs u_t-2\bs\omega_t\|^2 \\
\leq C(\|\bs u\|_2^2 + \|\bs\omega\|_2^2 + \|\bs m\|_2^2 + \|\bs H\|_2^2 + \|\bs m\|_n^4 +\|\nabla\bs u\|^4+1)
 (\|\bs u_t\|^2 + \|\bs\omega_t\|^2 + \|\bs m_t\|^2 + \|\bs H_t\|^2)\\
+\frac{1}{2}\|\partial_{tt}\bs H_e\|^2 + C\|\div\partial_t\bs H_e\|^2 .
\end{split}
\end{equation}
For any $t_F\in (0,T]$, integrating \eqref{eq:reg-t-pr-8} with time variable $t$ from $0$ to $t_F$, we have
\begin{equation}\label{eq:reg-t-pr-9}
\begin{split}
& \rho\|\bs u_t(t_F)\|^2 + \rho\kappa\|\bs\omega_t(t_F)\|^2 + \|\bs m_t(t_F)\|^2 + \mu_0\|\bs H_t(t_F)\|^2 
 + \int_0^{t_F}(\|\nabla\bs u_t\|^2 + \|\nabla\bs\omega_t\|^2 + \|\bs m_t\|_n^2 + \|\bs H_t\|_n^2)\dd t \\
\leq & \rho\|\bs u_t(0)\|^2 + \rho\kappa\|\bs\omega_t(0)\|^2 + \|\bs m_t(0)\|^2 + \mu_0\|\bs H_t(0)\|^2  + C\int_0^{t_F} (\|\partial_{tt}\bs H_e\|^2 + \|\div\partial_t\bs H_e\|^2)\dd t \\
&\quad + C\int_0^{t_F} (\|\bs u\|_2^2 + \|\bs\omega\|_2^2 + \|\bs m\|_2^2 + \|\bs H\|_2^2 + \|\bs m\|_n^4 +\|\nabla\bs u\|^4+1) (\|\bs u_t\|^2 + \|\bs\omega_t\|^2 + \|\bs m_t\|^2 + \|\bs H_t\|^2)\dd t.\\
\end{split}
\end{equation}

In view of  the Gr\"onwall inequality, Theorem \ref{the:reg-2} and  Assumption \ref{ass:solu}, to get the desired conclusion it suffices to prove that 
 the term 
\begin{equation}\label{bound-term}
\rho\|\bs u_t(0)\|^2 + \rho\kappa\|\bs\omega_t(0)\|^2 + \|\bs m_t(0)\|^2 + \mu_0\|\bs H_t(0)\|^2
\end{equation}
 is bounded. 
In fact, since $\bs u \in \bs H_0^1(\Omega)\cap\bs H_0(\div 0)$, we have
\begin{align*}
& \|\bs u_t(0)\| = \sup\limits_{\bs v \in \bs H_0(\div 0)}\frac{(\bs u_t(0),\bs v)}{\|\bs v\|},
\end{align*}
which yields
\begin{align*}
\|\bs u_t(0)\| & \lesssim \|\bs u_0\|_2\|\bs u_0\|_1 + \|\bs u_0\|_2 + \|\bs m_0\|_2\|\bs H(0)\|_1.
\end{align*}
Equations \eqref{eq:FHD-3} and \eqref{eq:FHD-4} imply
\begin{align*}
\|\bs\omega_t(0)\| & \lesssim \|\bs u_0\|_1\|\bs \omega_0\|_2 + \|\bs\omega_0\|_2 + \|\bs m_0\|_1 \|\bs H_0\|_1 + \|\bs u_0\|_1 + \|\bs \omega_0\|, \\
\|\bs m_t(0)\| & \lesssim \|\bs u_0\|_1\|\bs m_0\|_2 + \|\bs m_0\|_2 + \|\bs \omega_0\|_1\|\bs m_0\|_1 +\|\bs m_0\| +\|\bs H(0)\|.
\end{align*}
Taking $t = 0$ in \eqref{eq:reg-t-pr-5}, we obtain
$$
\|\bs H_t(0)\|^2 \lesssim \|\bs m_t(0)\|\|\bs H_t(0)\| + \|\partial_t\bs H_e(0)\| \|\bs H_t(0)\|.
$$
As a result, the term \eqref{bound-term} is bounded from above. This finishes the proof.
\end{proof}

\begin{theorem}
\label{the:reg-2-inf} Under  Assumptions \ref{ass:He} - \ref{ass:omega}, 
the  solution   $( \bs u,\ p, ~\bs\omega,~\bs m, \bs H)$ of problem \eqref{eq:FHD-1} - \eqref{eq:initial} satisfies
$$
\sup\limits_{0\leq t \leq T}(\|\bs u\|_2 + \|\bs \omega\|_2 + \|\bs m\|_2 +\|\bs H\|_2 +\|p\|_1) \leq C.
$$
\end{theorem}
\begin{proof}
Equations \eqref{eq:FHD-1} - \eqref{eq:FHD-4} imply
\begin{align*}
-(\eta+\zeta)\Delta \bs u + \nabla p & = -\rho(\bs u_t + (\bs u\cdot\nabla)\bs u) + \mu_0(\bs m\cdot\nabla)\bs H +2\zeta\curl\bs\omega, \\
\div\bs u & = 0,\\
-\eta^\prime \Delta\bs\omega - (\eta^\prime + \lambda^\prime)\nabla\div\bs\omega & = -\rho\kappa(\bs\omega_t + (\bs u\cdot\nabla)\bs\omega) + \mu_0\bs m\times \bs H + 2\zeta(\curl\bs u - 2\bs\omega) ,\\
-\sigma\Delta\bs m & =-\bs m_t - (\bs u\cdot\nabla)\bs m + \bs\omega \times \bs m - \frac{1}{\tau}(\bs m - \chi_0\bs H). 
\end{align*}
 Assumption \ref{ass:omega} gives
\begin{align}
\label{eq:reg2-inf-1}& \|\bs u\|_2 + \|\bs p\|_1  \leq C\|\bs u_t\| + \frac{1}{2}\|\bs u\|_2 + C\|\bs u\|_1^3 + C\|\bs m\|_n\|\bs H\|_n^{1/2}\|\bs H\|_2^{1/2} + C\|\bs \omega\|_1,\\
\label{eq:reg2-inf-2}& \|\bs\omega\|_2  \leq  C\|\bs\omega_t\|  + \frac{1}{2}\|\bs\omega\|_2 + C\|\bs u\|_1^2\|\bs\omega\|_1 + C\|\bs m\|_n\|\bs H\|_n + C\|\bs u\|_1 + C\|\bs \omega\|,\\
\label{eq:reg2-inf-3}& \|\bs m\|_2  \leq C\|\bs m_t\| + \frac{1}{4}\|\bs m\|_2 + \|\bs u\|_1^2\|\bs m\|_n + C\|\bs\omega\|_1\|\bs m\|_n +C\|\bs m\| + C\|\bs H\|.
\end{align}
Equations \eqref{eq:FHD-5} and \eqref{eq:FHD-6} indicate
$$
-\Delta \bs H = \grad\div\bs m + \frac{1}{\mu_0} \grad\div\bs H_e,
$$
Assumption \ref{ass:omega} yields
\begin{equation}\label{eq:reg2-inf-4}
\|\bs H\|_2 \leq C(\|\bs m\|_2+\|\bs H_e\|_2).
\end{equation}
Adding inequalities \eqref{eq:reg2-inf-1} - \eqref{eq:reg2-inf-3} together and using the inequality \eqref{eq:reg2-inf-4}, we obtain
\begin{align*}
\|\bs u\|_2 + \|\bs \omega\|_2 + \|\bs m\|_2 +\|p\|_1  \leq & C(\|\bs u_t\| +\|\bs\omega_t\| + \|\bs m_t\| + \|\bs u\|_1^3 +\|\bs\omega\|_1+\|\bs u\|_1+\|\bs m\|) \\
& + C(\|\bs H\| + \|\bs m\|_n^2\|\bs H\|_n + \|\bs m\|_n\|\bs H\|_n^{1/2}\|\bs H_e\|_2^{1/2}) \\
& + C( \|\bs u\|_1^2\|\bs\omega\|_1 + \|\bs m\|_n\|\bs H\|_n + \|\bs u\|_1^2\|\bs m\|_n + \|\bs\omega\|_1\|\bs m\|_n),
\end{align*}
which plus Theorems \ref{the:reg-2} and \ref{the:reg-t} leads to
$$
\sup\limits_{0\leq t \leq T}(\|\bs u\|_2 + \|\bs \omega\|_2 + \|\bs m\|_2+\|p\|_1 ) \leq C.
$$
Using this estimate and \eqref{eq:reg2-inf-4}, we finally get the desired result.
\end{proof}

\begin{theorem}\label{the:reg-t-ener}
Under  Assumptions \ref{ass:He} and \ref{ass:solu}, \ref{ass:omega}  {and the conditions 
$$
\bs u_0\in (H^3(\Omega))^3,\quad\bs \omega_0 \in (H^3(\Omega))^3,\quad \bs m_0 \in (H^3(\Omega))^3,
$$} the solution   components $ \bs u,~\bs\omega,~\bs m,$ and $\bs H$ of problem \eqref{eq:FHD-1} - \eqref{eq:initial} satisfy
\begin{equation}
\label{eq:reg-t-ener}
\begin{split}
& \sup\limits_{0\leq t\leq T} (\|\nabla\bs u_t\|^2 + \|\nabla\bs\omega_t\|^2 + \|\bs m_t\|_n^2 + \|\bs H_t\|_n^2)
+\int_0^T (\|\bs u_{tt}\|^2 + \|\bs\omega_{tt}\|^2 + \|\bs m_{tt}\|^2 + \|\bs H_{tt}\|^2)\dd t\leq C.
\end{split}
\end{equation}
\end{theorem}
\begin{proof}
Testing \eqref{eq:FHD-t-1}, \eqref{eq:FHD-t-3} and \eqref{eq:FHD-t-4} with $\bs u_{tt}$, $\bs \omega_{tt}$ and $\bs m_{tt}$, respectively, and adding the resulting equations together, we get
\begin{equation}
\label{eq:reg-t-ener-pr-1}
\begin{split}
& \frac{1}{2}\frac{\dd}{\dd t}(\eta\|\nabla\bs u_t\|^2 + \eta^\prime\|\nabla\bs\omega_t\|^2 + (\eta^\prime+\lambda^\prime)\|\div\bs\omega_t\|^2 +\sigma\|\bs m_t\|_n^2 +\frac{1}{\tau}\|\bs m_t\|^2) 
+ \frac{\zeta}{2}\frac{\dd}{\dd t}\|\curl\bs u_t - 2\bs\omega_t\|^2 +\rho\|\bs u_{tt}\|^2 \\
&\qquad\qquad\qquad + \rho\kappa\|\bs \omega_{tt}\|^2 + \|\bs m_{tt}\|^2 \\
= &  \mu_0((\bs m_t\cdot\nabla)\bs H,\bs u_{tt}) + \mu_0((\bs m\cdot\nabla)\bs H_t,\bs u_{tt}) -\rho((\bs u_t\cdot\nabla)\bs u,\bs u_{tt}) - \rho((\bs u\cdot\nabla)\bs u_t,\bs u_{tt}) + \mu_0(\bs m_t\times \bs H,\bs\omega_{tt}) \\
&\qquad\qquad + \mu_0(\bs m\times\bs H_t,\bs\omega_{tt})  -\rho\kappa((\bs u_t\cdot\nabla)\bs\omega,\bs\omega_{tt}) - \rho\kappa((\bs u\cdot\nabla)\bs\omega_t,\bs\omega_{tt}) + (\bs\omega_t \times \bs m,\bs m_{tt})  + (\bs\omega\times \bs m_t,\bs m_{tt}) \\
&\qquad\qquad + \frac{\chi_0}{\tau}(\bs H_t,\bs m_{tt}) - ((\bs u_t\cdot\nabla)\bs m,\bs m_{tt})  - ((\bs u\cdot\nabla)\bs m_t,\bs m_{tt}).
\end{split}
\end{equation}
Testing \eqref{eq:FHD-t-4} with $\mu_0\bs H_{tt}$, we have
\begin{equation}\label{eq:reg-t-ener-pr-2}
\begin{split}
\mu_0(\bs m_{tt},\bs H_{tt}) + \mu_0((\bs u_t\cdot\nabla)\bs m,\bs H_{tt}) + \mu_0((\bs u\cdot\nabla)\bs m_t,\bs H_{tt}) 
+ \mu_0\sigma(\div\bs m_t,\div\bs H_{tt}) = \mu_0(\bs\omega_t\times \bs m,\bs H_{tt}) \\
+ \mu_0(\bs \omega\times\bs m_t,\bs H_{tt})
 -\frac{\mu_0}{\tau}(\bs m_t,\bs H_{tt}) + \frac{\mu_0\chi_0}{2\tau}\frac{\dd }{\dd t}\|\bs H_t\|^2.
\end{split}
\end{equation}
Differentiating \eqref{eq:FHD-6} with $t$, we obtain
\begin{equation}\label{eq:FHD-6-t}
\mu_0(\div\bs H_t + \div\bs m_t) = -\div \partial_t\bs H_e,
\end{equation}
and
\begin{equation}\label{eq:FHD-6-tt}
\mu_0(\div\bs H_{tt} + \div\bs m_{tt}) = -\div \partial_{tt}\bs H_e.
\end{equation}
Equation \eqref{eq:FHD-5} implies there exists a unique $\varphi \in H_0^1(\Omega)$, such that $\bs H= \nabla\varphi$. Testing \eqref{eq:FHD-6-t} with $\div\bs H_{tt}$ and $\varphi_{tt}$, and \eqref{eq:FHD-6-tt} with $\varphi_{tt}$, we arrive at
\begin{equation}\label{eq:reg-t-ener-pr-3}
\frac{\mu_0}{2}\frac{\dd}{\dd t}\|\bs H_t\|_n^2 + \mu_0(\div\bs m_t,\div\bs H_{tt}) = (\grad\div\partial_t\bs H_e,\bs H_{tt}),
\end{equation}
\begin{equation}
\label{eq:FHD-6-t-1}
\frac{\mu_0}{2}\frac{\dd}{\dd t}\|\bs H_t\|^2 + \mu_0(\bs m_t,\bs H_{tt})= -(\partial_t\bs H_e,\bs H_{tt}),
\end{equation}
and
\begin{equation}\label{eq:reg-t-ener-pr-4}
\mu_0\|\bs H_{tt}\|^2 + \mu_0(\bs m_{tt},\bs H_{tt}) =- (\partial_{tt}\bs H_e,\bs H_{tt}).
\end{equation}
Using \eqref{eq:reg-t-ener-pr-2}, \eqref{eq:reg-t-ener-pr-3}, \eqref{eq:FHD-6-t-1} and \eqref{eq:reg-t-ener-pr-4}, we get
\begin{equation}\label{eq:reg-t-ener-pr-5}
\begin{split}
\frac{1}{2}\frac{\dd}{\dd t}(\mu_0\sigma\|\bs H_t\|_n^2 + \frac{\mu_0(1+\chi_0)}{\tau}\|\bs H_t\|^2) + \mu_0\|\bs H_{tt}\|^2 = \mu_0((\bs u_t\cdot\nabla)\bs m,\bs H_{tt})
+ \mu_0((\bs u\cdot\nabla)\bs m_t,\bs H_{tt})  -\mu_0(\bs\omega_t\times \bs m,\bs H_{tt})  \\
- \mu_0(\bs\omega \times \bs m_t,\bs H_{tt}) 
-(\partial_{tt}\bs H_e,\bs H_{tt}) + \sigma(\grad\div\partial_t\bs H_e,\bs H_{tt}) -\frac{1}{\tau}(\partial_t\bs H_e,\bs H_{tt}).
\end{split}
\end{equation}
Adding \eqref{eq:reg-t-ener-pr-1} and \eqref{eq:reg-t-ener-pr-5} together, and using the inequalities
\begin{align*}
\mu_0|((\bs m_t\cdot\nabla)\bs H,\bs u_{tt})&+(\bs m\cdot\nabla)\bs H_t,\bs u_{tt})|  \leq \mu_0(\|\bs m_t\|_{L^6}\|\nabla\bs H\|_{L^3}  + \|\bs m\|_{L^\infty}\|\nabla\bs H_t\|)\|\bs u_{tt}\| \\
& \leq \frac{\rho}{4}\|\bs u_{tt}\|^2 + C\|\bs H\|_n\|\bs H\|_2\|\bs m_t\|_n^2    + C\|\bs m\|_2^2\|\bs H_t\|_n^2,\\
\rho|((\bs u_t\cdot\nabla)\bs u,\bs u_{tt}) & + ((\bs u\cdot\nabla)\bs u_t,\bs u_{tt})|  \leq\frac{\rho}{4}\|\bs u_{tt}\|^2 + C\|\bs u\|_2^2\|\nabla\bs u_t\|^2,\\
\mu_0|(\bs m_t\times\bs H,\bs \omega_{tt}) & + (\bs m\times \bs H_t,\bs\omega_{tt})| \leq \frac{\rho\kappa}{4}\|\bs\omega_{tt}\|^2 + C\|\bs H\|_n^2\|\bs m_t\|_n^2 + C \|\bs m\|_n^2 \|\bs H_t\|_n^2,\\
\rho\kappa|((\bs u_t\cdot\nabla)\bs\omega,\bs\omega_{tt}) & + ((\bs u\cdot\nabla)\bs\omega_t,\bs\omega_{tt}) \leq\frac{\rho\kappa}{4}\|\bs\omega_{tt}\|^2 + C\|\bs\omega\|_2^2\|\nabla\bs u_t\|^2 + C \|\bs u\|_2^2\|\nabla\bs\omega_t\|^2,\\
|(\bs\omega_t\times\bs m,\bs m_{tt}) & + (\bs \omega\times\bs m_t,\bs m_{tt}) \leq \frac{1}{6}\|\bs m_{tt}\|^2 + C\|\bs m\|_n^2\|\nabla\bs\omega_t\|^2 + C\|\nabla\bs\omega\|^2\|\bs m_t\|_n^2,\\
\frac{\chi_0}{\tau}& |(\bs H_t,\bs m_{tt})| \leq \frac{1}{6}\|\bs m_{tt}\|^2 + C\|\bs H_t\|_n^2, \\
|((\bs u_t\cdot\nabla)\bs m,\bs m_{tt}) & + ((\bs u\cdot\nabla)\bs m_t,\bs m_{tt})| \leq \frac{1}{6}\|\bs m_{tt}\|^2 + C\|\bs m\|_2^2\|\nabla\bs u_t\|^2 + C \|\bs u\|_2^2\|\bs m_t\|_n^2,\\
\mu_0|((\bs u_t\cdot\nabla)\bs m,\bs H_{tt}) & + ((\bs u\cdot\nabla)\bs m_t,\bs H_{tt})| \leq \frac{\mu_0}{6}\|\bs H_{tt}\|^2 + C\|\bs m\|_2^2\|\nabla\bs u_t\|^2 + C\|\bs u\|_2^2\|\bs m_t\|_n^2,\\
\mu_0|(\bs\omega_t\times\bs m,\bs H_{tt}) & + (\bs\omega\times \bs m_t,\bs H_{tt})| \leq \frac{\mu_0}{6}\|\bs H_{tt}\|^2 + C\|\bs m\|_2^2\|\nabla\bs\omega_t\|^2 + C\|\bs\omega\|_2^2\|\bs m_t\|_n^2,\\
|(\partial_{tt}\bs H_e,\bs H_{tt}) & - \sigma(\grad\div\partial_t\bs H_e,\bs H_{tt}) + \frac{1}{\tau}(\partial_t\bs H_e,\bs H_{tt})| \leq \frac{\mu_0}{6}\|\bs H_{tt}\|^2 + C(\|\partial_{tt}\bs H_e\|^2 + \|\partial_t\bs H_e\|_2^2),
\end{align*}
we have
\begin{align*}
& \frac{\dd}{\dd t}(\eta\|\nabla\bs u_t\|^2 + \eta^\prime\|\nabla\bs\omega_t\|^2 + (\eta^\prime+\lambda^\prime)\|\div\bs\omega_t\|^2 +\sigma\|\bs m_t\|_n^2 + \mu_0\sigma\|\bs H_t\|_n^2)
\\
&\quad
+ \frac{\dd }{\dd t}\left ( \frac{1}{\tau}\|\bs m_t\|^2 + \zeta\|\curl\bs u_t - 2\bs\omega_t\|^2 + \frac{\mu_0(1+\chi_0)}{\tau}\|\bs H_t\|^2 \right)
+\rho\|\bs u_{tt}\|^2 + \rho\kappa\|\bs\omega_{tt}\|^2 + \|\bs m_{tt}\|^2 + \mu_0\|\bs H_{tt}\|^2 \\
\leq & C(\|\partial_{tt}\bs H_e\|^2 + \|\partial_t\bs H_e\|_2^2)  + C(\|\bs u\|_2^2 + \|\bs\omega\|_2^2 + \|\bs m\|_2^2 + \|\bs H\|_2^2 +1)(\|\nabla\bs u_t\|^2 + \|\nabla\bs\omega_t\|^2 + \|\bs m_t\|_n^2 + \|\bs H_t\|_n^2). 
\end{align*}
Integrating the above inequality with respect to  $t$ from $0$ to any $t_F \in (0,T]$, 
we have
\begin{align*}
& (\eta\|\nabla\bs u_t(t_F)\|^2 + \eta^\prime\|\nabla\bs\omega_t(t_F)\|^2  +\sigma\|\bs m_t(t_F)\|_n^2 + \mu_0\sigma\|\bs H_t(t_F)\|_n^2)
\\
&\qquad +\int_0^{t_F}(\rho\|\bs u_{tt}\|^2 + \rho\kappa\|\bs\omega_{tt}\|^2 + \|\bs m_{tt}\|^2 + \mu_0\|\bs H_{tt}\|^2)\dd t \\
\leq & C(\|\nabla\bs u_t(0)\|^2 + \|\nabla\bs\omega_t(0)\|^2 + \|\bs m_t(0)\|_n^2 +\|\bs H_t(0)\|_n^2) + C\int_0^T(\|\partial_{tt}\bs H_e\|^2 + \|\partial_t\bs H_e\|_2^2) \dd t \\
& + C\int_0^{t_F}(\|\bs u\|_2^2 + \|\bs\omega\|_2^2 + \|\bs m\|_2^2 + \|\bs H\|_2^2 +1)(\|\nabla\bs u_t\|^2 + \|\nabla\bs\omega_t\|^2 + \|\bs m_t\|_n^2 + \|\bs H_t\|_n^2) \dd t.
\end{align*}
{Testing \eqref{eq:FHD-1}, \eqref{eq:FHD-3} and \eqref{eq:FHD-4} at time $t = 0$ with $\mathcal A_1(\bs u_t(0))$, $-\Delta\bs\omega_t(0)$ and $\mathcal A_2(\bs m_t(0))$, respectively, we obtain
\begin{align*}
\|\nabla\bs u_t(0)\|^2 & \lesssim (\|\bs u_0\|_3 + \|\bs\omega_0\|_2 + \|\bs m_0\|_2\|\bs H_0\|_2)\|\nabla\bs u_t(0)\|, \\
\|\nabla\bs\omega_t(0)\|^2 & \lesssim (\|\bs\omega_0\|_3 + \|\bs u_0\|_2\|\bs\omega_0\|_2 + \|\bs H_0\|_2\|\bs m_0\|_2)\|\nabla\bs\omega_t(0)\|,\\
\|\bs m_t(0)\|_n^2 & \lesssim (\|\bs m_0\|_3 + \|\bs u_0\|_2\|\bs m_0\|_2 + \|\bs\omega_0\|_2 \|\bs m_0\|_2 + \|\bs H_0\|_1)\|\bs m_t\|_n.
\end{align*}}
%
%
Equation \eqref{eq:FHD-6-t} implies
$$
\|\nabla\bs u_t(0)\| +\|\nabla\bs\omega_t(0)\| + \|\bs m_t(0)\|_n + \|\bs H_t(0)\|_n \leq C.
$$
Finally, applying the Gr\"onwall inequality and using Theorems \ref{the:reg-2} - \ref{the:reg-2-inf}, we obtain the desired result.
\end{proof}


\begin{theorem}\label{eq:reg-t-2-norm}
Under the same conditions as in Theorem \ref{the:reg-t-ener}, the solution   $( \bs u,\ p, ~\bs\omega,~\bs m, \bs H)$ of problem \eqref{eq:FHD-1} - \eqref{eq:initial}  satisfies
$$
\int_0^T(\|\bs u_t\|_2^2+\|p_t\|_1^2 + \|\bs\omega_t\|_2^2 + \|\bs m_t\|^2_2 +\|\bs H_t\|_2^2 )\dd t \leq C.
$$
\end{theorem}
\begin{proof}
Equations \eqref{eq:FHD-t-1} - \eqref{eq:FHD-t-4} and   Assumption \ref{ass:omega} imply
\begin{align*}
\|\bs u_t\|_2 + \|\bs p_t\|_1 & \lesssim \|\bs u_{tt}\| + \|\nabla\bs u_t\|\|\bs u\|_2 + \|\bs m_t\|_n\|\bs H\|_2 + \|\bs m\|_2\|\bs H_t\|_n + \|\bs\omega_t\|_1,\\
\|\bs\omega_t\|_2 & \lesssim  \|\bs\omega_{tt}\| + \|\bs u_t\|_1\|\bs \omega\|_2 + \|\bs u\|_2\|\bs\omega_t\|_1 + \|\bs H\|_n\|\bs m_t\|_n + \|\bs m\|_n\|\bs H_t\|_n  +\|\bs u_t\|_1 + \|\bs \omega_t\|, \\
\|\bs m_t\|_2 &\lesssim \|\bs m_{tt}\| + \|\bs u_t\|_1\|\bs m\|_2 + \|\bs u\|_2\|\bs m_t\|_n + \|\bs\omega_t\|_1\|\bs m\|_n  + \|\bs \omega\|_1\|\bs m_t\|_n + \|\bs m_t\| + \|\bs H_t\|,
\end{align*}
which, together with Theorems \ref{the:reg-2} - \ref{the:reg-t-ener}, give
$$
\int_0^T(\|\bs u_t\|_2^2+\|p_t\|_1^2 + \|\bs\omega_t\|_2^2 + \|\bs m_t\|_2^2 )\dd t \leq C.
$$
Equation \eqref{eq:FHD-6-t} and \eqref{eq:FHD-5} indicate
$$
-\Delta\bs H_{t} = \grad\div\bs m_t + \frac{1}{\mu_0}\grad\div\partial_t\bs H_e,
$$
which plus Assumption \ref{ass:omega} yields
$$
\|\bs H_t\|_2^2 \lesssim \|\bs m_t\|_2^2 + \|\partial_t\bs H_e\|_2^2.
$$
This leads to the desired result.
\end{proof}

By Theorems \ref{the:reg-2} - \ref{eq:reg-t-2-norm} and the definitions of   $\bs z$, $\bs k$ and $\tilde p$, we have the following a priori estimates of these auxiliary variables.
\begin{theorem}\label{the:reg-aux}
Under the same conditions as in Theorem \ref{the:reg-t-ener}, the variables $\bs z$, $\bs k$ and $\tilde p$ defined in \eqref{eq:kcurlm} satisfy
$$
\int_0^T (\|\bs z\|_1^2 +\|\bs k\|_1^2 + \|\tilde p\|_1^2) \dd t \leq C
$$
and
$$
\sup\limits_{0 \leq t \leq T} (\|\bs k\|_1 + \|\tilde p\|_1) \leq C.
$$
\end{theorem}
\begin{proof}
The definition of $\bs z$ implies
\begin{align*}
\int_0^T\|\bs z\|_1^2\dd t & \leq C\int_0^T \left( \int_\Omega (|\nabla\bs u|^2|\bs m|^2 +  |\bs u|^2 |\nabla\bs m|^2)\dd\Omega\right)\dd t  \leq C\sup\limits_{0\leq t \leq T}(\|\bs m\|_2^2 + \|\bs u\|_2^2)\int_0^T(\|\nabla\bs u\|^2 + \|\bs m\|_n^2) \dd t,
\end{align*}
which, together with   Theorems \ref{the:ener} and \ref{the:reg-2-inf}, yields the desired result for $\bs z$.

The definition of $\bs k$ and Theorems \ref{the:reg-2} and \ref{the:reg-2-inf} indicate
$$
\int_0^T\|\bs k\|_1^2\dd t \leq \int_0^T\|\bs m\|_2^2 \dd t  \leq C
$$
and 
$$
\sup\limits_{0\leq t \leq T}\|\bs k\|_1 \leq \sup\limits_{0\leq t \leq T}\|\bs m\|_2 \leq C.
$$

The definition of $\tilde p$ shows
\begin{align*}
\sup\limits_{0\leq t\leq T}\|\tilde p\|_1  \leq C\sup\limits_{0\leq t\leq T}(\|p\|_1 + \|\bs m\|_2\|\bs H\|_n + \|\bs m\|_n\|\bs H\|_2)
\end{align*}
and 
\begin{align*}
\int_0^T\|\tilde p_t\|_1^2\dd t & \leq C \int_0^T\big(\|p_t\|_1^2 + \|\bs m_t\|_2^2\|\bs H\|_n^2 + \|\bs m_t\|_n^2\|\bs H\|_2^2  + \|\bs m\|_2^2\|\bs H_t\|_n^2 + \|\bs m\|_n^2\|\bs H_t\|_2^2\big)\dd t \\
& \leq C \sup\limits_{0\leq t \leq T}(1+ \|\bs H\|_n^2 + \|\bs H\|_2^2 + \|\bs m\|_2^2 + \|\bs m\|_n^2) 
\int_0^T (\|p_t\|_1^2 + \|\bs m_t\|_2^2 + \|\bs m_t\|_n^2 +\|\bs H_t\|_n^2 + \|\bs H_t\|_2^2)\dd t,
\end{align*}
which, together with Theorems  \ref{the:reg-2}, \ref{the:reg-t}, \ref{the:reg-2-inf} and \ref{eq:reg-t-2-norm}, imply the desired results for $\tilde p$.
\end{proof}

\section{Semi-discrete finite element scheme}\label{sec:semi}
In this section, we will first introduce   finite element spaces and some preliminary results, 
then give  the semi-discrete finite element scheme   of the FHD model \eqref{eq:FHD-1} - \eqref{eq:initial} and  carry out   analyses of  stability,  solvability and  convergence.

\subsection{Finite element spaces}
We consider some $H^1$-, $H(\curl)$- and $H(\div)$- conforming finite element spaces which will be used in the spatial discretization of the weak problem \eqref{eq:FHD-weak-1} - \eqref{eq:FHD-weak-8}.


We introduce the following finite element spaces:
\begin{itemize}
\item $\bs S_h = (S_h)^3 \subset  \bs S$, where $S_h$ is a space of continuous piecewise  polynomials ~(cf. \citet{Ciarlet1978}) for the velocity $\bs u$, with $S_h|_K = \mathcal P_1(K) + B_K$ for any $K \in \mathcal T_h$. Here $B_K$ is the bubble function space  
 defined as
\begin{align*}
B_K = \text{span}\{\lambda_1^K\lambda_2^K\lambda_3^K\lambda_4^K:
\lambda_i^K~~(i = 1,\dots,4)&~\text{is the barycentric coordinates of }K\}.
\end{align*}
\item $L_h \subset H^1(\Omega) \cap W$ is the continuous linear element space~(cf. \citet{Ciarlet1978}) for the pressure variable $p$, with $L_h|_K = \mathcal P_1(K)$ for any $K \in \mathcal T_h$.
\item $\bs \Sigma_h = (\Sigma_h)^3 \subset \bs S$, where $\Sigma_h$ is the continuous linear element space   for the angular momentum variable $\bs\omega$, with $\Sigma_h|_K = \mathcal P_1(K)$ for any $K \in \mathcal T_h$.
\item $\bs U_h \subset \bs U$ is the lowest order edge-element space ($NE_0$, cf.  \citet{Nedelec1980,Nedelec1986}) for the new variables $\bs z$ and $\bs k$, with $\mathcal P_0(K)^3 \subset \bs U_h|_K$ for any $K \in \mathcal T_h$.
\item $\bs V_h \subset V$ is the lowest order face-element space ($RT_0$, cf.    \citet{Raviart;Thomas1977,Nedelec1980,Nedelec1986,Brezzi;Douglas;Duran;Fortin1987,Brezzi;Douglas;Marini1985}) for the interior magnetization $\bs m$ and the demagnetizing field $\bs H$, with $\mathcal P_0(K)^3 \subset \bs V_h|_K$.
\item $W_h\subset W$ is the piecewise constant space 
for the new variable $\phi$, with $W_h|_K = \mathcal P_0(K)$ for any $K \in \mathcal T_h$.
\end{itemize}
The above choice of finite element spaces satisfy the   inf-sup condition~(cf. \citet{Arnoldmini1984}):
\begin{equation}\label{eq:inf-sup}
\sup\limits_{\bs v_h \in \bs S_h}\frac{(q_h,\div\bs v_h)}{\|\bs v_h\|_1} \gtrsim \|q_h\|\qquad\forall~q_h \in L_h
\end{equation}
 and the commutative exact sequence~(cf. \citet{Hiptmair2002}):
 \begin{equation}\label{eq:exact_sq}
\begin{CD}
H_0^1@>{\grad}>> \bs U @>{\curl}>> \bs V @>{\div}>> W \\
@VV \pi_h V @ VV \pi ^{c}_h V  @VV \pi ^d_h V @ VV Q_h V \\\
\Sigma_{h} @>{\grad}>> \bs U_{h} @>{\curl}>>  \bs V_{h} @>{\div}>> W_{h}
\end{CD}
\end{equation}
Here $\pi_{h}:H_{0}^{1}(\Omega)\rightarrow \Sigma_{h}$, $\pi_{h}^{c}: \bs U \rightarrow \bs U_{h}$ and $\pi_{h}^{d}:\bs V\rightarrow \bs V_{h}$ are the classical interpolation operators, and $Q_{h}:W\rightarrow W_{h}$ is the $L^{2}$ orthogonal projection operator. 
Note that the diagram \eqref{eq:exact_sq} also indicates that 
$$ \grad \Sigma_{h} \subset \bs U_{h}, \quad  \curl \bs U_{h} \subset \bs V_{h}, \quad \div  \bs V_{h} = W_{h}.$$

We define   three discrete weak operators, $$\div_{h}:\bs U_{h} \rightarrow \Sigma_{h},\quad  \curl_{h}:\bs V_{h} \rightarrow \bs U_{h}, \quad \grad_{h}:\ W_{h}\rightarrow \bs V_{h},$$ as the adjoint operators of $-\grad$, $\curl$ and $-\div$, respectively, i.e.,  for  any $\bs u_{h}\in \bs U_{h} $, 
  $\div_{h}\bs u_{h}\in \Sigma_{h}$ satisfies 
\begin{equation}\label{eq:weak-div}
(\div_{h}\bs u_{h},s_{h}) := -(\bs u_{h},\grad s_{h})\qquad\forall~~s_{h} \in \Sigma_{h};
\end{equation}
for any $\bs v_{h}\in \bs V_{h} $,  $\curl_{h} \bs v_{h} \in \bs U_{h}$ satisfies
\begin{equation}\label{eq:weak-curl}
(\curl_{h}\bs v_{h},\bs u_{h}) := (\bs v_{h},\curl\bs u_{h})\qquad\forall~~\bs u_{h} \in \bs U_{h}; 
\end{equation}
and for any $w_{h}\in W_{h} $,
$\grad_{h} w_{h} \in \bs V_{h}$ satisfies 
\begin{equation}
\label{eq:weak-grad}
(\grad_{h} w_{h},\bs v_{h}) := -(w_{h},\div\bs v_{h})\qquad\forall~~ v_{h} \in \bs V_{h}.
\end{equation}
Thus, we have the following reversed ordering exact sequence:
\begin{equation}\label{eq:exact-seq-re}
\begin{CD}
0@<{}<<\Sigma_h @<{\div_h}<< \bs U_h @<{\curl_h}<< \bs V_h @<{\grad_h}<< W_h@<{}<<0
\end{CD}
\end{equation} 

Introduce the null spaces of the differential operators
$$
\mathfrak Z_{h}^{c}: = \bs U_{h} \cap \ker(\curl)\quad\text{and}\quad \mathfrak Z_{h}^{d}: = \bs V_{h}\cap \ker(\div),
$$
and the null spaces of the weak differential operators
\begin{equation}\label{Khd}
\mathfrak K_{h}^{c} := \bs U_{h}\cap \ker(\div_{h})\quad\text{and}\quad \mathfrak K_{h}^{d} := \bs V_{h} \cap \ker(\curl_{h}).
\end{equation}
Similarly, we use the notations $\mathfrak Z^{c}$ and $\mathfrak Z^{d}$ to denote the null spaces in the continuous level, i.e.
$$
\mathfrak Z^{c} := \bs U \cap \ker(\curl),\quad  \mathfrak Z^{d}: = \bs V\cap \ker(\div),$$
 and  also define 
$$
\mathfrak K^{c} := \bs U \cap (\mathfrak Z^{c})^{\bot} \quad\text{and}\quad \mathfrak K^{d} := \bs V \cap (\mathfrak Z^{d})^{\bot}.
$$
We have the following Hodge decompositions (cf. \citet{Arnold;Falk;Winther2006,Arnold;Falk;Winther2010}):
\begin{align}
\label{eq:hodge-dis-1}
\bs U_{h} & = \mathfrak Z_{h}^{c} \oplus ^{\bot}\mathfrak K_{h}^{c}= \grad\Sigma_{h} \oplus^{\bot} \div_{h}\bs U_{h},\\ 
\label{eq:hodge-dis-2}
\bs V_{h} &= \mathfrak Z_{h}^{d} \oplus^{\bot}\mathfrak K_{h}^{d} = \curl\bs U_{h} \oplus^{\bot} \grad_{h}W_{h},
\end{align}
where $\oplus^\bot$ stand for the $L^2$ orthogonal decomposition.
There hold the following discrete Poincar\'e inequalities (cf. \citet{Arnold;Falk;Winther2006,Arnold;Falk;Winther2010,Chen2016MultiGrid}):
\begin{equation}\label{eq:P-1}
\|\curl \bs u_{h}\| \gtrsim  \|\bs u_{h}\| \qquad \forall~~\bs u_{h} \in \mathfrak K_{h}^{c},
\end{equation}
\begin{equation}\label{eq:P-2}
\|\div\bs v_{h}\| \gtrsim  \|\bs v_{h}\|\qquad\forall~~\bs v_{h} \in \mathfrak K_{h}^{d}.
\end{equation}
\begin{equation}\label{eq:P-3}
\|\grad_{h} w_{h} \|\gtrsim  \|w_{h}\|\qquad\forall~~w_{h} \in W_{h},
\end{equation}
\begin{equation}\label{eq:P-4}
\|\curl_{h} \bs v_{h}\| \gtrsim \|\bs v_{h}\|\qquad\forall~~\bs v_{h} \in \mathfrak Z_{h}^{d}.
\end{equation}

For any  $s\in H_{0}^{1}(\Omega)$, $\bs u \in \bs U$ and $\bs v \in \bs V$, we define $P_{h}^{g} s \in \Sigma_{h}$, $P_{h}^{c} \bs u \in \mathfrak K_{h}^{c}$ and $P_{h}^{d} \bs v \in \mathfrak K_{h}^{d}$ respectively by 
\begin{equation}\label{eq:pj}
(\grad P_{h}^{g}s,\grad \sigma_{h}) = (\grad s,\grad\sigma_{h})\qquad \forall~~\sigma_{h} \in \Sigma_{h},
\end{equation}
\begin{equation}\label{eq:project-1}
(\curl P_{h}^{c}\bs u,\curl \bs\phi_{h}) = (\curl \bs u,\curl\bs \phi_{h})\qquad\forall~~\bs\phi_{h} \in \mathfrak K_{h}^{c}
\end{equation}
and
\begin{equation}\label{eq:project-2}
(\div P_{h}^{d} \bs v,\div \bs\psi_{h} ) = (\div\bs v,\div\bs\psi_{h})\qquad\forall~~\bs\psi_{h} \in \mathfrak K_{h}^{d}.
\end{equation}

The Hodge decomposition on the continuous level implies that for any $\bs u \in \bs U$ and $\bs v \in \bs V$, there exist $u_{1} \in H_{0}^{1}(\Omega)$, $\bs u_{2} \in \mathfrak K^{c}$, $\bs v_{1} \in \mathfrak K^{c}$ and $\bs v_{2} \in \mathfrak K^{d}$ such that
$$
\bs u = \grad u_{1} \oplus^{\bot} \bs u_{2}\quad \text{and}\quad\bs v = \curl\bs v_{1} \oplus^{\bot} \bs v_{2}.
$$
We introduce two projection-based quasi-interpolation operators,  $I_{h}^{c}:\bs U \rightarrow \bs U_{h}$ and $I_{h}^{d}:\bs V \rightarrow \bs V_{h}$,   defined by
\begin{equation}\label{eq:inter-1}
I_{h}^{c} \bs u = \grad P_{h}^{g}\bs u_{1} \oplus^{\bot} P_{h}^{c}\bs u_{2}
\end{equation}
and
\begin{equation}
\label{eq:inter-2}
I_{h}^{d}\bs v = \curl P_{h}^{c}\bs v_{1} \oplus^{\bot} P_{h}^{d}\bs v_{2}.
\end{equation}
 The quasi-interpolation operators have the following properties; cf. \citep[Lemmas 3.3-3.5]{Wu2020} and \citep[Lemma 2.10]{WuXie2023}.
\begin{lemma} \label{lem:Ih}
The projection-based quasi-interpolation operators $I_h^c$ and $I_h^d$ have the following properties:
\begin{enumerate}
\item For any $\bs u \in \bs U$ and $\bs v \in \bs V$, there hold
\begin{align*}
& (I_{h}^{c}\bs u,\grad s_{h}) = (\bs u,\grad s_{h})\qquad\qquad\qquad\forall~~s_{h} \in \Sigma_{h},\\
& (I_{h}^{d} \bs v,\curl\bs\phi_{h}) = (\bs v,\curl\bs\phi_{h})\,\qquad\qquad\qquad\forall~~\bs\phi_{h} \in \bs U_{h},\\
& (\curl I_{h}^{c}\bs u,\curl\bs\psi_{h}) = (\curl\bs u,\curl\bs\psi_{h})\,\,\,\qquad\forall~~\bs\psi_{h} \in \bs U_{h},\\
& (\div I_{h}^{d}\bs v,\div\bs\phi_{h}) = (\div\bs v,\div\bs\phi_{h})\qquad\qquad\forall~~\bs\phi_{h} \in \bs V_{h}.
\end{align*}

\item For any $\bs u \in \bs U$ and $\bs v\in \bs V$, there hold
$$
\|\curl  I_{h}^{c}\bs u\| \leq \|\curl\bs u\|,\qquad
\|\div I_{h}^{d}\bs v\| \leq \|\div\bs v\|.
$$

\item For any $\bs u \in  H_{0}(\curl) \cap  H(\div)$ and $\bs v \in  H_{0}(\div) \cap  H(\curl)$, there hold
  $$
  \div_{h}  I_{h}^{c} \bs u = Q_{h}^{g} \div\bs u, \quad \curl_{h}  I_{h}^{d} \bs v = Q_{h}^{c} \curl \bs v,
  $$
  where $Q_{h}^{g}:L^{2}(\Omega) \rightarrow \Sigma_{h}$ and $Q_{h}^{c}: \bs L^{2}(\Omega) \rightarrow \bs U_{h}$ are $L^{2}$ orthogonal projection operators. Therefore,
$$
\|\div_{h} I_{h}^{c} \bs u \|\leq \|\div\bs u\|, \quad\|\curl_{h} I_{h}^{d}\bs v\| \leq \|\curl\bs v\|.
$$ 

\item For any $\bs u \in \bs U\cap  H^{l+1}(\Omega)$ and $\bs v \in \bs V \cap  H^{l+1}(\Omega)$ with $l \geq 0$, there hold
\begin{align*}
\|\bs u -  I_{h}^{c}\bs u\|  \lesssim h^{r}\|\bs u\|_{r}\qquad\text{for}\quad 1\leq r \leq l+1,\\
\|\bs v -  I_{h}^{d}\bs v\|  \lesssim h^{r}\|\bs v\|_{r}\qquad\text{for}\quad 1\leq r \leq l+1.
\end{align*}
Furthermore, if $\curl\bs u\in H^{l+1}(\Omega)$ and $\div\bs v \in H^{l+1}(\Omega)$, then there hold
\begin{align*}
\|\curl(I -  I_{h}^{c})\bs u\| \lesssim h^{r} \|\curl\bs v\|_{r}\qquad\text{for}\quad 1\leq r \leq l+1,\\
\|\div(I -  I_{h}^{d})\bs v\| \lesssim h^{r} \|\div\bs v\|_{r}\qquad\text{for}\quad 1\leq r \leq l+1.
\end{align*}
\item For any $\bs u \in L^\infty(\Omega) \cap \bs U$ and $\bs v \in L^\infty(\Omega) \cap \bs V$, there hold
$$
\|I_h^c\bs u\|_{L^\infty} \lesssim \|\bs u\|_{L^\infty}
$$
and
$$
\|I_h^d\bs v\|_{L^\infty} \lesssim \|\bs v\|_{L^\infty}.
$$
\end{enumerate}
\end{lemma}
For any $\bs u_h \in \bs U_h$ and $\bs v_h\in \bs V_h$, we also use the notation $\|\cdot\|_n$ to denote the norm with 
\begin{align*}
\|\bs u_h\|_n = (\|\curl\bs u_h\|^2 + \|\div_h\bs u_h\|^2)^{1/2} \quad\text{and}\quad
\|\bs v_h\|_n = (\|\curl_h\bs v_h\|^2 + \|\div\bs v_h\|^2)^{1/2}.
\end{align*}
The Poincar\'e inequalities \eqref{eq:P-1} - \eqref{eq:P-4} imply
$$
\|\bs u_h\| \lesssim\|\bs u_h\|_n,\quad\text{and}\quad\|\bs v_h\| \lesssim \|\bs v_h\|_n.
$$
We also have the following inequalities (see e.g. ~\citet{Girault1986,WuXie2023}):
\begin{align}
\label{eq:H-01}  & \|\bs v\|_{L^3} \leq \|\bs v\|^{1/2}\|\bs v\|_{L^6}^{1/2} \lesssim \|\bs v\|^{1/2}\|\nabla\bs v\|^{1/2} &\forall~\bs v \in \bs S \text{ or } \bs S_h, \\
\label{eqH01-2} & \|\bs v\|_{L^4} \leq \|\bs v\|^{1/2}\|\bs v\|_{L^6}^{3/4} \lesssim \|\bs v\|^{1/2}\|\nabla\bs v\|^{3/4}  & \forall~\bs v \in \bs S \text{ or } \bs S_h, \\
\label{eq:UV} & \|\bs v\|_{L^3} \leq \|\bs v\|^{1/2}\|\bs v\|_{L^6}^{1/2} \lesssim \|\bs v\|^{1/2}\|\bs v\|_n^{1/2} & \forall~\bs v \in \bs U\cap\bs V,~\bs U_h\text{ or }\bs V_h,\\
 \label{eq:UV-2} & \|\bs v\|_{L^4} \leq \|\bs v\|^{1/2}\|\bs v\|_{L^6}^{3/4} \lesssim \|\bs v\|^{1/2}\|\bs v\|_n^{3/4} & \forall~\bs v \in \bs U\cap\bs V,~\bs U_h\text{ or }\bs V_h.
\end{align}

%

\subsection{Semi-discrete scheme}
The semi-discrete scheme of the FHD model based on the weak form \eqref{eq:FHD-weak-1} - \eqref{eq:FHD-weak-8} reads:
Find $\bs u_h:[0,T]\rightarrow\bs S_h$, $\tilde p_h:[0,T]\rightarrow L_h$, $\bs\omega_h:[0,T]\rightarrow\bs \Sigma_h$, $\bs m_h:[0,T]\rightarrow \bs V_h$, $\bs z_h:[0,T]\rightarrow\bs U_h$, $\bs k_h:[0,T]\rightarrow \bs U_h$, $\bs H_h:[0,T]\rightarrow\bs V_h$ and $\phi_h:[0,T]\rightarrow W_h$ such that
\begin{align}
\label{eq:FHD-weak-dis-1}
\begin{split}
 \rho(\partial_t\bs u_h,\bs v_h) + \rho b(\bs u_h,\bs u_h,\bs v_h) + \eta(\nabla\bs u_h,\nabla\bs v_h) + \zeta (\curl\bs u_h,\curl\bs v_h) - \mu_0 c(\bs v_h,\bs m_h,\bs H_h)
  & \\ 
  - (\tilde p_h,\div\bs v_h)
  - \frac{\mu_0}{2}(\bs v_h \times \bs k_h,\bs H_h) 
  - \frac{\mu_0}{2}(\bs m_h\times \curl\bs v_h,\bs H_h) - 2\zeta(\bs\omega_h,\curl\bs v_h) &= 0
\end{split} & \forall ~\bs v_h \in \bs S_h,\\
\label{eq:FHD-weak-dis-2}
 (\div\bs u_h,q_h) & = 0 &\forall~q_h\in L_h, \\
\label{eq:FHD-weak-dis-3}
\begin{split}
\rho\kappa(\partial_t\bs\omega_h,\bs s_h) + \eta^\prime(\nabla\bs\omega_h,\nabla\bs s_h) + (\eta^\prime + \lambda^\prime)(\div\bs\omega_h,\div\bs s_h) 
+ \rho\kappa b(\bs u_h,\bs\omega_h,\bs s_h) &\\
-\mu_0(\bs m_h\times \bs H_h,\bs s_h) - 2\zeta(\curl\bs u_h,\bs s_h) 
+ 4\zeta(\bs \omega_h,\bs s_h)& = 0
\end{split} & \forall~\bs s_h \in \bs \Sigma_h,\\
\label{eq:FHD-weak-dis-4}
\begin{split}
(\partial_t\bs m_h,\bs F_h) - c(\bs u_h,\bs m_h,\bs F_h) + \sigma(\div\bs m_h,\div\bs F_h) + \sigma(\curl\bs k_h,\bs F_h)- \frac{1}{2}(\curl\bs z_h,\bs F_h)&  \\
-\frac{1}{2}(\bs m_h\times\curl\bs u_h,\bs F_h) 
- \frac{1}{2}(\bs u_h \times \bs k_h,\bs F_h) 
 - (\bs\omega_h\times \bs m_h,\bs F_h) &\\
 + \frac{1}{\tau}(\bs m_h,\bs F_h) 
-\frac{\chi_0}{\tau}(\bs H_h,\bs F_h) & = 0
\end{split}& \forall~\bs F_h\in \bs V_h,\\
\label{eq:FHD-weak-dis-5} (\bs z_h,\bs\Lambda_h) - (\bs u_h\times \bs m_h,\bs\Lambda_h) & = 0 & \forall~\bs \Lambda_h \in \bs U_h, \\
\label{eq:FHD-weak-dis-6} (\bs k_h,\bs \Theta_h) - (\bs m_h,\curl\bs\Theta_h) & = 0 & \forall~\bs \Theta_h \in \bs U_h, \\
\label{eq:FHD-weak-dis-7} (\bs H_h,\bs G_h) + (\phi_h,\div\bs G_h) & = 0 &\forall~\bs G_h \in \bs V_h, \\
\label{eq:FHD-weak-dis-8}\mu_0(\div\bs H_h+\div\bs m_h,r_h) + (\div\bs H_e,r_h) & = 0 &\forall~ r_h \in W_h,
\end{align}
with the initial data 
\begin{equation}
\label{eq:initial-semi}
\bs u_h(\cdot,0) = \pi_h\bs u_0,\quad\bs\omega_h(\cdot,0) = \pi_h\bs \omega_0,\quad \bs m_h(\cdot,0) = I_h^d\bs m_0.
\end{equation}

In what follows let us show the energy-stability of this semi-discrete scheme. 
 To this end, we define the energy $\mathcal E_h(t)$ and the dissipated energy $\mathcal F_h(t)$ 
 at time $t \in [0,T]$  as 
$$
 \mathcal E_h(t) = \rho\|\bs u_h(\cdot,t)\|^2 + \rho\kappa\|\bs\omega_h(\cdot,t)\|^2 + \|\bs m_h(\cdot,t)\|^2 + \mu_0\|\bs H_h(\cdot,t)\|^2
 $$
 and 
 \begin{align*}
 \mathcal F_h(t) =  \eta\|\nabla\bs u_h\|^2 + \eta^\prime \|\nabla\bs\omega_h\|^2 + (\eta^\prime+\lambda^\prime)\|\div\bs\omega_h\|^2 + \frac{1}{\tau}\|\bs m_h\|^2 
 + \frac{\chi_0+(1+\mu_0)\chi_0}{\tau}\|\bs H_h\|^2
  + \sigma\|\div\bs m_h\|^2\\
   + \sigma\|\bs k_h\|^2 
  + \mu_0\sigma\|\div\bs H_h\|^2 
   +\zeta\|\curl\bs u_h - 2\bs\omega_h\|^2,
 \end{align*}
respectively. Since the system \eqref{eq:FHD-weak-dis-1} - \eqref{eq:FHD-weak-dis-8} inherits the structure of the weak formulation \eqref{eq:FHD-weak-1} - \eqref{eq:FHD-weak-8}, by a similar proof as that of Theorem \ref{the:ener}, we easily have the following energy estimate.
 \begin{theorem}
 \label{the:ener-semi-dis}
 Given $\bs H_e \in \bs H^1(\bs H(\div))$, assume $\bs u_h:[0,T]\rightarrow\bs S_h$, $\tilde p_h:[0,T]\rightarrow L_h$, $\bs\omega_h:[0,T]\rightarrow\bs \Sigma_h$, $\bs m_h:[0,T]\rightarrow \bs V_h$, $\bs z_h:[0,T]\rightarrow\bs U_h$, $\bs k_h:[0,T]\rightarrow \bs U_h$, $\bs H_h:[0,T]\rightarrow\bs V_h$ and $\phi_h:[0,T]\rightarrow W_h$ solve the scheme \eqref{eq:FHD-weak-dis-1} - \eqref{eq:FHD-weak-dis-8}. Then the energy inequality
 \begin{equation}\label{eq:ener-semi}
 \mathcal E_h(t) +C_1\int_0^t\mathcal F_h(s) \dd s \leq \mathcal E_h(0)  + C_2\int_0^t (\|\bs H_e(\cdot,s)\|_{d}^2  + \|\partial_t\bs H_e(\cdot,s)\|^2)\dd s
\end{equation}
 holds for all $t \in [0,T]$, where $C_1$ and $C_2$ are positive constants depending only on $\chi_0$, $\mu_0$ and $\tau$.
 \end{theorem}

Now we turn to show the solvability of  the semi-discrete scheme \eqref{eq:FHD-weak-dis-1} - \eqref{eq:FHD-weak-dis-8}. The choice of finite element spaces $\bs V_h$ and $W_h$ implies that, for a given $\bs m_h \in \bs V_h$, there exists a unique pair $(\bs H_h(\bs m_h),\phi_h(\bs m_h))\in \bs V_h \times W_h$ satisfying the equations \eqref{eq:FHD-weak-dis-7} - \eqref{eq:FHD-weak-dis-8}. Equations \eqref{eq:FHD-weak-dis-5} and \eqref{eq:FHD-weak-dis-6} imply $\bs z_h = Q_h^c(\bs u_h \times \bs m_h)$ and $\bs k_h = \curl_h\bs m_h$. 

We introduce the subspace 
$$
\mathfrak K_h^s = \{ \bs v_h \in \bs S_h:~(\div\bs v_h,q_h) = 0,\quad\forall~	q_h\in L_h\}
$$
and denote 
$$
\mathbb V_h = \mathfrak K_h^s \times \bs\Sigma_h \times \bs V_h \times \mathfrak K_h^d,
$$
where $\mathfrak K_h^d$ is given in \eqref{Khd}. 
Define functionals $\mathbb A(\cdot,\cdot):~\mathbb V_h \times \mathbb V_h\rightarrow \mathbb R$ and $\mathbb B(\cdot,\cdot):\mathbb V_h \times \mathbb V_h \rightarrow \mathbb R$ as follows: for  
$$
\mathcal U_h: = (\bs u_h,\bs\omega_h,\bs m_h,\bs H_h)  \in \mathbb V_h,\qquad \mathcal V_h: = (\bs v_h,\bs s_h,\bs F_h,\bs A_h)  \in \mathbb V_h,
$$
\begin{align*}
\mathbb A(\mathcal U_h,\mathcal V_h)  := &\eta(\nabla\bs u_h,\nabla\bs v_h)  + \eta^\prime (\nabla\bs\omega_h,\nabla\bs s_h)  + (\eta^\prime + \lambda^\prime) (\div\bs\omega_h,\div\bs s_h)  + \sigma(\div\bs m_h,\div\bs F_h) \\
& + \sigma(\curl_h\bs m_h,\curl_h\bs F_h) + \frac{1}{\tau}(\bs m_h,\bs F_h)  + \mu_0\sigma(\div\bs H_h,\div \bs A_h) + \zeta(\curl\bs u_h,\curl\bs v_h)  \\
&-2\zeta(\bs\omega_h,\curl\bs v_h) - 2\zeta(\curl\bs u_h,\bs s_h)  + 4\zeta(\bs \omega_h,\bs s_h)    + \frac{\mu_0(1+\chi_0)}{\tau}(\bs H_h,\bs A_h) + \frac{\chi_0}{\tau}(Q_{\mathfrak K_h^d}\bs m_h,\bs F_h),
\end{align*}
\begin{align*}
\mathbb B(\mathcal U_h,\mathcal V_h) := & \rho b(\bs u_h,\bs u_h,\bs v_h) -\mu_0 c(\bs v_h,\bs m_h,\bs H_h) - \frac{\mu_0}{2}(\bs v_h \times \curl_h\bs m_h,\bs H_h)  - \frac{\mu_0}{2}(\bs m_h\times \curl\bs v_h,\bs H_h)  \\
&+ \rho\kappa b(\bs u_h,\bs\omega_h,\bs s_h)- \mu_0(\bs m_h\times \bs H_h,\bs s_h)  - c(\bs u_h,\bs m_h,\bs F_h) -\frac{1}{2}(\bs m_h\times \curl\bs u_h,\bs F_h) \\
&- \frac{1}{2}(\bs u_h\times \curl_h\bs m_h,\bs F_h)   - \frac{1}{2}(\bs u_h\times \bs m_h,\curl_h\bs F_h) -(\bs\omega_h\times \bs m_h,\bs F_h)+ \mu_0c(\bs u_h,\bs m_h,\bs A_h)\\
& + \frac{\mu_0}{2}(\bs m_h\times \curl\bs u_h,\bs A_h)  + \frac{\mu_0}{2}(\bs u_h \times \curl_h\bs m_h,\bs A_h)  + \mu_0(\bs \omega_h \times \bs m_h,\bs A_h).
\end{align*}
Here $Q_{\mathfrak K_h^d}:(L^2(\Omega))^3\rightarrow \mathfrak K_h^d$ is the $L^2$ orthogonal projection operator. 

We define 
$$
\interleave\mathcal V_h\interleave = \left(\eta\|\nabla\bs v_h\|^2 + \eta^\prime\|\nabla\bs s_h\|^2 + \sigma\|\bs m_h\|_n^2 + \mu_0\sigma\|\div\bs H_h\|^2\right)^{1/2}.
$$
It is easy to see that $\mathbb A(\cdot,\cdot)$ is a bilinear form with
\begin{align}
\label{eq:A-coer}\mathbb A(\mathcal V_h,\mathcal V_h) & \geq \interleave \mathcal V_h\interleave^2 & \forall~\mathcal V_h \in \mathbb V_h, \\
\label{eq:A-continu} \mathbb A(\mathcal U_h,\mathcal V_h) & \lesssim \interleave \mathcal U_h\interleave \interleave \mathcal V_h\interleave &\forall~\mathcal U_h,~\mathcal V_h \in \mathbb V_h,
\end{align}
and $\mathbb B(\cdot,\cdot)$ is skew-symmetric operator, i.e.,
$$
\mathbb B(\mathcal V_h,\mathcal V_h) = 0\qquad\qquad\forall~\mathcal V_h \in \mathbb V_h.
$$

Denote 
$$\mathbb D = \hbox{diag}(\rho,\rho\kappa,1,\mu_0)$$
and
$$
\mathcal L(\mathcal V_h) := -\frac{\chi_0}{\mu_0\tau}(Q_{\mathfrak K_h^d}\bs H_e,\bs F_h)-\sigma(\div\bs H_e,\div\bs A_h)-(\partial_t\bs H_e+\frac{1}{\tau}\bs H_e,\bs A_h). 
$$
Then we introduce an auxiliary problem: Find $\mathcal U_h = (\bs u_h,\bs\omega_h,\bs m_h,\bs H_h) \in \mathbb V_h$, such that
\begin{equation}\label{eq:au-semi}
(\mathbb D\partial_t\mathcal U_h ,\mathcal V_h) + \mathbb A(\mathcal U_h,\mathcal V_h) + \mathbb B(\mathcal U_h,\mathcal V_h)= \mathcal L(\mathcal V_h)\qquad\forall~~\mathcal V_h \in \mathbb V_h,
\end{equation}
with the initial value
$$
\mathcal U_h(\cdot,0) = (\bs u_h(\cdot,0),\bs\omega_h(\cdot,0),\bs m_h(\cdot,0),\bs H_h(\cdot,0))  \in \mathbb V_h,
$$
where $\bs H_h(\cdot,0) \in \mathfrak K_h^d$ is determined by the saddle point problem
\begin{equation}\label{eq:H-initial}
\left\{
\begin{array}{ll}
(\bs H_h(\cdot,0),\bs G_h) + (\phi_h(\cdot,0),\div\bs G_h) = 0 & \forall~\bs G_h \in \bs V_h, \\
\mu_0(\div\bs H_h(\cdot,0) +\div\bs m_h(\cdot,0),r_h) = -(\div\bs H_e(\cdot,0),r_h) & \forall~r_h\in W_h. 
\end{array} 
\right.
\end{equation}
We have the following Lemma.
\begin{lemma}\label{lem:eq-sol-semi}
The semi-discrete form \eqref{eq:FHD-weak-dis-1} - \eqref{eq:FHD-weak-dis-8} and the auxiliary problem \eqref{eq:au-semi} are equivalent in the following sense:
\begin{itemize}
\item if $\bs u_h:[0,T]\rightarrow\bs S_h$, $\tilde p_h:[0,T]\rightarrow L_h$, $\bs\omega_h:[0,T]\rightarrow\bs \Sigma_h$, $\bs m_h:[0,T]\rightarrow \bs V_h$, $\bs z_h:[0,T]\rightarrow\bs U_h$, $\bs k_h:[0,T]\rightarrow \bs U_h$, $\bs H_h:[0,T]\rightarrow\bs V_h$ and $\phi_h:[0,T] \rightarrow W_h$ solve \eqref{eq:FHD-weak-dis-1} - \eqref{eq:FHD-weak-dis-8}, then $\mathcal U_h =(\bs u_h,\bs\omega_h,\bs m_h,\bs H_h)^\intercal \in \mathbb V_h$ solves \eqref{eq:au-semi}.
\item if $\mathcal U_h =(\bs u_h,\bs\omega_h,\bs m_h,\bs H_h) \in \mathbb V_h$ solves \eqref{eq:au-semi}, then $\bs u_h$, $\bs\omega_h$, $\bs m_h$, $\bs z_h = Q_h^c(\bs u_h \times \bs m_h)$, $\bs k_h = \curl_h\bs m_h$, a unique $\phi_h \in W_h$ satisfying $\bs H_h = \grad_h\phi_h$, and a unique $\tilde p_h \in L_h$ solve \eqref{eq:FHD-weak-dis-1} - \eqref{eq:FHD-weak-dis-8}.  
\end{itemize}
\end{lemma}
\begin{proof}
We first assume that $\bs u_h\in \bs S_h$, $\tilde p_h\in L_h$, $\bs\omega_h\in \bs \Sigma_h$, $\bs m_h\in \bs V_h$, $\bs z_h\in \bs U_h$, $\bs k_h\in \bs U_h$, $\bs H_h \in \bs V_h$, and $\phi_h\in W_h$ solve \eqref{eq:FHD-weak-dis-1} - \eqref{eq:FHD-weak-dis-8}. Equations \eqref{eq:FHD-weak-dis-2}, \eqref{eq:FHD-weak-dis-5}, \eqref{eq:FHD-weak-dis-6} and \eqref{eq:FHD-weak-dis-7} imply that
$$
\bs u_h \in \mathfrak K_h^s,\quad \bs z_h = Q_h^c(\bs u_h\times\bs m_h),\quad \bs k_h = \curl_h\bs m_h,\quad \bs H_h = \grad_h\phi_h \in \mathfrak K_h^d.
$$
Substitute them into \eqref{eq:FHD-weak-dis-1}, \eqref{eq:FHD-weak-dis-3} and \eqref{eq:FHD-weak-dis-4}, then we get
\begin{align}
\label{eq:exist-dis-pr-1}
\begin{split}
\rho(\partial_t\bs u_h,\bs v_h) + \rho b(\bs u_h,\bs u_h,\bs v_h) + \eta(\nabla\bs u_h,\nabla\bs v_h) + \zeta(\curl\bs v_h,\curl\bs v_h)
 - \mu_0 c(\bs v_h,\bs m_h,\bs H_h)&\\
  - \frac{\mu_0}{2}(\bs v_h\times\curl_h\bs m_h,\bs H_h)
   -\frac{\mu_0}{2}(\bs m_h\times\curl\bs v_h,\bs H_h) - 2\zeta(\bs\omega_h,\curl\bs v_h) &= 0 ,
   \end{split}\\
 \label{eq:exist-dis-pr-2} 
 \begin{split}
 \rho\kappa(\partial_t\bs\omega_h,\bs s_h) +\eta^\prime(\nabla\bs\omega_h,\nabla\bs s_h) + (\eta^\prime+\lambda^\prime)(\div\bs\omega_h,\div\bs s_h) + \rho\kappa b(\bs u_h,\bs \omega_h,\bs s_h)& \\
 -\mu_0(\bs m_h\times\bs H_h,\bs s_h) -2\zeta(\curl\bs u_h,\bs s_h) + 4\zeta(\bs\omega_h,\bs s_h) & = 0,
 \end{split}\\
 \label{eq:exist-dis-pr-3}
 \begin{split}
 (\partial_t\bs m_h,\bs F_h) - c(\bs u_h,\bs m_h,\bs F_h) + \sigma(\div\bs m_h,\div\bs F_h) + \sigma(\curl_h\bs m_h,\curl_h\bs F_h) & \\
 -\frac{1}{2} (\bs m_h\times\curl\bs u_h,\bs F_h) - \frac{1}{2}(\bs u_h\times \curl_h\bs m_h,\bs F_h) - \frac{1}{2}(\bs u_h\times \bs m_h,\curl_h\bs F_h) & \\
 -(\bs \omega_h \times \bs m_h,\bs F_h) + \frac{1}{\tau}(\bs m_h,\bs F_h) - \frac{\chi_0}{\tau}(\bs H_h,\bs F_h) & = 0,
 \end{split}
\end{align}
 for any $\bs v_h \in \mathfrak K_h^s$, $\bs s_h\in \bs\Sigma_h$ and $\bs F_h \in \bs V_h$. Equation \eqref{eq:FHD-weak-dis-8} and the exact sequence \eqref{eq:exact-seq-re} imply
\begin{equation}\label{eq:exist-dis-pr-4}
\mu_0(\bs H_h + \bs m_h,\bs A_h) = -(I_h^d\bs H_e,\bs A_h)\qquad\forall~\bs A_h \in \mathfrak K_h^d
\end{equation}
and
\begin{equation}\label{eq:exist-dis-pr-5}
\mu_0(\partial_t\bs H_h + \partial_t\bs m_h,\bs A_h) = -(I_h^d\partial_t\bs H_e,\bs A_h) \qquad\forall~\bs A_h \in \mathfrak K_h^d.
\end{equation}
Taking $\bs F_h = \bs A_h \in \mathfrak K_h^d$ in \eqref{eq:exist-dis-pr-3}, substituting $(\partial_t\bs m_h,\bs A_h)$ in \eqref{eq:exist-dis-pr-5}, and using \eqref{eq:exist-dis-pr-4} and \eqref{eq:FHD-weak-dis-8} with $r_h = \div\bs A_h \in W_h$, we get
\begin{equation}\label{eq:exist-dis-pr-6}
\begin{split}
\mu_0(\partial_t\bs H_h,\bs A_h) + \mu_0c(\bs u_h,\bs m_h,\bs A_h) +\sigma\mu_0 (\div\bs H_h,\div\bs A_h) 
+ \frac{\mu_0}{2}(\bs m_h\times\curl\bs u_h,\bs A_h) + \frac{\mu_0}{2}(\bs u_h\times \curl_h\bs m_h,\bs A_h)  \\
+ \mu_0(\bs\omega_h \times \bs m_h,\bs A_h)
+ \frac{(1+\chi_0)\mu_0}{\tau}(\bs H_h,\bs A_h) = -\frac{1}{\tau}(I_h^d\bs H_e,\bs A_h)
- \sigma(\div\bs H_e,\div\bs A_h) -  (I_h^d\partial_t\bs H_e,\bs A_h).
\end{split}
\end{equation}
Equation \eqref{eq:exist-dis-pr-4} also implies
$$
\bs H_h =- Q_{\mathfrak K_h^d}\left( \bs m_h  - \frac{1}{\mu_0}I_h^d \bs H_e\right).
$$
Substitute the above equation into the term $(\bs H_h,\bs F_h)$ of \eqref{eq:exist-dis-pr-3}, and we see that $\mathcal U_h = (\bs u_h,\bs\omega_h,$ $\bs m_h,\bs H_h)  \in \mathbb V_h$ solves \eqref{eq:au-semi}.

On the other hand, we assume that $\mathcal U_h = (\bs u_h,\bs\omega_h,\bs m_h,\bs H_h)  \in \mathbb V_h$ solves \eqref{eq:au-semi}, then it is easy to check that $\bs u_h\in \mathfrak K_h^s$, $\bs z_h = Q_h^c(\bs u_h\times\bs m_h)$ and $\bs k_h = \curl_h\bs m_h$ satisfy \eqref{eq:FHD-weak-dis-2}, \eqref{eq:FHD-weak-dis-5} and \eqref{eq:FHD-weak-dis-6}, respectively.

 Taking $\mathcal V_h = (\bs 0,\bs 0,\bs 0,\bs A_h)  \in \mathbb V_h$  in \eqref{eq:au-semi}, we obtain \eqref{eq:exist-dis-pr-6}.
Taking $\mathcal V_h = (\bs 0,\bs 0,\bs A_h,\bs 0)  \in \mathbb V_h$ in \eqref{eq:au-semi} with $\bs A_h \in \mathfrak K_h^d$, we get
\begin{equation*}
\begin{split}
(\partial_t\bs m_h,\bs A_h) - c(\bs u_h,\bs m_h,\bs A_h) + \sigma(\div\bs m_h,\div\bs A_h) 
 -\frac{1}{2} (\bs m_h\times\curl\bs u_h,\bs A_h) - \frac{1}{2}(\bs u_h\times \curl_h\bs m_h,\bs A_h) 
 -(\bs \omega_h \times \bs m_h,\bs A_h) &\\+ \frac{1}{\tau}(\bs m_h,\bs A_h)
  - \frac{\chi_0}{\tau}(\bs H_h,\bs A_h) 
  = -\frac{\chi_0}{\mu_0\tau}(\bs H_e, \bs A_h). &
\end{split}
\end{equation*}
Substitute this into \eqref{eq:exist-dis-pr-6}, and  we obtain 
\begin{equation}\label{eq:exist-dis-pr-8}
\begin{split}
\mu_0(\partial_t(\bs H_h+\bs m_h),\bs A_h) + \sigma\mu_0(\div(\bs H_h+\bs m_h),\div\bs A_h)
 + \frac{(1+\chi_0)\mu_0}{\tau}(\bs H_h+\bs m_h,\bs A_h) = -\frac{1+\chi_0}{\tau}(\bs H_e,\bs A_h) \\
 - \sigma(\div\bs H_e,\div\bs A_h) - (\partial_t\bs H_e,\bs A_h).
\end{split}
\end{equation}
 We easily see that \eqref{eq:exist-dis-pr-8} has a unique solution in $\mathfrak K_h^d$ with the initial condition $\bs H_h(\cdot,0) + Q_{\mathfrak K_h^d}\bs m_h(\cdot,0)$. Also we can    check that if $\bs H_h$ and $\bs m_h$ satisfy \eqref{eq:FHD-weak-dis-8}, then $\bs H_h +Q_{\mathfrak K_h^d}\bs m_h \in \mathfrak K_h^d$ solves \eqref{eq:exist-dis-pr-8}. Thus,   \eqref{eq:exist-dis-pr-8} is equivalent to \eqref{eq:FHD-weak-dis-8}. 
 
 Since $\bs H_h \in \mathfrak K_h^d$, we have \eqref{eq:FHD-weak-dis-7}. 
 
 Taking $\mathcal V_h = (\bs s_h,\bs 0,\bs 0,\bs 0)^\intercal\in \mathbb V_h$ in \eqref{eq:au-semi} and using the inf-sup condition \eqref{eq:inf-sup}, we see that there exists a unique $p_h\in L_h$, such that $\bs u_h$ and $p_h\in L_h$ satisfy \eqref{eq:FHD-weak-dis-1} - \eqref{eq:FHD-weak-dis-2}. 
 
 Taking $\mathcal V_h = (\bs 0,\bs s_h,\bs 0,\bs 0)^\intercal \in \mathbb V_h$ in \eqref{eq:au-semi}, we obtain \eqref{eq:FHD-weak-dis-3}. 
 
 Taking $\mathcal V_h = (\bs 0,\bs 0,\bs F_h,\bs 0)^\intercal \in \mathbb V_h$ in \eqref{eq:au-semi} and applying \eqref{eq:FHD-weak-dis-8}, we obtain \eqref{eq:FHD-weak-dis-4}.
 
 As a result, the desired conclusion follows.
\end{proof}

In the reset of this subsection, we will focus on the existence of solution to the auxiliary problem \eqref{eq:au-semi}. The definitions of the bilinear form $\mathbb A(\cdot,\cdot)$ and the linear form $\mathcal L(\cdot)$ imply they are Lipschitz continuous  on the space $L^2(\mathbb V_h)$. We define a subspace of $L^2(\mathbb V_h)$ as follows: 
$$
\mathcal L^2(\mathbb V_h) = \{ \mathcal U_h = (\bs u_h,\bs\omega_h,\bs m_h,\bs H_h)  \in L^2(\mathbb V_h):~\mathcal U_h \text{ satisfies }\eqref{eq:ener-semi} \text{ with }\bs k_h = \curl_h\bs m_h\}.
$$
We claim that $\mathcal B(\cdot,\cdot)$ is   Lipschitz continuous on   $\mathcal L^2(\mathbb V_h)$. In fact, for any $\mathcal U_i = (\bs u_i,\bs\omega_i,\bs m_i,\bs H_i)^\intercal \in \mathcal L^2(\mathbb V_h)$ ($i = 1,2$) and $\mathcal V_h = (\bs v_h,\bs s_h,\bs F_h,\bs A_h)^\intercal \in \mathbb V_h$, we have
\begin{align*}
\mathbb B(\mathcal U_1 - \mathcal U_2,\mathcal V_h) \lesssim &(\|\nabla\bs u_1\| +\|\nabla\bs u_2\|)\|\nabla(\bs u_1 - \bs u_2)\| \|\nabla\bs v_h\|  + \|\div\bs H_1\| \|\bs m_1-\bs m_2\|_n \|\nabla\bs v_h\| \\
&+ \|\bs m_2\|_n\|\div(\bs H_1 - \bs H_2)\|\|\nabla\bs v_h\|  + (\|\nabla\bs \omega_1\| +\|\nabla\bs u_2\|)\|\nabla(\bs\omega_1 - \bs\omega_2)\| \|\nabla\bs s_h\| \\
& +\|\bs H_1\|\|\bs m_1 - \bs m_2\|_n \|\nabla\bs s_h\|  + \|\bs m_2\|\|\div(\bs H_1 - \bs H_2)\|\|\nabla\bs s_h\| + \|\bs m_1\|_n\|\nabla(\bs u_1 - \bs u_2)\|\|\bs F_h\|_n \\
& + (\|\nabla\bs u_2\|+\|\bs \omega_2\|)\|\bs m_1 - \bs m_2\|_n\|\bs F_h\|_n +\|\bs m_1\|\|\nabla(\bs\omega_1 - \bs\omega_2)\|\|\bs F_h\|_n\\
& +\left[ (\|\nabla\bs u_2\| +\|\bs\omega_2\|)\|\bs m_1 - \bs m_2\|_n + \|\bs m_1\|\|\nabla(\bs\omega_1 - \bs\omega_2)\|\right]\|\div\bs A_h\| \\
& + \|\bs m_1\|_n \|\nabla(\bs u_1 - \bs u_2)\|\|\div\bs A_h\|\big)\\
\lesssim & (\|\nabla\bs u_1\| + \|\nabla\bs u_2\| + \|\nabla\bs\omega_1 \| + \|\bs m_1\|_n +\|\bs m_2\|_n +\|\div\bs H_1\|)\interleave \mathcal U_1 - \mathcal U_2\interleave \interleave \mathcal V_h \interleave,
\end{align*}
which means $\mathbb B(\cdot,\cdot)$ is Lipschitz continuous with respect to the first variable with Lipschitz constant
\begin{align*}
Lip & = C\left( \int_0^T(\|\nabla\bs u_1\|^2 + \|\bs u_2\|^2 +\|\bs\omega_1\|^2 + \|\bs m_1\|_n^2+\|\bs m_2\|_n^2 + \|\div\bs H_1\|^2)\dd s\right)^{1/2} \\
& \leq C\left( \mathcal E_h(0) + \int_0^T(\|\bs H_e\|_{d}^2 +\|\partial_t\bs H_e\|)\dd s \right)^{1/2}.
\end{align*}
We have the following existence result.
\begin{theorem}
Under  Assumptions \ref{ass:He} and \ref{ass:initial-data}, there exist $\bs u_h \in L^2(\bs V_h)$, $\tilde p_h \in L^2(L_h)$, $\bs\omega_h \in L^2(\bs\Sigma_h)$, $\bs m_h \in L^2(\bs V_h)$, $\bs z_h \in\bs U_h$, $\bs k_h\in \bs U_h$, $\bs H_h \in L^2(\bs V_h)$ and $\phi_h \in W_h$ solve the semi-discrete problem \eqref{eq:FHD-weak-dis-1} - \eqref{eq:FHD-weak-dis-8}. 
\end{theorem}
 \begin{proof}
By Lemma \ref{lem:eq-sol-semi}, we only need to prove that the auxiliary problem \eqref{eq:au-semi} has a unique solution. For any $t,~s\in (0,T]$, we can rewrite the auxiliary problem \eqref{eq:au-semi} as
\begin{equation}
\label{eq:pr-sol-1}
(\mathbb D\mathcal U_h(t),\mathcal V_h) = (\mathbb D\mathcal U_h(s),\mathcal V_h) + \int_s^t \left( \mathbb A(\mathcal U_h(l),\mathcal V_h) + \mathbb B(\mathcal U_h(l),\mathcal V_h) + \mathcal L(\mathcal V_h)\right)\dd l,
\end{equation}
where we recall that
 $
\mathcal U_h = (\bs u_h,\bs\omega_h,\bs m_h,\bs H_h)  \in \mathbb V_h$ and $ \mathcal V_h = (\bs v_h,\bs s_h,\bs F_h,\bs A_h)  \in \mathbb V_h.
$

Define $\mathbb T:~\mathbb V_h \to \mathbb V_h$ as
$$
(\mathbb D\mathbb T \mathcal U_h(t),\mathcal V_h) = (\mathbb D\mathcal U_h(s),\mathcal V_h) + \int_s^t((\mathbb A(\mathcal U_h(l)),\mathcal V_h) + \mathbb B(\mathcal U_h(l),\mathcal V_h) + \mathcal L(\mathcal V_h))\dd l.
$$
It is easy to check that $\mathbb T$ is a contraction map when $|t-s|$ is small enough.
Thus, the fixed point theorem implies \eqref{eq:pr-sol-1} has a local solution $\mathcal U_h $ in $L^2(\mathbb V_h)$, and Lemma \ref{lem:eq-sol-semi} implies $\bs u_h$, $\bs\omega_h$, $\bs m_h$, $\bs z_h = Q_h^c(\bs u_h \times \bs m_h)$, $\bs k_h = \curl_h\bs m_h$, a unique $\phi_h \in W_h$ satisfying $\bs H_h = \grad_h\phi_h$, and a unique $\tilde p_h \in L_h$   also solve the semi-discrete problem \eqref{eq:FHD-weak-dis-1} - \eqref{eq:FHD-weak-dis-8} and hence they form a local solution of this problem. 
Finally, the   energy estimate in Theorem \ref{the:ener-semi-dis}  means the local solution of \eqref{eq:FHD-weak-dis-1} - \eqref{eq:FHD-weak-dis-8} is a global one.
\end{proof}

Now we turn to discuss the uniqueness of the semi-discrete solution. We make the following assumption.
\begin{assum}\label{ass:semi}
There exist constants $h_0>0$ and $M_0>0$ such that for any $0<h \leq h_0$ and any solution $(\bs u_h,~\tilde p_h,~\bs\omega_h,~\bs m_h,~\bs H_h,~\bs z_h,~\bs k_h,~\phi_h)$ of \eqref{eq:FHD-weak-dis-1} - \eqref{eq:FHD-weak-dis-8}, there holds
$$
\int_0^T(\|\nabla\bs u_h\|^4 + \|\nabla\bs\omega_h\|^4 +\|\bs m_h\|_n^4  +\|\bs H_h\|_n^4)\dd s \leq M_0.
$$ 
\end{assum}

\begin{theorem}\label{the:unique-solu-semi-dis}
Under  Assumption \ref{ass:semi}, the semi-discrete scheme \eqref{eq:FHD-weak-dis-1} - \eqref{eq:FHD-weak-dis-8} has a unique solution when the spatial mesh size satisfies   $0<h \leq h_0$.
\end{theorem}
\begin{proof}
Let $(\bs u_{hi},~\tilde p_{hi},~\bs\omega_{hi},~\bs m_{hi},~\bs H_{hi},~\bs z_{hi},~\bs k_{hi},~\bs\phi_{hi})$ ($i = 1,~2$) be two solutions of the semi-discrete scheme \eqref{eq:FHD-weak-dis-1} - \eqref{eq:FHD-weak-dis-8} with the same initial data $(\bs u_{h}(\cdot,0),~\bs\omega_{h}(\cdot,0),~\bs m_h(\cdot,0))$. Set $d_v: = v_{h1}  - v_{h2}$ with $v = \bs u,~\tilde p,~\bs\omega,~\bs m,~\bs H,~\bs z,~\bs k$ or $\phi$, we have
\begin{equation}\label{eq:pr-uni-1}
\begin{split}
\rho (\partial_td_u,\bs v_h)  + \eta(\nabla d_u,\nabla\bs v_h) + \zeta(\curl d_u,\curl\bs v_h) - (d_p,\div\bs v_h)   - 2\zeta(d_\omega,\curl\bs v_h) = -\rho d_{11}(\bs v_h) \\ + \mu_0 d_{12}(\bs v_h) + \frac{\mu_0}{2} d_{13}(\bs v_h)   + \frac{\mu_0}{2}d_{14}(\bs v_h),
\end{split}
\end{equation}
\begin{equation}\label{eq:pr-uni-2}
(\div d_u,q_h) = 0,
\end{equation}
\begin{equation}\label{eq:pr-uni-3}
\begin{split}
\rho\kappa (\partial_t d_\omega,\bs s_h) +\eta^\prime(\nabla d_\omega,\nabla\bs s_h) + (\eta^\prime + \lambda^\prime)(\div d_\omega,\div\bs s_h) - 2\zeta(\curl d_u,\bs s_h) 
+ 4\zeta(d_\omega,\bs s_h) = \mu_0d_{31}(\bs s_h)  
-\rho\kappa d_{32}(\bs s_h),
\end{split}
\end{equation}
\begin{equation}\label{eq:pr-uni-4}
\begin{split}
(\partial_t d_m,\bs F_h)  + \sigma(\div d_m,\div\bs F_h) + \sigma(\curl d_k,\bs F_h) - \frac{1}{2}(\curl d_z,\bs F_h) 
 + \frac{1}{\tau}(d_m,\bs F_h) - \frac{\chi_0}{\tau}(d_H,\bs F_h) \\= d_{41}(\bs F_h)   + \frac{1}{2}d_{42}(\bs F_h)  
 + \frac{1}{2} d_{43}(\bs F_h)   + d_{44}(\bs F_h) ,
\end{split}
\end{equation}
\begin{equation}\label{eq:pr-uni-5}
(d_z,\bs\Lambda_h) = d_{51}(\Lambda_h) ,
\end{equation}
\begin{equation}\label{eq:pr-uni-6}
(d_k,\bs\Theta_h) - (d_m,\curl\bs\Theta_h) = 0,
\end{equation}
\begin{equation}\label{eq:pr-uni-7}
(d_H,\bs G_h) +(d_\phi,\div\bs G_h) = 0,
\end{equation}
\begin{equation}\label{eq:pr-uni-8}
(\div d_H + \div d_m,r_h) = 0,
\end{equation}
  for any $\bs v_h\in \bs S_h$, $q_h \in L_h$, $\bs s_h \in \bs \Sigma_h$, $\bs F_h \in \bs V_h$, $\bs\Lambda_h\in \bs U_h$, $\bs\Theta_h \in \bs U_h$, $\bs G_h \in \bs V_h$ and $r_h \in W_h$,
with initial data
$$
d_u(\cdot,0) = \bs 0,\qquad d_\omega(\cdot,0) = \bs 0,\qquad d_m(\cdot,0) = \bs 0.
$$
Here
$$
\begin{array}{ll}
 d_{11}(\bs v_h) := b(\bs u_{h1},\bs u_{h1},\bs v_h) - b(\bs u_{h2},\bs u_{h2},\bs v_h), &\quad d_{12}(\bs v_h) := c(\bs v_h,\bs m_{h1},\bs H_{h1}) - c(\bs v_h,\bs m_{h2},\bs H_{h2}), \\
d_{13}(\bs v_h) := (\bs v_h\times \bs k_{h1},\bs H_{h1}) - (\bs v_h\times \bs k_{h2},\bs H_{h2}),&\quad d_{14}(\bs v_h): = (\bs m_{h1}\times\curl\bs v_h,\bs H_{h1}) - (\bs m_{h2}\times\curl\bs v_h,\bs H_{h2}),\\
 d_{31}(\bs s_h): =  (\bs m_{h1}\times \bs H_{h1},\bs s_h) -(\bs m_{h2}\times \bs H_{h2},\bs s_h),&\quad d_{32}(\bs s_h) := b(\bs u_{h1},\bs\omega_{h1},\bs s_h) - b(\bs u_{h2},\bs\omega_{h2},\bs s_h),\\
 d_{41}(\bs F_h) := c(\bs u_{h1},\bs m_{h1},\bs F_h) - c(\bs u_{h2},\bs m_{h2},\bs F_h),&\quad d_{42}(\bs F_h) := (\bs m_{h1} \times\curl\bs u_{h1},\bs F_h) - (\bs m_{h2} \times\curl\bs u_{h2},\bs F_h), \\
d_{43}(\bs F_h): = (\bs u_{h1}\times\bs k_{h1},\bs F_h) - (\bs u_{h2}\times\bs k_{h2},\bs F_h) ,&\quad d_{44}(\bs F_h): = (\bs \omega_{h1}\times \bs m_{h1},\bs F_h) - (\bs\omega_{h2} \times \bs m_{h2},\bs F_h),\\
 d_{51}(\bs\Lambda_h) := (\bs u_{h1}\times \bs m_{h1},\bs\Lambda_h) - (\bs u_{h2}\times \bs m_{h2},\bs\Lambda_h).&
\end{array}
$$

Taking $\bs v_h = d_u$ in \eqref{eq:pr-uni-1}, $\bs s_h = d_\omega$  in \eqref{eq:pr-uni-3} and $\bs F_h = d_m$  in \eqref{eq:pr-uni-4}, and adding the resulting equations together, we get
\begin{equation}\label{eq:pr-uni-9}
\begin{split}
 \frac{1}{2}\frac{\dd}{\dd t}(\rho\|d_u\|^2 + \rho\kappa\|d_\omega\|^2 + \|d_m\|^2) + \eta|d_u|_1^2 +\eta^\prime|d_\omega|_1^2 + \sigma\|\div d_m\|^2 + \sigma\|d_k\|^2 \\
+ (\eta^\prime+\lambda^\prime)\|\div d_\omega\|^2 + \frac{1}{\tau}\|d_m\|^2 + \frac{\chi_0}{\tau}\| d_H\|^2  +\zeta\|\curl d_u-2d_\omega\|^2 \\
= \frac{1}{2}(\curl d_z,d_m) -\rho d_{11}(d_u) 
 + \mu_0d_{12}(d_u)
+ \frac{\mu_0}{2} d_{13}(d_u) 
 + \frac{\mu_0}{2}d_{14}(d_u) 
 + \mu_0 d_{31}(d_\omega) \\
 -\rho\kappa d_{32}(d_\omega) 
 + d_{41}(d_m) 
 + \frac{1}{2}d_{42}(d_m) 
 + \frac{1}{2}d_{43}(d_m) 
 + d_{44}(d_m).
 \end{split}
\end{equation}
Taking $\bs F_h = d_H$ in \eqref{eq:pr-uni-4}, we have
\begin{align*}
\frac{1}{2}\frac{\dd}{\dd t}\|d_H\|^2 + \sigma\|\div d_H\|^2 + \frac{1+\chi_0}{\tau}\|d_H\|^2 = -d_{41}(d_H) 
 - \frac{1}{2}d_{42}(d_H) 
 - \frac{1}{2}d_{43}(d_H) 
 - d_{44}(d_H).
\end{align*}
Multiplying the above equation with $\mu_0$, and adding the resulting equation with \eqref{eq:pr-uni-9}, we obtain
\begin{align*}
& \frac{1}{2}\frac{\dd}{\dd t}(\rho\|d_u\|^2 + \rho\kappa\|d_\omega\|^2 + \|d_m\|^2 + \mu_0\|d_H\|^2 ) + \eta|d_u|_1^2 + \eta^\prime|d_\omega|_1^2 + \sigma\|\div d_m\|^2 
+ \sigma\|d_k\|^2 \\
&\quad + \sigma\mu_0 \|\div d_H\|^2 
 + (\eta^\prime + \lambda^\prime)\|\div d_\omega\|^2 + \frac{1}{\tau}\|d_m\|^2  
+ \frac{\chi_0 + (1+\chi_0)\mu_0}{\tau} \|d_H\|^2 +\zeta\|\curl d_u - 2d_\omega\|^2 \\
= & \frac{1}{2}(\curl d_z,d_m) -\rho d_{11}(d_u) 
 + \mu_0d_{12}(d_u)
+ \frac{\mu_0}{2} d_{13}(d_u) 
 + \frac{\mu_0}{2}d_{14}(d_u) 
 + \mu_0 d_{31}(d_\omega)  -\rho\kappa d_{32}(d_\omega) 
 + d_{41}(d_m) \\
 &
 + \frac{1}{2}d_{42}(d_m) 
 + \frac{1}{2}d_{43}(d_m) 
 + d_{44}(d_m)  - \mu_0d_{41}(d_H) -\frac{\mu_0}{2}d_{42}(d_H) - \frac{\mu_0}{2} d_{43}(d_H) - \mu_0d_{44}(d_H),
\end{align*}
which, together with the  inequalities
\begin{align*}
\rho |d_{11}(d_u)| & = \rho b(d_u,\bs u_{h2},d_u) \leq \frac{\eta}{16}|d_u|_1^2 + C\|\nabla\bs u_{h2}\|^4\|d_u\|^2, \\
\mu_0(d_{12}|d_u)-d_{41}(d_H)| & \leq \frac{\sigma}{18}\|d_m\|_n^2 + \frac{\eta}{16}|d_u|_1^2 + \frac{\mu_0\sigma}{8}\|\div d_H\|^2 + C\|\div\bs H_{h1}\|^4\|d_m\|^2\\
&\quad + C\|\bs H_{h1}\|_n^4\|d_u\|^2 + C\|\nabla\bs u_{h1}\|^4 (\|d_m\|^2+\|d_H\|^2),\\
\frac{\mu_0}{2}|d_{13}(d_u)-d_{43}(d_H)| & \leq \frac{\eta}{16}|d_u|_1^2+\frac{\sigma}{18}\|d_k\|^2 + \frac{\mu_0\sigma}{8}\|\div d_H\|^2  + C\|\bs H_{h2}\|_n^4\|d_u\|^2 + C\|\nabla\bs u_{h2}\|^4\|d_H\|^2,\\
\frac{\mu_0}{2}|d_{14}(d_u) - d_{42}(d_H)| &\leq \frac{\eta}{16}|d_u|_1^2 + \frac{\sigma}{18}\|d_m\|_n^2 + \frac{\mu_0\sigma}{8}\|\div d_H\|^2   + C(\|\nabla\bs u_{h1}\|^4 + \|\bs H_{h1}\|_n^4)\|d_m\|^2,\\
\mu_0|d_{31}(d_\omega)| & \leq \frac{\eta^\prime}{6}|d_\omega|_1^2 + \frac{\sigma}{18}\|d_m\|_n^2 + \frac{\mu_0\sigma}{8}\|\div d_H\|^2   + C(\|\bs H_{h1}\|^4 + \|\bs m_{h1}\|^4)\|d_\omega\|^2, \\
\rho\kappa |d_{32}(d_\omega)| & \leq \frac{\eta}{16}|d_u|_1^2 + \frac{\eta^\prime}{6}|d_\omega|_1^2 + C\|\bs\omega_{h1}\|^2(\|d_\omega\|^2 + \|d_u\|^2),\\
|d_{41}(d_m)| & \leq \frac{\eta}{16}|d_u|_1^2 + \frac{\sigma}{18}\|d_m\|_n^2  + C\|\bs m_{h1}\|_n^4\|d_u\|^2, \\
\frac{1}{2}|d_{42}(d_m)| & \leq \frac{\eta}{16}|d_u|_1^2 + \frac{\sigma}{18}\|d_m\|_n^2 + C\|\bs m_{h2}\|_n^4\|d_m\|^2, \\
\frac{1}{2}|(\curl d_z,d_m) - d_{43}(d_m)| & \leq \frac{\eta}{16}|d_u|_1^2 + \frac{\sigma}{18}\|d_k\|^2 + \frac{\sigma}{18}\|d_m\|_n^2  + C(\|\bs m_{h1}\|^4_n +\|\bs k_{h1}\|^4)\|d_u\|^2,\\
|d_{44}(d_m)| & \leq \frac{\eta^\prime}{6}|d_\omega|_1^2 + \frac{\sigma}{18}\|d_m\|_n^2 + C\|\bs m_{h1}\|^4 \|d_m\|^2,
\end{align*}
gives
\begin{align*}
&\frac{\dd}{\dd t}(\rho\|d_u\|^2 + \rho\kappa\|d_\omega\|^2 + \|d_m\|^2 + \mu_0\|d_H\|^2 ) + \eta|d_u|_1^2 + \eta^\prime|d_\omega|_1^2 + \sigma\|\div d_m\|^2 
+ \sigma\|d_k\|^2 + \sigma\mu_0 \|\div d_H\|^2\\
&\qquad   + (\eta^\prime + \lambda^\prime)\|\div d_\omega\|^2 + \frac{1}{\tau}\|d_m\|^2 + \frac{\chi_0 + (1+\chi_0)\mu_0}{\tau} \|d_H\|^2 +\zeta\|\curl d_u - 2d_\omega\|^2 \\
\leq & C(\|\nabla\bs u_{h2}\|^4 +\|\nabla\bs u_{h1}\|^4 + \|\bs H_{h1}\|_n^4  +\|\bs H_{h2}\|_n^4 + \|\bs \omega_{h,1}\|^4 + \|\bs m_{h,1}\|_n^4 +\|\bs m_{h2}\|_n^4 ) \\
& \times (\rho\|d_u\|^2 + \rho\kappa\|d_\omega\|^2 + \mu_0\|d_m\|^2 + \|d_H\|^2).
\end{align*}
Finally, using the Gr\"owall inequality 
 we get   
$$
d_u = 0,\quad d_\omega = 0,\quad d_m = 0,\quad d_H = 0.
$$
This completes the proof.
\end{proof}

 \subsection{Error estimates}
 To derive   error estimates for the semi-discrete formulation \eqref{eq:FHD-weak-dis-1} - \eqref{eq:FHD-weak-dis-8}, we introduce some notations:
 $$
 \begin{array}{llll}
  \xi_u: = \pi_h^S\bs u - \bs u_h, & \quad \xi_\omega := \pi_h^\Sigma \bs\omega - \bs\omega_h,& \quad\xi_m := I_h^d\bs m - \bs m_h,& \quad\xi_H := I_h^d\bs H - \bs H_h,\\
  \xi_p := Q_h^L p - p_h, & \quad \xi_z := I_h^c\bs z - \bs z_h,& \quad\xi_k := I_h^c\bs k - \bs k_h,& \quad\xi_\phi := Q_h^W\phi - \phi_h,\\
  \theta_u := \pi_h^S\bs u - \bs u,& \quad \theta_\omega := \pi_h^\Sigma\bs\omega - \bs\omega,& \quad \theta_m := I_h^d\bs m - \bs m,&\quad \theta_H := I_h^d\bs H - \bs H, \\
  \theta_p := Q_h^Lp - p,&\quad\theta_z := I_h^c\bs z - \bs z,&\quad \theta_k := I_h^c\bs k - \bs k,&\quad\theta_\phi := Q_h^W\phi - \phi.
 \end{array}
 $$
From the continuous formulations \eqref{eq:FHD-weak-1},\eqref{eq:FHD-weak-3} and \eqref{eq:FHD-weak-4} and the discrete  formulations \eqref{eq:FHD-weak-dis-1}, \eqref{eq:FHD-weak-dis-3} and \eqref{eq:FHD-weak-dis-4}, we easily obtain
\begin{equation}\label{eq:err-semi-1}
\begin{split}
\rho(\partial_t\xi_u,\bs v_h) + \eta(\nabla\xi_u,\nabla\bs v_h) + \zeta (\curl\xi_u,\curl\bs v_h) - (\xi_p,\div\bs v_h) 
- 2\zeta(\xi_\omega,\curl \bs v_h) = -\rho f_{11}(\bs v_h) \\+ \mu_0 f_{12}(\bs v_h) + \frac{\mu_0}{2} f_{13}(\bs v_h) 
+ \frac{\mu_0}{2} f_{14}(\bs v_h) + e_1(\bs v_h)
\end{split}\qquad \forall~\bs v_h \in \bs S_h,
\end{equation}
\begin{equation}\label{eq:err-semi-2}
\begin{split}
\rho\kappa(\partial_t\xi_\omega,\bs s_h) + \eta^\prime (\nabla\xi_\omega,\nabla\bs s_h) + (\eta^\prime + \lambda^\prime) (\div\xi_\omega,\div\bs s_h) - 2\zeta (\curl \xi_u,\bs s_h) 
+ 4\zeta(\xi_\omega,\bs s_h) \\= -\rho\kappa f_{21}(\bs s_h) + \mu_0 f_{22}(\bs s_h) + e_2(s_h)
\end{split}\qquad\qquad\forall~\bs s_h \in \bs\Sigma_h,
\end{equation}
and
\begin{equation}\label{eq:err-semi-3}
\begin{split}
(\partial_t\xi_m,\bs F_h) + \sigma(\div\xi_m,\div\bs F_h) + \sigma(\curl\xi_k,\bs F_h)
 - \frac{1}{2}(\curl\xi_z,\bs F_h) + \frac{1}{\tau}(\xi_m,\bs F_h) \\
 - \frac{\chi_0}{\tau}(\xi_H,\bs F_h) 
= f_{31}(\bs F_h) + \frac{1}{2} f_{32}(\bs F_h) + \frac{1}{2} f_{33}(\bs F_h) 
+ f_{34}(\bs F_h) + e_3(\bs F_h)
\end{split}\qquad\qquad\qquad\qquad\forall~\bs F_h \in \bs V_h,
\end{equation}
with
\begin{align*}
& f_{11}(\bs v_h) := b(\bs u,\bs u,\bs v_h) - b(\bs u_h,\bs u_h,\bs v_h),\qquad\qquad  f_{12}(\bs v_h) := c(\bs v_h,\bs m,\bs H) - c(\bs v_h,\bs m_h,\bs H_h), \\
 & f_{13}(\bs v_h) := (\bs v_h\times\bs k,\bs H) - (\bs v_h\times\bs k_h,\bs H_h),\quad\quad\,\,\, f_{14}(\bs v_h) := (\bs m\times\curl\bs v_h,\bs H) - (\bs m_h\times\curl\bs v_h,\bs H_h),\\
& f_{21}(\bs s_h): = b(\bs u,\bs\omega,\bs s_h) - b(\bs u_h,\bs\omega_h,\bs s_h),\qquad\quad\,\,\,\, f_{22}(\bs s_h):= (\bs m\times \bs H,\bs s_h) - (\bs m_h\times \bs H_h,\bs s_h),\\
& f_{31}(\bs F_h) := c(\bs u,\bs m,\bs F_h) - c(\bs u_h,\bs m_h,\bs F_h),\qquad\,\,\,  f_{32}(\bs F_h) := (\bs m\times \curl\bs u,\bs F_h) - (\bs m_h\times \curl\bs u_h,\bs F_h),\\
& f_{33}(\bs F_h) := (\bs u\times\bs k,\bs F_h) - (\bs u_h\times\bs k_h,\bs F_h),\qquad\,\, f_{34}(\bs F_h) := (\bs\omega\times\bs m,\bs F_h) - (\bs\omega_h\times\bs m_h,\bs F_h),\\
& e_1(\bs v_h) := \rho(\partial_t\theta_u,\bs v_h) + \eta(\nabla\theta_u,\nabla\bs v_h) + \zeta(\curl\theta_u,\curl\bs v_h) - (\theta_p,\div\bs v_h) - 2\zeta(\theta_\omega,\curl\bs v_h),\\
& e_2(\bs s_h) := \rho\kappa(\partial_t\theta_\omega,\bs s_h) + \eta^\prime (\nabla\theta_\omega,\nabla\bs v_h) + (\eta^\prime + \lambda^\prime) (\div\theta_\omega,\div\bs s_h) - 2\zeta(\curl\theta_\omega,\bs s_h) + 4\zeta(\theta_\omega,\bs s_h),\\
& e_3(\bs F_h) := (\partial_t\theta_m,\bs F_h) + \sigma(\curl\theta_k,\bs F_h)+  \frac{1}{\tau}(\theta_m,\bs F_h) - \frac{\chi_0}{\tau}(\theta_H,\bs F_h) - \frac{1}{2}(\curl\theta_z,\bs F_h).
\end{align*}

Define
\begin{align*}
& J_1 := e_1(\xi_u) + e_2(\xi_\omega) + e_3(\xi_m) - (\div\theta_u,\xi_p) + \frac{\chi_0}{\tau}(\theta_H,\xi_H)  +\frac{\chi_0}{\tau} (\theta_H,\xi_m)+\sigma(\theta_k,\xi_k),\\
& J_2 := - e_3(\xi_H) + \frac{\chi_0}{\tau}(\theta_H,\xi_H) + \frac{\chi_0}{\tau}(\theta_H,\xi_m),\\
& J_3 := \rho f_{11}(\xi_u) + \rho\kappa f_{21}(\xi_\omega),\\
& J_4 := \mu_0 \left(f_{12}(\xi_u) - f_{31}(\xi_H) + f_{22}(\xi_\omega) - f_{34}(\xi_H) \right),\\
& J_5 := \frac{\mu_0}{2} \left(f_{13}(\xi_u) - f_{33}(\xi_H) + f_{14}(\xi_u) - f_{32}(\xi_H) \right),\\
& J_6 := \frac{1}{2}(\curl\xi_z,\xi_m) + f_{31}(\xi_m) + \frac{1}{2}f_{32}(\xi_m) + \frac{1}{2} f_{33}(\xi_m) + f_{34}(\xi_m).
\end{align*}
We have the following identity of $\xi$ and $\theta$:
\begin{lemma}\label{lem:err-ID}
There holds
\begin{align*}
\frac{1}{2}\frac{\dd}{\dd t}(\rho\|\xi_u\|^2 + \rho\kappa\|\xi_\omega\|^2 + \|\xi_m\|^2 + \mu_0\|\xi_H\|^2) + \eta\|\nabla\xi_u\|^2 +\eta^\prime \|\nabla\xi_\omega\|^2 + \frac{1}{\tau}\|\xi_m\|^2
+\sigma\|\div\xi_m\|^2 + \sigma\|\xi_k\|^2 \\+ \mu_0\sigma\|\div\xi_H\|^2 + (\eta^\prime+\lambda^\prime) \|\div\xi_\omega\|^2  
 + \frac{\chi_0+(1+\chi_0)\mu_0}{\tau}\|\xi_H\|^2 + \zeta\|\curl\xi_u - 2\xi_\omega\|^2\\
 = J_1 + \mu_0J_2 - J_3 + J_4 + J_5 + J_6.  
\end{align*}
\end{lemma}
\begin{proof}
Taking $\bs v_h = \xi_u\in \bs S_h$ in \eqref{eq:err-semi-1}, $\bs s_h = \xi_\omega \in \bs\Sigma_h$ in \eqref{eq:err-semi-2}, and $\bs F_h = \xi_m \in \bs V_h$ in \eqref{eq:err-semi-3}, we obtain
\begin{equation}
\label{eq:ID-1}
\begin{split}
\frac{\rho}{2}\frac{\dd}{\dd t}\|\xi_u\|^2 + \eta\|\nabla\xi_u\|^2 + \zeta\|\curl\xi_u\|^2 - 2\zeta(\xi_\omega,\curl\xi_u) = -\rho f_{11}(\xi_u) + \mu_0f_{12}(\xi_u) + \frac{\mu_0}{2}f_{13}(\xi_u) \\
+ \frac{\mu_0}{2}f_{14}(\xi_u) + e_1(\xi_u) -(\div\theta_u,\xi_p),
\end{split}
\end{equation}
\begin{equation}\label{eq:ID-2}
\begin{split}
\frac{\rho\kappa}{2}\frac{\dd}{\dd t}\|\xi_\omega\|^2 + \eta^\prime \|\nabla\xi_\omega\|^2 + (\eta^\prime + \lambda^\prime)\|\div\xi_\omega\|^2 - 2\zeta(\curl\xi_u,\xi_\omega)
 + 4\zeta\|\xi_\omega\|^2 = -\rho\kappa f_{21}(\xi_\omega) + \mu_0f_{22}(\xi_\omega) + e_2(\xi_\omega),
\end{split}
\end{equation}
and
\begin{equation}
\label{eq:ID-3}
\begin{split}
\frac{1}{2}\frac{\dd}{\dd t}\|\xi_m\|^2 +\sigma\|\div\xi_m\|^2 + \sigma(\curl\xi_k,\xi_m) -\frac{1}{2}(\curl\xi_z,\xi_m) + \frac{1}{\tau}\|\xi_m\|^2 
 - \frac{\chi_0}{\tau}(\xi_H,\xi_m) 
\\= f_{31}(\xi_m) + \frac{1}{2}f_{32}(\xi_m) 
+ \frac{1}{2}f_{33}(\xi_m) + f_{34}(\xi_m) + e_3(\xi_m).
\end{split}
\end{equation}
Equations \eqref{eq:FHD-weak-7} and \eqref{eq:FHD-weak-dis-7} imply
\begin{equation}\label{eq:ID-4}
(\xi_H,\bs G_h) + (\xi_\phi,\div\bs G_h) = (\theta_H,\bs G_h)\qquad\forall~\bs G_h \in \bs V_h,
\end{equation}
and equations \eqref{eq:FHD-weak-8} and \eqref{eq:FHD-weak-dis-8} imply
\begin{equation}
\label{eq:ID-5}
(\div\xi_H,r_h) + (\div\xi_m,r_h) = 0\qquad\forall~r_h \in W_h.
\end{equation}
Taking $\bs G_h = \xi_H$ and $\bs G_h = \xi_m$ in \eqref{eq:ID-4}, and taking $r_h = \xi_\phi$ in \eqref{eq:ID-5}, we have
\begin{align*}
& \|\xi_H\|^2 + (\xi_\phi,\div\xi_H) = (\theta_H,\xi_H),\\
& (\xi_H,\xi_m) + (\xi_\phi,\div\xi_m) = (\theta_H,\xi_m),\\
& (\div\xi_H,\xi_\phi) + (\div\xi_m,\xi_\phi) = 0.
\end{align*}
These equations lead to 
$$
(\xi_m,\xi_H)= -(\xi_\phi,\div\xi_m) + (\theta_H,\xi_m)  = -\|\xi_H\|^2 + (\theta_H,\xi_H) + (\theta_H,\xi_m).
$$
Substituting the above equation into \eqref{eq:ID-3},  
adding the resulting equation with \eqref{eq:ID-1} and \eqref{eq:ID-2}, and using the fact that $(\curl\xi_k,\xi_m) = \|\xi_k\|^2 - (\theta_k,\xi_k)$, we have
\begin{equation}\label{eq:ID-6}
\begin{split}
\frac{1}{2}\frac{\dd }{\dd t}\left(\rho\|\xi_u\|^2 + \rho\kappa\|\xi_\omega\|^2 + \|\xi_m\|^2 \right) + \eta\|\nabla\xi_u\|^2 + \eta^\prime \|\nabla\xi_\omega\|^2 + (\eta^\prime + \lambda^\prime) \|\div\xi_\omega\|^2 
+\sigma\|\div\xi_m\|^2 + \sigma\|\xi_k\|^2 \\
+ \frac{1}{\tau} \|\xi_m\|^2 + \zeta \|\curl\xi_u - 2\xi_\omega\|^2 + \frac{\chi_0}{\tau} \|\xi_H\|^2 \\
 = J_1 + \frac{1}{2}(\curl\xi_z,\xi_m) - \rho f_{11}(\xi_u) + \mu_0f_{12}(\xi_u) + \frac{\mu_0}{2}f_{13}(\xi_u) 
 + \frac{\mu_0}{2} f_{14}(\xi_u) - \rho\kappa f_{21}(\xi_\omega) + \mu_0f_{22}(\xi_\omega) 
+ f_{31}(\xi_m) \\
+ \frac{1}{2}f_{32}(\xi_m) +\frac{1}{2}f_{33}(\xi_m) + f_{34}(\xi_m).
\end{split}
\end{equation}

Since $\curl_h\xi_H = \curl_hI_h^d\bs H - \curl_h\bs H_h = Q_h^d\curl\bs H - \curl_h\bs H_h = 0$, there exists $\psi_h \in W_h$ such that $\xi_H = \grad_h\psi_h$.
Differentiating \eqref{eq:ID-5} with respect to $t$, taking $r_h = \psi_h$ in the resulting equation and using integration by part, we get
\begin{equation}\label{eq:ID-007}
\frac{1}{2}\frac{\dd }{\dd t}\|\xi_H\|^2 = -(\partial_t\xi_m,\xi_H).
\end{equation}
Taking $\bs F_h = \xi_H$ in \eqref{eq:err-semi-3} and using the fact that $\curl_h\xi_H = 0$, we get
\begin{equation*}
\begin{split}
(\partial_t\xi_m,\xi_H) + \sigma(\div\xi_m,\div\xi_H) + \frac{1}{\tau}(\xi_m,\xi_H) - \frac{\chi_0}{\tau}\|\xi_H\|^2 = f_{31}(\xi_H) + \frac{1}{2}f_{32}(\xi_H) + \frac{1}{2}f_{33}(\xi_H) + f_{34}(\xi_H) + e_3(\xi_H).
\end{split}
\end{equation*}
Equations \eqref{eq:FHD-weak-8} and \eqref{eq:FHD-weak-dis-8} imply 
$$
\div\xi_m = - \div\xi_H,
$$
which plus \eqref{eq:ID-007} yields
\begin{equation}\label{eq:ID-7}
\begin{split}
\frac{1}{2}\frac{\dd}{\dd t}\|\xi_H\|^2  + \sigma\|\div\xi_H\|^2 + \frac{1+\chi_0}{\tau}\|\xi_H\|^2 
=  -f_{31}(\xi_H) - \frac{1}{2}f_{32}(\xi_H)  
- \frac{1}{2}f_{33}(\xi_H) 
- f_{34}(\xi_H) + J_2.
\end{split}
\end{equation}
Multiplying \eqref{eq:ID-7} with $\mu_0$ and adding the resulting equation with \eqref{eq:ID-6}, we finally obtain the desired result.
\end{proof}

 The following Lemmas show  the estimates of $J_i$ ($i =1,2,\dots,6$). Using the regularity Theorems \ref{the:reg-2} - \ref{the:reg-aux} and the properties of the interpolation (or quasi-interpolation) operators, we easily have the estimates of $J_1$ and $J_2$. 
\begin{lemma}\label{lem:J1J2}
Under  Assumptions \ref{ass:initial-data} - \ref{ass:omega}, there hold
\begin{align*}
J_1  \lesssim & h^2\left[ \|\bs u_t\|_2\|\xi_u\|+\|\bs\omega_t\|_2\|\xi_\omega\|\right]  + h^2 \|\bs\omega\|_2(\|\nabla\xi_u\| + \|\xi_\omega\|)  +h(\|\bs m\|_1+\|\bs m_t\|_1 +\|\bs H\|_1)\|\xi_m\| + h\|\bs H\|_1\|\xi_H\|\\
&  + h(\|\bs u\|_2 +\|\tilde p\|_1)\|\nabla\xi_u\|  + h\|\bs\omega\|_2\|\xi_\omega\| +h\|\bs\omega\|_2\|\nabla\xi_\omega\| + h(\|\bs k\|_1 +\|\bs z\|_1 )\|\xi_k\| + h\|\bs u\|_2\|\xi_p\|,\\
 J_2 \lesssim & h(\|\bs m_t\|_1 + \|\bs m\|_1 + \|\bs H\|_1 )\|\xi_H\|+h\|\bs H\|_1\|\xi_m\|).
\end{align*}
\end{lemma}
 
 Lemmas \ref{lem:J3} - \ref{lem:J6} give  the estimates of $J_i$ ($i = 3,4,5,6$).
 \begin{lemma}
 \label{lem:J3}
 Under   Assumptions \ref{ass:initial-data} - \ref{ass:omega}, there holds
 \begin{align*}
 J_3 \lesssim & \|\bs u\|_2\|\xi_u\|\|\nabla\xi_u\|  + \|\bs\omega\|_2( \|\xi_u\|\|\nabla\xi_\omega\| +\|\nabla\xi_u\|\|\xi_\omega\| ) \\
 & + h^2 \|\bs u\|_2 \|\bs\omega\|_2\|\nabla\xi_\omega\| + h^2\|\bs u\|_2^2 \|\nabla\xi_u\| + h\|\bs u\|_2^2\|\xi_u\| + h\|\bs u\|_2 \|\bs \omega\|_2 \|\xi_\omega\|.
 \end{align*}
 \end{lemma}
 \begin{proof}
 Some simple calculations show
 \begin{align*}
 f_{11}(\xi_u) 
 & = b(\xi_u,\pi_h^S\bs u,\xi_u) - b(\bs u,\theta_u,\xi_u) - b(\theta_u,\pi_h^S\bs u,\xi_u) \\
 &\lesssim \|\bs u\|_{2} \|\xi_u\|\|\nabla\xi_u\| + \|\bs u\|_{2}\|\xi_u\|\|\nabla\theta_u\|  + \|\bs u\|_{2}\|\nabla\xi_u\|\|\theta_u\| ,\\
f_{21}(\xi_\omega) & = -b(\bs u,\theta_\omega,\xi_\omega) - b(\theta_u,\pi_h^\Sigma\bs\omega,\xi_\omega) + b(\xi_u,\pi_h^\Sigma,\xi_\omega) \\
&\lesssim \|\bs u\|_{2}\|\xi_\omega\|\|\nabla\theta_\omega\| + \|\bs u\|_{2}\|\nabla\xi_\omega\|\|\theta_\omega\| + \|\bs\omega\|_{2} \|\xi_\omega\|\|\nabla\theta_u\| \\
& \quad + \|\bs\omega\|_{2}\|\nabla\xi_\omega\|\|\theta_u\| + \|\bs\omega\|_{2}\|\nabla\xi_u\|\|\xi_\omega\| + \|\bs\omega\|_{2} \|\xi_u\|\|\nabla\xi_\omega\|.
\end{align*}
Then the  desired result follows from the properties of the interpolation operators $\pi_h^S$ and $\pi_h^\Sigma$.
 \end{proof}

\begin{lemma}\label{lem:J4}
Under   Assumptions \ref{ass:initial-data} - \ref{ass:omega}, there holds 
\begin{align*}
J_4 & \lesssim h\|\bs m\|_2\|\bs H\|_2 \|\xi_u\|  + h\|\bs H\|_2\|\bs m\|_2 \|\nabla\xi_u\|  + h\|\bs u\|_2\|\bs m\|_1\|\div\xi_H\|+ h(\|\bs u\|_2\|+\|\bs\omega\|_2)\|\bs m\|_2\|\xi_H\|\\
& \quad   + h\|\bs H\|_2\|\bs m\|_2\|\xi_\omega\| + h^2\|\bs m\|_2\|\bs\omega\|_2\|\xi_H\|  + \|\bs H\|_2\|\nabla\xi_u\|\|\xi_m\| +\|\bs H\|_1\|\xi_m\|\|\nabla\xi_\omega\| \\
&\quad +\|\bs H\|_2\|\xi_u\|^{1/2}\|\nabla\xi_u\|^{1/2}\|\xi_m\| + \|\bs H\|_2\|\xi_u\|\|\div\xi_m\| + \|\bs u\|_2\|\xi_m\|\|\div\xi_H\| +\|\bs u\|_2 \|\xi_H\|\|\div\xi_m\| \\
&\quad  + \|\bs H\|_1\|\xi_m\|\|\nabla\xi_\omega\|+\|\bs\omega\|_2\|\xi_m\|\|\xi_H\|.
\end{align*}
\end{lemma}
\begin{proof}
Note that
\begin{align*}
& f_{12}(\xi_u) - f_{31}(\xi_H)  = -c(\xi_u,\theta_m,\bs H) + c(\xi_u,\xi_m,\bs H) + c(\xi_u,\xi_m,\theta_H) - c(\xi_u, I_h^d\bs m,\theta_H) \\
&\quad + c(\theta_u,\bs m,\xi_H) - c(\pi_h^S\bs u,\xi_m,\xi_H) + c(\pi_h^S\bs u,\theta_m,\xi_H)\\
 \lesssim & \|\bs H\|_{2}\|\theta_m\| \|\xi_u\|^{1/2}\|\nabla\xi_u\|^{1/2} + \|\bs H\|_{2} \|\xi_u\|\|\div\theta_m\|  + \|\bs H\|_{2}\|\xi_u\|^{1/2}\|\nabla\xi_u\|^{1/2}\|\xi_m\|  \\
 & + \|\bs H\|_{2}\|\xi_u\|\|\div\xi_m\|  + \|\nabla\xi_u\|\|\xi_m\|^{1/2}\|\xi_m\|_{1,h}^{1/2}\|\div\theta_H\|  + \|\nabla\xi_u\|\|\div\xi_m\|\|\theta_H\|_{L^3} \\
&+ \|\bs m\|_{2}\|\xi_u\|\|\div\theta_H\| + \|\bs m\|_{2} \|\xi_u\|^{1/2}\|\xi_u\|^{1/2}\|\theta_H\| + \|\bs m\|_{2} \|\div\xi_H\|\|\theta_u\| \\
& + \|\bs m\|_{2} \|\xi_H\|^{1/2}\|\xi_H\|_n^{1/2}\|\theta_u\|   + \|\bs u\|_{2}\|\xi_m\|\|\div\xi_H\|  + \|\bs u\|_{2}\|\xi_H\|\|\div\xi_m\|\\
&+ \|\bs u\|_{2}\|\div\xi_H\|\|\theta_m\| +\|\bs u\|_{2}\|\xi_H\|\|\div\theta_m\|
\end{align*}
and
\begin{align*}
f_{22}(\xi_\omega) - f_{34}(\xi_H) = & -(\theta_m\times \bs H,\xi_\omega) + (\xi_m\times\bs H,\xi_\omega) + (\xi_m\times \theta_H,\xi_\omega) - (I_h^d\bs m\times \theta_H,\xi_\omega) \\
&+ (\theta_\omega \times \bs m,\xi_H) - (\pi_h^\Sigma\bs\omega \times \xi_m,\xi_H) + (\pi_h^\Sigma\bs\omega \times\theta_m,\xi_H) \\
\lesssim & \|\bs H\|_{2}\|\xi_\omega\|\|\theta_m\| + \|\bs H\|_{1}\|\xi_m\|\|\nabla\xi_\omega\| + \|\xi_m\|^{1/2}\|\xi_m\|_n^{1/2}\|\nabla\xi_\omega\|\|\theta_H\| \\
& + \|\bs m\|_{2}(\|\xi_\omega\|\|\theta_H\| +\|\xi_H\| \|\theta_\omega\|) + \|\bs\omega\|_{2}(\|\xi_m\|\|\xi_H\|  + \|\xi_H\|\|\theta_m\|).
\end{align*}
Using the properties of the operators $\pi_h^\Sigma$ and $I_h^d$ and the inverse inequalities, we get the desired result.
\end{proof}

\begin{lemma}\label{lem:J5}
Under   Assumptions \ref{ass:initial-data} - \ref{ass:omega}, there holds
\begin{align*}
J_5 & \lesssim h\|\bs k\|_1\|\bs H\|_2 \|\xi_u\| + h(\|\bs k\|_1 + \|\bs m\|_2)\|\bs u\|_2 \|\xi_H\| + h\|\bs H\|_2(\|\bs m\|_2 + \|\bs k\|_1)\|\nabla\xi_u\| \\
&\quad + h\|\bs m\|_1\|\bs u\|_2\|\div\xi_H\| + h\|\bs u\|_2\|\xi_m\|\|\div\xi_H\| + \|\bs H\|_2\|\xi_u\|\|\xi_k\| + \|\bs H\|_2\|\xi_m\| \|\nabla\xi_u\| \\
&\quad + \|\bs u\|_2\|\xi_k\|\|\xi_H\| + \|\xi_m\|\|\nabla\xi_u\| + \|\bs u\|_2\|\xi_m\|\|\div\xi_H\|.
\end{align*}
\end{lemma}
\begin{proof}
Since 
\begin{align*}
f_{13}(\xi_u) - f_{33}(\xi_H) = & -(\xi_u\times \theta_k,\bs H) + (\xi_u\times \xi_k,\bs H) + (\xi_u\times \xi_k,\theta_H) - (\xi_u\times I_h^c\bs k,\theta_H) \\
& + (\theta_u\times \bs k,\xi_H) - (\pi_h^S\bs u\times \xi_k,\xi_H) + (\pi_h^S\bs u\times \theta_k,\xi_H) \\
\lesssim & \|\bs H\|_{2} \|\xi_u\|\|\theta_k\| + \|\bs H\|_{2}\|\xi_u\|\|\xi_k\| + \|\nabla\xi_u\|\|\xi_k\|\|\theta_H\|_{L^3}\\
& + \|I_h^c\bs k\|_{L^3}\|\nabla\xi_u\|\|\theta_H\| + \|\bs k\|_{1}\|\xi_H\|\|\nabla\theta_u\| + \|\bs u\|_{2} \|\xi_k\|\|\xi_H\| \\
& + \|\bs u\|_{2}\|\xi_H\|\|\theta_k\|
\end{align*}
and
\begin{align*}
f_{14}(\xi_u) - f_{32}(\xi_H) = & -(\theta_m\times\curl\xi_u,\bs H) + (\xi_m\times\curl\xi_u,\bs H) + (\xi_m\times \curl\xi_u,\theta_H)\\
& -(I_h^d\bs m\times\curl\xi_u,\theta_H) + (\theta_m\times \curl\bs u,\xi_H) - (\xi_m\times \curl\bs u,\xi_H) \\
& -(\xi_m\times\curl\theta_u,\xi_H) + (I_h^d\bs m\times\curl\theta_u,\xi_H) \\
\lesssim & \|\bs H\|_2\|\nabla\xi_u\|\|\theta_m\| + \|\bs H\|_2 \|\xi_m\|\|\nabla\xi_u\| + \|\xi_m\|_{n}\|\nabla\xi_u\|\|\theta_H\|_{L^3} \\
& + \|\bs m\|_{2}\|\nabla\xi_u\|\|\theta_H\| +\|\bs u\|_{2} \|\div\xi_H\|\|\theta_m\| + \|\bs u\|_{2}\|\xi_m\|\|\div\xi_H\| \\
& + \|\xi_m\|\|\div\xi_H\| \|\nabla\theta_u\|  + \|\bs m\|_{2}\|\xi_H\|\|\nabla\theta_u\|,
\end{align*}
using the properties of the operators $\pi_h^S$, $I_h^c$ and $I_h^d$ and the inverse inequalities, we get the desired result.
\end{proof}
\begin{lemma}\label{lem:J6}
Under   Assumptions \ref{ass:initial-data} - \ref{ass:omega}, there holds
\begin{align*}
J_6 & \lesssim h(\|\bs u\|_2\|\bs m\|_2+\|\bs z\|_1)\|\xi_k\| + h[\|\bs u\|_2(\|\bs m\|_2+\|\bs k\|_1)+ \|\bs m\|_2(\|\bs\omega\|_2 + \|\bs m\|_1)]\|\xi_m\| \\
& \quad + h\|\bs k\|_1\|\xi_z\| + h\|\bs u\|_2\|\bs m\|_1\|\xi_m\|_n  + \|\bs m\|_2\|\xi_u\|\|\div\xi_m\|\\
&\quad + (\|\bs k\|_1+\|\bs m\|_2)\|\nabla\xi_u\|\|\xi_m\| + \|\bs m\|_1\|\nabla\xi_\omega\|\|\xi_m\|.
\end{align*}
\end{lemma}
\begin{proof}
Equations \eqref{eq:FHD-weak-6} and \eqref{eq:FHD-weak-dis-6} imply that
$$
(\xi_k,\bs \Theta_h) = (\xi_m,\curl\bs\Theta_h) + (\theta_k,\bs\Theta_h)\qquad\forall~\bs\Theta_h \in \bs U_h.
$$
Taking $\bs\Theta_h = \xi_z \in \bs U_h$, we have
$$
(\xi_k,\xi_z) = (\xi_m,\curl\xi_z) + (\theta_k,\xi_z).
$$
Equations \eqref{eq:FHD-weak-5} and \eqref{eq:FHD-weak-dis-5} indicate
$$
(\xi_z,\bs\Lambda_h)  - (\bs u\times\bs m - \bs u_h \times\bs m_h,\bs\Lambda_h) = (\theta_z,\bs\Lambda_h)\qquad\forall~\bs\Lambda_h \in \bs U_h.
$$
Taking $\bs \Lambda_h = \xi_k$, we get
$$
(\xi_z,\xi_k) = (\bs u\times\bs m - \bs u_h \times\bs m_h,\xi_k) + (\theta_z,\xi_k).
$$
Thus, we obtain
\begin{align*}
& (\xi_m,\curl\xi_z) + f_{33}(\xi_m) =  -(\theta_u\times\bs m,\xi_k) + (\xi_u\times\bs m,\xi_k) - (\pi_h^S\bs u\times\theta_m,\xi_k)  - (\theta_u\times\bs k,\xi_m) \\
& + (\xi_u\times \bs k,\xi_m) + (\xi_u\times\theta_k,\xi_m)  - (\pi_h^S\bs u\times \theta_k,\xi_m) + (\theta_z,\xi_k) - (\theta_k,\xi_z) \\
\lesssim & \|\bs m\|_{2}\|\xi_k\|\|\theta_u\| + \|\bs m\|_{2}\|\xi_u\|\|\xi_k\| + \|\bs u\|_{2}\|\xi_k\|\|\theta_m\|  + \|\bs k\|_{1}\|\xi_m\|\|\nabla\theta_u\| + \|\bs k\|_{1}\|\nabla\xi_u\|\|\xi_m\| \\
&+ \|\nabla\xi_u\|\|\xi_m\|_n \|\theta_k\|  + \|\bs u\|_{2}\|\xi_m\|\|\theta_k\| + \|\xi_k\|\|\theta_z\| + \|\xi_z\|\|\theta_k\|.
\end{align*}
Some simple calculations show that
\begin{align*}
f_{31}(\xi_m) = & -c(\theta_u,\bs m,\xi_m) + c(\xi_u,\bs m,\xi_m) + c(\xi_u,\theta_m,\xi_m) - c(\pi_h^S\bs u,\theta_m,\xi_m) \\
\lesssim & \|\bs m\|_{2} \|\div\xi_m\| \|\theta_u\|+ \|\bs m\|_{2}\|\nabla\theta_u\|\|\xi_m\| + \|\bs m\|_{2}\|\nabla\xi_u\|\|\xi_m\| \\
& + \|\bs m\|_{2}\|\xi_u\|\|\div\xi_m\| + \|\nabla\xi_u\|\|\div\xi_m\|\|\theta_m\|_{L^3} + \|\nabla\xi_u\|\|\xi_m\|_{n}\|\div\theta_m\| \\
& + \|\bs u\|_{2} \|\div\xi_m\| \|\theta_m\| + \|\bs u\|_{2} \|\xi_m\|\|\div\theta_m\|,
\end{align*}
\begin{align*}
f_{32}(\xi_m) = & -(\theta_m\times \curl\bs u,\xi_m) + (I_h^d\bs m\times \curl\xi_u,\xi_m) - (I_h^d\bs m \times \curl\theta_u,\xi_m) \\
\lesssim & \|\bs u\|_{2}\|\xi_m\|_n\|\theta_m\| + \|\bs m\|_{2}\|\nabla\xi_u\|\|\xi_m\| + \|\bs m\|_{1}\|\xi_m\|_n\|\nabla\theta_u\|,
\end{align*}
and
\begin{align*}
f_{34}(\xi_m) = & -(\theta_\omega\times\bs m,\xi_m) + (\xi_\omega\times \bs m,\xi_m) + (\xi_\omega\times \theta_m,\xi_m) - (\pi_h^\Sigma\bs\omega\times\theta_m,\xi_m)\\
\lesssim & \|\bs m\|_{2}\|\xi_m\|\|\theta_\omega\| + \|\bs m\|_{2}\|\xi_\omega\|\|\xi_m\| + \|\nabla\xi_\omega\|\|\xi_m\|_n\|\theta_m\| + \|\bs\omega\|_{2} \|\xi_m\|\|\theta_m\|.
\end{align*}
Then the desired result follows from the inverse inequality and the properties of the operators $\pi_h^S$, $\pi_h^\Sigma$, $I_h^c$ and $I_h^d$.
\end{proof}

Lemmas \ref{lem:xiz} - \ref{lem:xip} give the estimates of $\xi_z$ and $\xi_p$.
\begin{lemma}\label{lem:xiz}
Under   Assumptions \ref{ass:initial-data} - \ref{ass:omega}, there holds
\begin{align*}
\|\xi_z\| \lesssim&  h^2\|\bs m\|_2\|\bs u\|_2  + \|\bs m\|_2\|\xi_u\|+ \|\bs m\|_1\|\nabla\xi_u\|+ h(\|\bs u\|_2\|\bs m\|_1 + \|\bs z\|_1)  + \|\nabla\xi_u\|\|\xi_m\|_{n} +\|\bs u\|_2 \|\xi_m\| .
\end{align*}
\end{lemma}
\begin{proof}
Equations \eqref{eq:FHD-weak-5} and \eqref{eq:FHD-weak-dis-5} imply
$$
(\xi_z,\bs\Lambda_h) = (\bs u \times\bs m - \bs u_h\times\bs m_h,\bs\Lambda_h) + (\theta_z,\bs\Lambda_h)\qquad\forall~\bs\Lambda_h \in \bs U_h.
$$
Taking $\bs\Lambda_h = \xi_z\in\bs U_h$ in the above equation, we have
\begin{align*}
\|\xi_z\|^2 
& = -(\theta_u \times \bs m,\xi_z) + (\xi_u\times \bs m,\xi_z) + (\xi_u\times \theta_m,\xi_z) - (\pi_h^S\bs u\times\theta_m,\xi_z) - (\xi_u\times\xi_m,\xi_z)  + (\pi_h^S\bs u\times\xi_m,\xi_z) + (\theta_z,\xi_z)\\
& \lesssim \|\bs m\|_{2}\|\xi_z\|\|\theta_u\| + \|\bs m\|_{2}\|\xi_u\|\|\xi_z\| + \|\nabla\xi_u\|\|\xi_z\|_{L^3}\|\theta_m\| \\
&\quad +\|\bs u\|_{2}\|\xi_z\|\|\theta_m\| + \|\nabla\xi_u\|\|\xi_m\|_{n}\|\xi_z\| + \|\bs u\|_{2}\|\xi_m\|\|\xi_z\| +\|\theta_z\|\|\xi_z\| \\
& \lesssim h^2\|\bs m\|_2\|\bs u\|_2\|\xi_z\|  + \|\bs m\|_2\|\xi_u\|\|\xi_z\|+ \|\bs m\|_1\|\nabla\xi_u\|\|\xi_z\| + h\|\bs u\|_2\|\bs m\|_1\|\xi_z\| \\
&\quad  + \|\nabla\xi_u\|\|\xi_m\|_{n}\|\xi_z\| +\|\bs u\|_2 \|\xi_m\|\|\xi_z\| + h\|\bs z\|_1\|\xi_z\|.
\end{align*}
The desired result follows by cancelling $\|\xi_z\|$ from the above inequality.
\end{proof}
\begin{lemma}\label{lem:xip}
Under the Assumptions \ref{ass:initial-data} - \ref{ass:omega}, for any $t \in [0,T]$, we have the following estimate for $\xi_p$
\begin{align*}
\int_0^t\|\xi_p&(\cdot,s)\|\dd s  \lesssim  \|\xi_u(\cdot,t)\| + h^2\int_0^t (\|\bs u_t\|_2+\|\bs\omega\|_2)\dd s + h \int_0^t(\|\bs u\|_2 + \|\bs p\|_1)\dd s \\
& + h\int_0^t(\|\bs u\|_2^2 + \|\bs m\|_2\|\bs H\|_2 +\|\bs k\|_1 \|\bs H\|_1)\dd s \\
&  + h \int_0^t (\|\bs u\|_2\|\nabla\xi_u\| +\|\bs H\|_2\|\xi_k\| +\|\bs H\|_2\|\xi_m\|_n)\dd s\\
& + \int_0^t(\|\bs u\|_2\|\nabla\xi_u\| + \|\bs H\|_2\|\xi_m\|_n +\|\bs m\|_2\|\div\xi_H\| + \|\bs H\|_1\|\xi_k\| + \|\bs k\|_1\|\xi_H\| )\dd s\\
&   +\int_0^t(\|\xi_u\|^{1/2}\|\nabla\xi_u\|^{3/2} + \|\xi_m\|^{1/2}\|\xi_m\|_{n}^{1/2} \|\div\xi_H\| +  \|\xi_k\|\|\xi_H\|^{1/2}\|\div\xi_H\|^{1/2} )\dd s.
\end{align*}
\end{lemma}
\begin{proof}
The inf-sup condition \eqref{eq:inf-sup} and equation\eqref{eq:ID-1} imply that for any $0\leq t \leq T$,
\begin{align*}
\int_0^t\|\xi_p&(\cdot,s)\|\dd s  \lesssim \sup\limits_{\bs v_h \in \bs S_h} \int_0^t \frac{(\xi_p(\cdot,s),\nabla\cdot\bs v_h)}{\|\bs v_h\|_1}\dd s \\
& \lesssim \|\xi_u(\cdot,t)\| + \int_0^t(\|\nabla\xi_u(\cdot,s)\| + \|\xi_\omega(\cdot,s)\|)\dd s + h^2\int_{0}^t(\|\partial_t\bs u(\cdot,s)\|_2 + \|\bs\omega(\cdot,s)\|_2)\dd s \\
&\quad + h\int_0^t(\|\bs u(\cdot,s)\|_2 + \|p(\cdot,s)\|_1)\dd s + T_1
\end{align*}
with 
$$
T_1: = \sup\limits_{\bs v_h \in \bs S_h}\int_0^t \frac{|f_{11}(\bs v_h) + f_{12}(\bs v_h) + f_{13}(\bs v_h) + f_{14}(\bs v_h)|}{\|\bs v_h\|_1} \dd s,
$$
where $f_{11}(\bs v_h),~f_{12}(\bs v_h),~f_{13}(\bs v_h),~f_{14}(\bs v_h)$ are as same as those in \eqref{eq:ID-1}.
For any $\bs v_h \in \bs S_h$, using the properties of the (quasi-) interpolation operators and the inverse inequalities, we have 
\begin{align*}
|f_{11}(\bs v_h) |& \leq |b(\theta_u,\bs u,\bs v_h)| + |b(\xi_u,\bs u,\bs v_h)| + |b(\xi_u,\theta_u,\bs v_h)|  + |b(\pi_h^S\bs u,\theta_u,\bs v_h)| + |b(\xi_u,\xi_u,\bs v_h)| + |b(\pi_h^S\bs u,\xi_u,\bs v_h)| \\
& \lesssim \|\nabla\theta_u\|\|\nabla\bs u\|\|\bs v_h\|_{1} + \|\theta_u\|\|\bs u\|_{2}\|\bs v_h\|_1 + \|\xi_u\|\|\nabla\bs v_h\|\|\bs u\|_2 + \|\xi_u\|\|\nabla\bs v_h\|\|\bs u\|_{2} \\
&\quad+ \|\nabla\xi_u\|\|\nabla\theta_u\|\|\nabla\bs v_h\| +\|\nabla\xi_u\|\|\nabla\bs v_h\| \|\nabla\theta_u\| + \|\bs u\|_{2}\|\nabla\theta_u\|\|\bs v_h\| + \|\bs u\|_{2}\|\nabla\bs v_h\| \|\theta_u\|\\
&\quad + \|\xi_u\|_{L^3} \|\nabla\xi_u\|\|\nabla\bs v_h\|  +\|\bs u\|_{2}\|\nabla\xi_u\|\|\bs v_h\| + \|\bs u\|_{2}\|\nabla\bs v_h\|\|\xi_u\| \\
& \lesssim h(\|\bs u\|_2^2 +\|\bs u\|_2\|\nabla\xi_u\|)\|\bs v_h\|_1+ (\|\bs u\|_2\|\nabla\xi_u\| +\|\xi_u\|^{1/2}\|\nabla\xi_u\|^{3/2})\|\bs v_h\|_1 ,
\end{align*}
\begin{align*}
|f_{12}(\bs v_h)| \leq & |c(\bs v_h,\theta_m,\bs H)| + |c(\bs v_h,\xi_m,\bs H)| + |c(\bs v_h,\xi_m,\xi_H)| + |c(\bs v_h,I_h^d\bs m,\xi_H)|  + |c(\bs v_h,\xi_m,\theta_H)| +|c(\bs v_h,I_h^d\bs m,\theta_H)| \\
\lesssim & \|\nabla\bs v_h\|(\|\theta_m\|\|\bs H\|_2 + \|\bs H\|_{1}\|\div\theta_m\| +\|\bs H\|_{1}\|\div\xi_m\|) \\
& + \|\nabla\bs v_h\|(\|\xi_m\|\|\bs H\|_2 + \|\xi_m\|_{L^3}\|\div\xi_H\| + \|\xi_H\|_{L^3}\|\div\xi_m\|)\\
& + \|\nabla\bs v_h\|(\|\bs m\|_{2}\|\div\xi_H\| + \|\xi_H\|\|\bs m\|_2  + \|\xi_m\|_{L^3}\|\div\theta_H\|) \\
& + \|\nabla\bs v_h\| ( \|\theta_H\|_{L^3}\|\div\xi_m\| + \|\bs m\|_{2} \|\div\theta_H\| + \|\theta_H\|\|\bs m\|_2)\\
\lesssim & h\|\bs m\|_2\|\bs H\|_2 \|\nabla\bs v_h\| + \|\bs H\|_2\|\xi_m\|\|\bs v_h\|_1 +\|\bs H\|_1\|\div\xi_m\|\|\bs v_h\|_1\\
& + \|\xi_m\|^{1/2}\|\xi_m\|_n^{1/2}\|\div\xi_H\|\|\bs v_h\|_1 + \|\bs m\|_2\|\div\xi_H\|\|\bs v_h\|_1,
\end{align*}
\begin{align*}
|f_{13}(\bs v_h)| \leq & |(\bs v_h\times\theta_k,\bs H)| + |(\bs v_h\times\xi_k,\bs H)| + |(\bs v_h\times\xi_k,\xi_H)| + |(\bs v_h\times I_h^c\bs k,\xi_H)| + |(\bs v_h\times\xi_k,\theta_H)| + |(\bs v_h\times I_h^c\bs k,\theta_H)| \\
\leq & \| \bs v_h\|_{1} (\|\theta_k\|\|\bs H\|_{1} + \|\xi_k\|\|\bs H\|_{1} + \|\xi_k\|\|\xi_H\|_{L^3} + \|\bs k\|_{1}\|\xi_H\|)  + \|\bs v_h\|_{1} (\|\xi_k\|\|\theta_H\|_{L^3} + \|\bs k\|_{1}\|\theta_H\|)\\
\lesssim & h(\|\bs k\|_1\|\bs H\|_1 + \|\bs H\|_2\|\xi_k\|) \|\bs v_h\|_1 + (\|\bs H\|_1\|\xi_k\| + \|\bs k\|_1\|\xi_H\|)\|\bs v_h\|_1 + \|\xi_k\|\|\xi_H\|^{1/2}\|\xi_H\|_n^{1/2},
\end{align*}
and
\begin{align*}
|f_{14}(\bs v_h)| \leq & |(\theta_m\times\curl\bs v_h,\bs H)| + |(\xi_m\times\curl\bs v_h,\bs H)| + |(\xi_m\times\curl\bs v_h,\xi_H)|x + |(I_h^d\bs m\times\curl\bs v_h,\xi_H)| \\
& + |(\xi_m\times\curl\bs v_h,\theta_H)| + |(I_h^d\bs m\times\curl\bs v_h,\xi_H)|\\
\leq &\|\bs v_h\|_1 ( \|\bs H\|_{2}\|\theta_m\| + \|\bs H\|_{2}\|\xi_m\| + \|\xi_m\|_{L^3}\|\div\xi_H\|) \\
& + \|\bs v_h\|_1 (\|\bs m\|_{2}\|\xi_H\| + \|\xi_m\|_{n}\|\theta_H\|_{L^3} + \|\bs m\|_{2}\|\xi_H\|)\\
\lesssim & h(\|\bs m\|_1\|\|\bs H\|_2 + \|\bs H\|_2\|\xi_m\|_n)\|\bs v_h\|_1 + (\|\bs H\|_2\|\xi_m\| +\|\bs m\|_2\|\xi_H\|)\|\bs v_h\|_1 \\
& +\|\xi_m\|^{1/2}\|\xi_m\|_n^{1/2}\|\div\xi_H\|\|\bs v_h\|_1.
\end{align*}
As a result, the desired result follows.
\end{proof}

In what follows let us estimate the terms  $\xi_u, \xi_\omega, \xi_m$ and $ \xi_H$. To this end, 
for any $t \in [0,T]$ we define 
$$
\mathfrak E_h(t): = \rho\|\xi_u(\cdot,t)\|^2 + \rho\kappa\|\xi_\omega(\cdot,t)\|^2 + \|\xi_m(\cdot,t)\|^2 + \mu_0\|\xi_H(\cdot,t)\|^2
$$
and
\begin{align*}
\mathfrak F_h(t) : = &\eta\|\nabla\xi_u(\cdot,t)\|^2  + \eta^\prime\|\nabla\xi_\omega(\cdot,t)\|^2 + \frac{1}{\tau}\|\xi_m(\cdot,t)\|^2 + (\eta^\prime + \lambda^\prime)\|\div\xi_\omega(\cdot,t)\|^2 +\sigma\|\div\xi_m\|^2 + \sigma\|\xi_k\|^2 \\
& + \sigma\mu_0\|\div\xi_H\|^2 + \frac{\chi_0+(1+\chi_0)\mu_0}{\tau}\|\xi_H(\cdot,t)\|^2  +\zeta\|\curl\xi_u(\cdot,t) - 2\xi_\omega(\cdot,t)\|^2.
\end{align*}
We have the following conclusion.
\begin{lemma}\label{lem:xi-estimate}
Under   Assumptions \ref{ass:initial-data} - \ref{ass:omega}, there holds
$$
\sup\limits_{0\leq t\leq T}\left(\mathfrak E_h(t) + \int_0^t\mathfrak F_h(s)\dd s \right) \lesssim h^2.
$$
\end{lemma}
\begin{proof}
By Lemmas \ref{lem:err-ID} - \ref{lem:J6}, we have
\begin{equation}\label{eq:err-pr-1}
\frac{1}{2}\frac{\dd}{\dd t}\mathfrak E_h(t)  + \mathfrak F_h(t) \leq Ch^2 \mathfrak T_1(t) + Ch\mathfrak T_2(t) + C\mathfrak T_3(t),
\end{equation}
where
\begin{align*}
\mathfrak T_1(t) & = \|\bs u_t\|_2\|\xi_u\| + (\|\bs\omega_t\|_2+\|\bs\omega\|_2)\|\xi_\omega\| + (\|\bs u\|_2^2 + \|\bs\omega\|_2)\|\nabla\xi_u\| + \|\bs u\|_2\|\bs\omega\|_2\|\nabla\xi_\omega\|  + \|\bs m\|_2\|\bs\omega\|_2\|\xi_H\| ,\\
\mathfrak T_2(t) & =(\|\bs u\|_2^2 +\|\bs m\|_2\|\bs H\|_2 +\|\bs k\|_1\|\bs H\|_2) \|\xi_u\| + \|\bs u\|_2\|\xi_p\| +\|\bs u\|_2\|\bs m\|_1 \|\xi_m\|_n\\
&\quad + (\|\bs u\|_2+\|\tilde p\|_1+\|\bs H\|_2(\|\bs m\|_2+\|\bs k\|_1))\|\nabla\xi_u\| \\
&\quad + (\|\bs\omega\|_2+\|\bs u\|_2\|\bs\omega\|_2+\|\bs H\|_2\|\bs m\|_2)\|\xi_\omega\| +\|\bs\omega\|_2 \|\nabla\xi_\omega\| \\
&\quad +(\|\bs m\|_1+\|\bs m_t\|_1 +\|\bs H\|_1 +\|\bs u\|_2(\|\bs m\|_2 + \|\bs k\|_1)  + \|\bs m\|_2(\|\bs\omega\|_2+\|\bs m\|_1)) \|\xi_m\| \\
&\quad  + (\|\bs m_t\|_1 + \|\bs m\|_1 + \|\bs H\|_1 +(\|\bs u\|_2+\|\bs\omega\|_2)\|\bs m\|_2 +\|\bs k\|_1\|\bs u\|_2) \|\xi_H\|\\
& \quad   + (\|\bs k\|_1+\|\bs z\|_1 +\|\bs u\|_2\|\bs m\|_2) \|\xi_k\| + \|\bs k\|_1\|\xi_z\|+ \|\bs u\|_2\|\bs m\|_1 \|\div\xi_H\|,\\
\mathfrak T_3(t) & = \|\bs u\|_2 \|\xi_u\|\|\nabla\xi_u\| + \|\bs\omega\|_2 \|\xi_u\|\|\nabla\xi_\omega\| + \|\bs\omega\|_2\|\nabla\xi_u\|\|\xi_\omega\| \\
&\quad +(\|\bs H\|_1+\|\bs m\|_1)\|\xi_m\|\|\nabla\xi_\omega\| + \|\bs H\|_2\|\xi_u\|\|\xi_k\| + \|\bs u\|_2\|\xi_k\|\|\xi_H\|\\
& \quad + (\|\bs H\|_2+\|\bs k\|_1 +\|\bs m\|_2)\|\nabla\xi_u\|\|\xi_m\| + (\|\bs H\|_2+\|\bs m\|_2)\|\xi_u\|\|\xi_m\|_n\\
&\quad + (\|\bs u\|_2 +\|\bs\omega\|_2)\|\xi_m\|\|\div\xi_H\| +\|\bs u\|_2 \|\xi_H\|\|\div\xi_m\| .
\end{align*}
For any $t \in [0,T]$, using Lemma \ref{lem:xip} and the inverse inequality, we get
\begin{align*}
h\int_0^t&\mathfrak \|\xi_p\|\dd s  \lesssim h\|\xi_u(\cdot,t)\| + h^3\int_0^t(\|\bs u_t\|_2 +\|\bs\omega\|_2) \dd s  + h^2\int_0^t (\|\bs u\|_2 + \|\tilde p\|_1) \dd s \\
& + h^2\int_0^t (\|\bs u\|_2 + \|\bs m\|_2\|\bs H\|_2 + \|\bs k\|_1\|\bs H\|_1)\dd s  + h^2 \int_0^t (\|\bs u\|_2\|\nabla\xi_u\| + \|\bs H\|_2\|\xi_k\| +\|\bs H\|_2\|\xi_m\|_n)\dd s \\
& + h\int_0^t(\|\bs u\| \|\nabla\xi_u\| + \|\bs H\|_2\|\xi_m\|_n + \|\bs m\|_2\|\div\xi_H\| +\|\bs H\|_1\|\xi_k\| +\|\bs k\|_1\|\xi_H\| )\dd s \\
&   +h^{1/2} \int_0^t( \|\xi_u\|\|\nabla\xi_u\|+ \|\xi_m\| \|\div\xi_H\| + \|\xi_k\| \|\xi_H\| )\dd s.
\end{align*}
Integrating \eqref{eq:err-pr-1} on the interval $(0,t)$, we obtain
\begin{align*}
\frac{1}{2}\mathfrak E_h(t) & + \int_0^t\mathfrak F_h(s)\dd s \leq C h^2\int_0^t\mathfrak T_1(s)\dd s + Ch\int_0^t\mathfrak T_2(s)+ C\int_0^t\mathfrak T_3(s)\dd s.
\end{align*}
Using the inverse inequalities, the Young's inequality and the regularity results Theorems \ref{the:reg-2} - \ref{the:reg-aux}, we have the following inequalities:
\begin{align*}
Ch^2\int_0^t\mathfrak T_1(s)\dd s  \leq &C h^2
  + C \int_0^t(\|\xi_u\|^2 + \|\xi_\omega\|^2 +\|\xi_m\|^2 + \|\xi_H\|^2)\dd s,
\end{align*}
\begin{align*}
Ch \int_0^t \mathfrak T_2(s)\dd s \leq & C h^2 + \frac{\rho}{4}\|\xi_u(\cdot,t)\|^2+ C\int_0^t(\|\xi_u\|^2+\|\xi_\omega\|^2 +\|\xi_m\|^2 + \|\xi_H\|^2)\dd s\\
& + \int_0^t \left( \frac{\eta}{4}\|\nabla\xi_u\|^2 + \frac{\eta^\prime}{4}\|\nabla\xi_\omega\|^2 + \frac{\sigma}{8}\|\xi_m\|_n^2 + \frac{\sigma}{8}\|\xi_k\|^2 + \frac{\sigma\mu_0}{4}\|\div\xi_H\|^2\right)\dd s,
\end{align*}
and
\begin{align*}
C\int_0^t \mathfrak T_3(s)\dd s \leq & \int_0^t \left( \frac{\eta}{4}\|\nabla\xi_u\|^2 + \frac{\eta^\prime}{4}\|\nabla\xi_\omega\|^2 + \frac{\sigma}{8}\|\xi_m\|_n^2 + \frac{\sigma}{8}\|\xi_k\|^2 + \frac{\sigma\mu_0}{4}\|\div\xi_H\|^2\right)\dd s \\
& + C \int_0^t(\|\xi_u\|^2 + \|\xi_\omega\|^2 +\|\xi_m\|^2 + \|\xi_H\|^2)\dd s.
\end{align*}
Thus, we obtain
$$
\mathfrak E_h(t) + \int_0^t\mathfrak F_h(s)\dd s \leq Ch^2 + C\int_0^t(\mathfrak E_h(s)+ \int_0^s \mathfrak F_h(r)\dd r)\dd s,
$$
%
which, together with the Gr\"onwall's inequality, 
 implies the desired result.
\end{proof}

Define the errors 
\begin{align*}
\hbox{err}_1(t) := &\|\bs u(\cdot,t) - \bs u_h(\cdot,t)\|^2 + \|\bs\omega(\cdot,t) - \bs\omega_h(\cdot,t)\|^2 +  \|\bs m(\cdot,t) - \bs m_h(\cdot,t)\|^2  + \|\bs H(\cdot,t) - \bs H_h(\cdot,t)\|^2
\end{align*}
and
\begin{align*}
\hbox{err}_2(t):= & \|\nabla(\bs u -\bs u_h)(\cdot,t)\|^2 +\|\nabla(\bs \omega - \bs\omega_h)(\cdot,t)\|^2 + \|\div(\bs m - \bs m_h)(\cdot,t)\|^2 \\
& + \|\div(\bs H - \bs H_h)(\cdot,t)\|^2 + \|(\bs k - \bs k_h)(\cdot,t)\|^2 +\|(\bs z - \bs z_h)(\cdot,t)\|^2.
\end{align*}
We finish this section by giving the following error estimates for the semi-discrete scheme \eqref{eq:FHD-weak-dis-1} - \eqref{eq:FHD-weak-dis-8}:
\begin{theorem}\label{the:err-semi}
Under   Assumptions \ref{ass:initial-data} - \ref{ass:omega}, there hold
$$
\sup\limits_{0\leq t \leq T} \left(\hbox{err}_1(t) + \int_0^t \hbox{err}_2(s)\dd s \right)\lesssim h^2
$$
and
$$
\int_0^t \|(p - p_h)(\cdot,s)\|\dd s \lesssim h.
$$
\end{theorem}
\begin{proof}
For any $s \in [0,t]$, the triangular inequality implies
\begin{align*}
\hbox{err}_1(s) \lesssim & \|\theta_u(\cdot,s)\| ^2 + \|\theta_\omega(\cdot,s)\|^2 + \|\theta_m(\cdot,s)\|^2 + \|\theta_H(\cdot,s)\|^2 + \mathfrak E_h(s),\\
\hbox{err}_2(s) \lesssim & \|\nabla\theta_u(\cdot,s)\|^2 + \|\theta_\omega(\cdot,s)\|^2 +\|\div\theta_m(\cdot,s)\|^2 + \|\div\theta_H(\cdot,s)\|^2  +\|\theta_k(\cdot,s)\|^2 + \|\theta_z\|^2 + \mathfrak F_h(s),
\end{align*}
and
$$
\|(p-p_h)(\cdot,s)\| \leq \|\theta_p(\cdot,s)\| + \|\xi_p(\cdot,s)\|.
$$
Then Lemmas \ref{lem:xi-estimate}, \ref{lem:xiz} and \ref{lem:xip} and   Theorem \ref{the:reg-2} - \ref{the:reg-aux} yield the desired results.
\end{proof}

 \section{Fully discrete finite element scheme}\label{sec:ful}
 In this section, we give a fully scheme for the ferrofluid model \eqref{eq:FHD-1} - \eqref{eq:initial}, based on the semi-discretization \eqref{eq:FHD-weak-dis-1} - \eqref{eq:FHD-weak-dis-8}. We will show that this full discretization is also energy-stable and has optimal convergence.
 
 \subsection{Full discretization}
For any $1\leq n \leq N$, we denote by $I_n$ the interval $(t_{n-1},t_n)$ and by $v^n$ the numerical approximation of $v(t_n)$ for any quantity $v(t)$ (both scalar or vector). We also set 
$$
\delta_tv^n := \frac{v^n-v^{n-1}}{\dt}.
$$
The fully discrete scheme of the Rosensweig ferrofluid flow model \eqref{eq:FHD-1} - \eqref{eq:initial} reads: Given $\bs u_h^{n-1} \in \bs S_h$, $\bs \omega_h^{n-1} \in \bs \Sigma_h$, $\bs m_h^{n-1} \in \bs V_h$ for $n = 1,2,\dots,N$, find $\bs u_h^n \in \bs S_h$, $\tilde p_h^n\in L_h$, $\bs \omega_h^n \in \bs\Sigma_h$, $\bs m_h^n \in \bs V_h$, $\bs z_h^n \in \bs U_h$, $\bs k_h^n \in \bs U_h$, $\bs H_h^n \in \bs V_h$ and $\phi_h^n\in W_h$ such that
\begin{align}
\label{eq:FHD-weak-disf-1}
\begin{split}
 \rho(\delta_t\bs u_h^n,\bs v_h) + \rho b(\bs u_h^n,\bs u_h^n,\bs v_h) + \eta(\nabla\bs u_h^n,\nabla\bs v_h) + \zeta (\curl\bs u_h^n,\curl\bs v_h)
  - \mu_0 c(\bs v_h,\bs m_h^n,\bs H_h^n)& \\ 
  - (\tilde p_h^n,\div\bs v_h)
  - \frac{\mu_0}{2}(\bs v_h \times \bs k_h^n,\bs H_h^n) 
  - \frac{\mu_0}{2}(\bs m_h^n\times \curl\bs v_h,\bs H_h^n) - 2\zeta(\bs\omega_h^n,\curl\bs v_h) &= 0
\end{split} & \forall ~\bs v_h \in \bs S_h,\\
\label{eq:FHD-weak-disf-2}
 (\div\bs u_h^n,q_h) & = 0 &\forall~q_h\in L_h, \\
\label{eq:FHD-weak-disf-3}
\begin{split}
\rho\kappa(\delta_t\bs\omega_h^n,\bs s_h) + \eta^\prime(\nabla\bs\omega_h^n,\nabla\bs s_h) + (\eta^\prime + \lambda^\prime)(\div\bs\omega_h^n,\div\bs s_h) 
+ \rho\kappa b(\bs u_h^n,\bs\omega_h^n,\bs s_h) &\\-\mu_0(\bs m_h^n\times \bs H_h^n,\bs s_h) - 2\zeta(\curl\bs u_h^n,\bs s_h) 
+ 4\zeta(\bs \omega_h^n,\bs s_h)& = 0
\end{split} & \forall~\bs s_h \in \bs \Sigma_h,\\
\label{eq:FHD-weak-disf-4}
\begin{split}
(\delta_t\bs m_h^n,\bs F_h) - c(\bs u_h^n,\bs m_h^n,\bs F_h) +\sigma(\div\bs m_h^n,\div\bs F_h) + \sigma(\curl\bs k_h^n,\bs F_h) - \frac{1}{2}(\curl\bs z_h^n,\bs F_h)&\\
 -\frac{1}{2}(\bs m_h^n\times\curl\bs u_h^n,\bs F_h) 
- \frac{1}{2}(\bs u_h^n \times \bs k_h^n,\bs F_h) 
 - (\bs\omega_h^n\times \bs m_h^n,\bs F_h) & \\
 + \frac{1}{\tau}(\bs m_h^n,\bs F_h) 
-\frac{\chi_0}{\tau}(\bs H_h^n,\bs F_h) & = 0
\end{split}& \forall~\bs F_h\in \bs V_h,\\
\label{eq:FHD-weak-disf-5} (\bs z_h^n,\bs\Lambda_h) - (\bs u_h^n\times \bs m_h^n,\bs\Lambda_h) & = 0 & \forall~\bs \Lambda_h \in \bs U_h, \\
\label{eq:FHD-weak-disf-6} (\bs k_h^n,\bs \Theta_h) - (\bs m_h^n,\curl\bs\Theta_h) & = 0 & \forall~\bs \Theta_h \in \bs U_h, \\
\label{eq:FHD-weak-disf-7} (\bs H_h^n,\bs G_h) + (\phi_h^n,\div\bs G_h) & = 0 &\forall~\bs G_h \in \bs V_h, \\
\label{eq:FHD-weak-disf-8}\mu_0(\div\bs H_h^n+\div\bs m_h^n,r_h) + (\div\bs H_e^n,r_h) & = 0 &\forall~ r \in W_h.
\end{align}

\subsection{Existence of solutions} We first introduce an auxiliary problem: Given $\mathcal U_h^{n-1} = (\bs u_h^{n-1},\bs\omega_h^{n-1},\bs m_h^{n-1},\bs H_h^{n-1})  \in \mathbb V_H$ ($1\leq n \leq N$), find $\mathcal U_h^n = (\bs u_h^n,\bs \omega_h^n,\bs m_h^n,\bs H_h^n)  \in \mathbb V_h$, such that
\begin{equation}\label{eq:au-f}
(\mathbb D\delta_t\mathcal U_h^n,\mathcal V_h) + \mathbb A(\mathcal U_h^n,\mathcal V_h) + \mathbb B(\mathcal U_h^n,\mathcal V_h) = \mathcal L^n(\mathcal V_h)\qquad\forall~\mathcal V_h \in \mathbb V_h,
\end{equation}
with initial value
$$
\mathcal U_h(\cdot,0) = (\bs u_h(\cdot,0),\bs\omega(\cdot,0),\bs m_h(\cdot,0),\bs H_h(\cdot,0))  \in \mathbb V_h,
$$
where $\bs H_h(\cdot,0)\in \mathfrak K_h^d$ is   determined by the saddle point problem \eqref{eq:H-initial}. Here $\mathcal L^n$ is defined by
$$
\mathcal L^n(\mathcal V_h) :=  -\frac{\chi_0}{\mu_0\tau}(Q_{\mathfrak K_h^d}\bs H_e^n,\bs F_h) -\sigma(\div\bs H_e^n,\div\bs A_h)-(\delta_t\bs H_e^n+\frac{1}{\tau}\bs H_e^n,\bs A_h) .
$$
We have the following Lemma.
\begin{lemma}\label{lem:eq-sol-f}
The fully-discrete form \eqref{eq:FHD-weak-disf-1} - \eqref{eq:FHD-weak-disf-8} and the auxiliary problem \eqref{eq:au-f} are equivalent in the following sense:
\begin{itemize}
\item if $\bs u_h^n\in \bs S_h$, $p_h^n\in L_h$, $\bs\omega_h^n\in \bs \Sigma_h$, $\bs m_h^n\in \bs V_h$, $\bs z_h^n\in \bs U_h$, $\bs k_h^n\in  \bs U_h$, $\bs H_h^n\in\bs V_h$ and $\phi_h^n\in  W_h$ solve \eqref{eq:FHD-weak-disf-1} - \eqref{eq:FHD-weak-disf-8}, then $\mathcal U_h^n =(\bs u_h^n,\bs\omega_h^n,\bs m_h^n,\bs H_h^n)^\intercal \in \mathbb V_h$ solve \eqref{eq:au-f}.
\item if $\mathcal U_h^n =(\bs u_h^n,\bs\omega_h^n,\bs m_h^n,\bs H_h^n)  \in \mathbb V_h$ solve \eqref{eq:au-f}, then $\bs u_h^n$, $\bs\omega_h^n$, $\bs m_h^n$, $\bs z_h = Q_h^c(\bs u_h^n \times \bs m_h^n)$, $\bs k_h^n = \curl_h\bs m_h^n$, a unique $\phi_h^n \in W_h$ satisfying $\bs H_h^n = \grad_h\phi_h$, and a unique $p_h^n \in L_h$ solve \eqref{eq:FHD-weak-disf-1} - \eqref{eq:FHD-weak-disf-8}.  
\end{itemize}
\end{lemma}
\begin{proof}
The proof is almost similar as that of Lemma \ref{lem:eq-sol-semi}, where we only need to change the partial derivative $\partial_t$ to the finite difference $\delta_t$.
\end{proof}
In the following theorem,  we show that the auxiliary problem \eqref{eq:au-f} has solutions.
\begin{theorem}
The fully discrete scheme \eqref{eq:FHD-weak-disf-1} - \eqref{eq:FHD-weak-disf-8} admits at least one solution, 
and there holds 
$$
\|\mathbb D^{1/2}\mathcal U_h^n\|^2 + \dt\interleave \mathcal U_h^n\interleave ^2 \lesssim \|\mathbb D^{1/2}\mathcal U_h^{n-1}\|^2 + \dt (\|\bs H_e^n\|_d^2 + \|\delta_t\bs H_e^n\|^2),
$$
where $\mathcal U_h^n = (\bs u_h^n,\bs\omega_h^n,\bs m_h^n,\bs H_h^n) $.
\end{theorem}
\begin{proof}
By Lemma \ref{lem:eq-sol-f}, we only need to prove the auxiliary problem \eqref{eq:au-f} has solutions. For any $\mathcal U_h^n,~\mathcal V_h \in \mathbb V_h$, we define $\Phi:\mathbb V_h \rightarrow \mathbb V_h$ as
$$
(\Phi(\mathcal U_h^n),\mathcal V_h) := (\mathbb D\mathcal U_h^n,\mathcal V_h) + \dt\mathbb A(\mathcal U_h^n,\mathcal V_h) + \dt\mathbb B(\mathcal U_h^n,\mathcal V_h) -\dt\mathcal L^n(\mathcal V_h) -(\mathbb D\mathcal U_h^{n-1},\mathcal V_h).
$$
We will show $\Phi$ is continuous mapping first. In fact, for any fixed $\mathcal U_{1} = (\bs u_1,\bs\omega_1,\bs m_1,\bs H_1)  \in\mathbb V_h\setminus \{0\}$ and    any $\epsilon >0$, there exists $\delta = \frac{\epsilon}{2C \interleave \mathcal U_1\interleave }>0$ such that for any $\mathcal U_2 = (\bs u_2,\bs \omega_2,\bs m_2,\bs H_2)  \in \mathbb V_h$ satisfying $\interleave \mathcal U_2 - \mathcal U_1\interleave < \delta $, there holds
\begin{align*}
\mathbb B(\mathcal U_2 - \mathcal U_1,\mathcal V_h) & \leq C (\|\nabla\bs u_1\| + \|\nabla\bs\omega_1\| + \|\bs m_1\|_n + \|\div\bs H_1\|) \interleave \mathcal U_2 - \mathcal U_1\interleave \interleave \mathcal V_h \interleave   < \epsilon 
\end{align*}
 for any $\mathcal V_h \in \mathbb V_h$. This means $\mathbb B(\cdot,\cdot)$ is continuous with respect to the first variable on $\mathbb V_h$. Thus $\Phi$ is continuous on the space $\mathbb V_h$.  
The definitions of $\mathbb A(\cdot,\cdot)$ and $\mathbb B(\cdot,\cdot)$ imply
\begin{align*}
(\Phi(\mathcal U_h^n),\mathcal U_h^n) & = \|\mathbb D^{1/2}\mathcal U_h^n\|^2 + \dt\interleave \mathcal U_h^n\interleave^2 - \dt\mathcal L^n(\mathcal U_h^n) -(\mathbb D\mathcal U_h^{n-1},\mathcal U_h^n) \\
& \geq \left(\|\mathbb D^{1/2}\|\mathcal U_h^n\| + \dt^{1/2}\interleave\mathcal U_h^n\interleave -\dt^{1/2}\|\mathcal L^n\|_{\mathbb V_h^\prime} - \|\mathbb D^{1/2}\mathcal U_h^{n-1}\| \right) \left (\|\mathbb D^{1/2} \mathcal U_h^n\| + \dt^{1/2}\interleave \mathcal U_h^n\interleave\right),
\end{align*}
where 
$$
\|\mathcal L^n\|_{\mathbb V_h^\prime} := \sup\limits_{\mathcal V_h \in \mathbb V_h}\frac{\mathcal L^n(\mathcal V_h)}{\interleave \mathcal V_h\interleave} \leq C(\|\bs H_e^n\|_d + \|\delta_t\bs H_e^n\|).
$$
Taking $c = \dt^{1/2}\|\mathcal L^n\|_{\mathbb V_h^\prime} + \|\mathbb D^{1/2}\mathcal U_h^{n-1}\|$, we get
$$
\Phi(\mathcal U_h^n),\mathcal U_h^n) \geq 0\qquad\forall~\mathcal U_h^n \in \mathbb V_h \text{ with } \|\mathbb D^{1/2}\mathcal U_h^n\| + \dt^{1/2}\interleave \mathcal U_h^n\interleave  = c.
$$
So by \citet[Chapter IV, Corollary 1.1]{Girault1986}, there exists an element $\mathcal U_h^n \in \mathbb V_h$ such that
$$
(\Phi(\mathcal U_h^n),\mathcal V_h) = 0\qquad\text{with }\|\mathbb D^{1/2}\mathcal U_h^n\| + \dt^{1/2}\interleave \mathcal U_h^n\interleave  \leq c
$$
This completes the proof.
\end{proof}

Let us  turn to the uniqueness of the fully discrete solution. We make an assumption first.
\begin{assum}\label{ass:f-sol}
There exist constants $h_0>0$ and $M_0>0$ such that for any $0<h\leq h_0$ and any solution $(\bs u_h^n,\bs\omega_h^n,\bs m_h^n,\bs H_h^n,\bs k_h^n,\bs z_h^n,p_h^n,\phi_h^n)$ ($0\leq n \leq N$) of \eqref{eq:FHD-weak-disf-1} - \eqref{eq:FHD-weak-disf-8}, there holds
$$
\dt \sum\limits_{n=1}^N (\|\nabla\bs u_h^n\|^4 + \|\nabla\bs\omega_h^n\|^4 + \|\bs m_h^n\|_n^4 + \|\bs H_h^n\|_n^4) \leq M_0.
$$
\end{assum}

\begin{theorem}\label{the:unique-solu-f-dis}
Under   Assumption \ref{ass:f-sol}, the fully discrete scheme \eqref{eq:FHD-weak-disf-1} - \eqref{eq:FHD-weak-disf-8} admits a unique solution.
\end{theorem}
\begin{proof}
Let $\bs u_i^n \in\bs S_h$, $\tilde p_i^n\in L_h$, $\bs \omega_i^n \in \bs\Sigma_h$, $\bs m_i^n \in \bs V_h$, $\bs z_i^n \in \bs U_h$, $\bs k_i^n \in \bs U_h$, $\bs H_i^n \in \bs V_h$ and $\phi_i^n\in W_h$ ($i = 1,~2$) solve \eqref{eq:FHD-weak-disf-1} - \eqref{eq:FHD-weak-disf-8}. For any quantity $v$ ($v = \bs u,~\tilde p,~\bs\omega,~\bs m,~\bs z,~\bs k,~\bs H,~\phi$),   denote $d_v^n: = v_1^n - v_2^n$. Then we have
\begin{equation}\label{eq:pr-unif-1}
\begin{split}
\rho (\delta_td_u^n,\bs v_h)  + \eta(\nabla d_u^n,\nabla\bs v_h) + \zeta(\curl d_u^n,\curl\bs v_h) - (d_p^n,\div\bs v_h)    - 2\zeta(d_\omega^n,\curl\bs v_h)\\ = -\rho d_{11}^n(\bs v_h) + \mu_0 d_{12}^n(\bs v_h) + \frac{\mu_0}{2} d_{13}^n(\bs v_h)   + \frac{\mu_0}{2}d_{14}^n(\bs v_h),
\end{split}
\end{equation}
\begin{equation}\label{eq:pr-unif-2}
(\div d_u^n,q_h) = 0,
\end{equation}
\begin{equation}\label{eq:pr-unif-3}
\begin{split}
\rho\kappa (\delta_t d_\omega^n,\bs s_h) +\eta^\prime(\nabla d_\omega^n,\nabla\bs s_h) + (\eta^\prime + \lambda^\prime)(\div d_\omega^n,\div\bs s_h) - 2\zeta(\curl d_u^n,\bs s_h) 
+ 4\zeta(d_\omega^n,\bs s_h) = \mu_0d_{31}^n(\bs s_h)  
-\rho\kappa d_{32}^n(\bs s_h),
\end{split}
\end{equation}
\begin{equation}\label{eq:pr-unif-4}
\begin{split}
(\delta_t d_m^n,\bs F_h)  + \sigma(\div d_m^n,\div\bs F_h) + \sigma(\curl d_k^n,\bs F_h) - \frac{1}{2}(\curl d_z^n,\bs F_h) 
 + \frac{1}{\tau}(d_m^n,\bs F_h) - \frac{\chi_0}{\tau}(d_H^n,\bs F_h)\\
  = d_{41}^n(\bs F_h)   + \frac{1}{2}d_{42}^n(\bs F_h)  
 + \frac{1}{2} d_{43}^n(\bs F_h)   + d_{44}^n(\bs F_h) ,
\end{split}
\end{equation}
\begin{equation}\label{eq:pr-unif-5}
(d_z^n,\bs\Lambda_h) = d_{51}^n(\Lambda_h) ,
\end{equation}
\begin{equation}\label{eq:pr-unif-6}
(d_k^n,\bs\Theta_h) - (d_m^n,\curl\bs\Theta_h) = 0,
\end{equation}
\begin{equation}\label{eq:pr-unif-7}
(d_H^n,\bs G_h) +(d_\phi^n,\div\bs G_h) = 0,
\end{equation}
\begin{equation}\label{eq:pr-unif-8}
(\div d_H^n + \div d_m^n,r_h) = 0,
\end{equation}
  for all $\bs v_h\in \bs S_h$, $q_h \in L_h$, $\bs s_h \in \bs \Sigma_h$, $\bs F_h \in \bs V_h$, $\bs\Lambda_h\in \bs U_h$, $\bs\Theta_h \in \bs U_h$, $\bs G_h \in \bs V_h$ and $r_h \in W_h$,
with initial data
$$
d_u^0 = \bs 0,\qquad d_\omega^0 = \bs 0,\qquad d_m^0 = \bs 0.
$$
Here
\begin{align*}
& d_{11}^n(\bs v_h) := b(\bs u_{1}^n,\bs u_{1}^n,\bs v_h) - b(\bs u_{2}^n,\bs u_{2}^n,\bs v_h),\qquad\qquad
  d_{12}^n(\bs v_h): = c(\bs v_h,\bs m_{1}^n,\bs H_{1}^n) - c(\bs v_h,\bs m_{2}^n,\bs H_{2}^n), \\
&d_{13}^n(\bs v_h): = (\bs v_h\times \bs k_{1}^n,\bs H_{1}^n) - (\bs v_h\times \bs k_{2}^n,\bs H_{2}^n),\qquad\,\,\,
 d_{14}^n(\bs v_h) := (\bs m_{1}^n\times\curl\bs v_h,\bs H_{1}^n) - (\bs m_{2}^n\times\curl\bs v_h,\bs H_{2}^n),\\
& d_{31}^n(\bs s_h) :=  (\bs m_{1}^n\times \bs H_{1}^n,\bs s_h) -(\bs m_{2}^n\times \bs H_{2}^n,\bs s_h),\qquad\,\, d_{32}^n(\bs s_h): = b(\bs u_{1}^n,\bs\omega_{1}^n,\bs s_h) - b(\bs u_{2}^n,\bs\omega_{2}^n,\bs s_h),\\
& d_{41}^n(\bs F_h) := c(\bs u_{1}^n,\bs m_{1}^n,\bs F_h) - c(\bs u_{2}^n,\bs m_{2}^n,\bs F_h),\qquad\,\,\,\, d_{42}^n(\bs F_h) := (\bs m_{1}^n \times\curl\bs u_{1}^n,\bs F_h) - (\bs m_{2}^n \times\curl\bs u_{2}^n,\bs F_h), \\
&d_{43}^n(\bs F_h) := (\bs u_{1}^n\times\bs k_{1}^n,\bs F_h) - (\bs u_{2}^n\times\bs k_{2}^n,\bs F_h) ,\qquad\,\,\,\, d_{44}^n(\bs F_h): = (\bs \omega_{1}^n\times \bs m_{1}^n,\bs F_h) - (\bs\omega_{2}^n \times \bs m_{2}^n,\bs F_h),\\
& d_{51}^n(\bs\Lambda_h) := (\bs u_{1}^n\times \bs m_{1}^n,\bs\Lambda_h) - (\bs u_{2}^n\times \bs m_{2}^n,\bs\Lambda_h).
\end{align*}
Taking $\bs v_h = d_u^n$ in \eqref{eq:pr-unif-1}, $\bs s_h = d_\omega^n$ in \eqref{eq:pr-unif-3}, and $\bs F_h = d_m^n$ and $d_H^n$ in \eqref{eq:pr-unif-4}, and using   similar techniques as in the proof of Theorem \ref{the:unique-solu-semi-dis}, we get
\begin{align*}
\rho(\delta_td_u^n,d_u^n) + \rho\kappa(\delta_td_\omega^n,d_\omega^n) + (\delta_td_m^n,d_m^n) + \mu_0(\delta_td_H^n,d_H^n ) + \eta|d_u^n|_1^2 + \eta^\prime|d_\omega^n|_1^2 
+ \sigma\|\div d_m^n\|^2 
+ \sigma\|d_k^n\|^2\\ 
+ \sigma\mu_0 \|\div d_H^n\|^2  + (\eta^\prime + \lambda^\prime)\|\div d_\omega^n\|^2 
 + \frac{1}{\tau}\|d_m^n\|^2 
+ \frac{\chi_0 + (1+\chi_0)\mu_0}{\tau} \|d_H^n\|^2 +\zeta\|\curl d_u^n - 2d_\omega^n\|^2 \\
= \frac{1}{2}(\curl d_z^n,d_m^n) -\rho d_{11}^n(d_u^n) 
 + \mu_0d_{12}^n(d_u^n)
+ \frac{\mu_0}{2} d_{13}^n(d_u^n) 
 + \frac{\mu_0}{2}d_{14}^n(d_u^n) 
 + \mu_0 d_{31}^n(d_\omega^n) 
 -\rho\kappa d_{32}^n(d_\omega^n) 
 + d_{41}^n(d_m^n) \\
 + \frac{1}{2}d_{42}^n(d_m^n) 
 + \frac{1}{2}d_{43}^n(d_m^n) 
 + d_{44}^n(d_m^n) 
 - \mu_0d_{41}^n(d_H^n) -\frac{\mu_0}{2}d_{42}(d_H^n) - \frac{\mu_0}{2} d_{43}^n(d_H^n) - \mu_0d_{44}^n(d_H^n).
\end{align*}
Using   similar estimates of $d_{ij}$ as in the proof of Theorem \ref{the:unique-solu-semi-dis} and the inequalities
\begin{align*}
(d_u^n,d_u^{n-1}) & \leq \frac{1}{2} (\|d_u^n\|^2 + \|d_u^{n-1}\|^2), \qquad \qquad\qquad
(d_\omega^n,d_\omega^{n-1})  \leq \frac{1}{2} (\|d_\omega^n\|^2 + \|d_\omega^{n-1}\|^2), \\
(d_m^n,d_m^{n-1}) & \leq \frac{1}{2} (\|d_m^n\|^2 + \|d_m^{n-1}\|^2), \qquad\qquad\qquad
(d_H^n,d_H^{n-1})  \leq \frac{1}{2} (\|d_H^n\|^2 + \|d_H^{n-1}\|^2), 
\end{align*}
we obtain
\begin{align*}
& (\rho\|d_u^n\|^2 + \rho\kappa\|d_\omega^n\|^2 + \|d_m^n\|^2 + \mu_0\|d_H^n\|^2 ) - (\rho\|d_u^{n-1}\|^2 + \rho\kappa\|d_\omega^{n-1}\|^2   + \|d_m^{n-1}\|^2 + \mu_0\|d_H^{n-1}\|^2 ) 
 + \dt \eta|d_u^n|_1^2 \\
&\quad + \dt \eta^\prime|d_\omega^n|_1^2 + \dt \sigma\|\div d_m^n\|^2  + \dt \sigma\|d_k^n\|^2 + \dt \sigma\mu_0 \|\div d_H^n\|^2  + \dt (\eta^\prime + \lambda^\prime)\|\div d_\omega^n\|^2 + \frac{\dt}{\tau}\|d_m^n\|^2  \\
 &\qquad
+ \dt\frac{\chi_0 + (1+\chi_0)\mu_0}{\tau} \|d_H^n\|^2  +\dt \zeta\|\curl d_u^n - 2d_\omega^n\|^2 \\
\leq & C\dt (\|\nabla\bs u_{2}^n\|^4 +\|\nabla\bs u_{1}^n\|^4 + \|\bs H_{1}^n\|_n^4  +\|\bs H_{2}^n\|_n^4  
+ \|\bs \omega_{1}^n\|^4 + \|\bs m_{1}^n\|_n^4 +\|\bs m_{2}^n\|_n^4 ) 
(\|d_u^n\|^2 + \|d_\omega^n\|^2 + \|d_m^n\|^2 + \|d_H^n\|^2).
\end{align*}
For any $1\leq L \leq N$, summing over $n$ from $1$ to $L$, we have
$$
(\rho\|d_u^L\|^2 + \rho\kappa\|d_\omega^L\|^2 + \|d_m^L\|^2 + \mu_0\|d_H^L\|^2 ) \leq C\dt\sum\limits_{n=1}^L (\|d_u^n\|^2 + \|d_\omega^n\|^2 + \|d_m^n\|^2 + \|d_H^n\|^2).
$$
Then the discrete Gr\"owall's inequality implies
$$
d_u^n = 0,\quad d_\omega^n = 0,\quad d_m^n = 0,\quad d_H^n = 0.
$$
For any given $\bs u_h^n$, $\bs m_h^n$ and $\bs H_h^n$, equations \eqref{eq:FHD-weak-disf-5} - \eqref{eq:FHD-weak-disf-8} give unique $\bs z_h^n$, $\bs k_h^n$ and $\phi_h^n$, which means that the solution of the full discrete scheme \eqref{eq:FHD-weak-disf-1} - \eqref{eq:FHD-weak-disf-8} is unique.
\end{proof}

\subsection{Energy estimates}
For the fully discrete scheme\eqref{eq:FHD-weak-disf-1} - \eqref{eq:FHD-weak-disf-8}, define the energy
$$
\mathcal E_h^n := \rho\|\bs u_h^n\|^2 + \rho\kappa\|\bs\omega_h^n\|^2 + \|\bs m_h^n\|^2 + \mu_0\|\bs H_h^n\|^2
$$ 
and the dissipated energy
\begin{align*}
\mathcal F_h^n = & \eta\|\nabla\bs u_h^n\|^2 + \eta^\prime \|\nabla\bs\omega_h^n\|^2 + (\eta^\prime+\lambda^\prime)\|\div\bs\omega_h^n\|^2 + \frac{1}{\tau}\|\bs m_h^n\|^2 + \frac{1+(1+\mu_0)\chi_0}{\tau}\|\bs H_h^n\|^2\\
 &  + \sigma\|\div\bs m_h^n\|^2 + \sigma \|\bs k_h^n\|^2  + \mu_0\sigma\|\div\bs H_h^n\|^2 +\zeta\|\curl\bs u_h^n - 2\bs\omega_h^n\|^2.
\end{align*}
Then we have the following energy estimate.
\begin{theorem}\label{the:ener-dis-f}
For any given $\bs u_h^{n-1} \in \bs S_h$, $\bs \omega_h^{n-1} \in \bs \Sigma_h$, $\bs m_h^{n-1} \in \bs V_h$, let $\bs u_h^n \in \bs S_h$, $\tilde p_h^n\in L_h$, $\bs \omega_h^n \in \bs\Sigma_h$, $\bs m_h^n \in \bs V_h$, $\bs z_h^n \in \bs U_h$, $\bs k_h^n \in \bs U_h$, $\bs H_h^n \in \bs V_h$ and $\phi_h^n\in W_h$ ($1\leq n \leq N$) solve the fully discretization scheme \eqref{eq:FHD-weak-disf-1} - \eqref{eq:FHD-weak-disf-8}, then  
\begin{align*}
\mathcal E_h^n + 2\dt \mathcal F_h^n \leq \mathcal E_h^{n-1} + C\dt \left(\|\bs H_e^n\|_{d} + \|\delta_t\bs H_e^n\| \right).
\end{align*}
\end{theorem}
\begin{proof}
Taking $\bs v_h = \bs u_h^n$ 
in \eqref{eq:FHD-weak-disf-1} and $q = p_h^n$ in  \eqref{eq:FHD-weak-disf-2}, we get
\begin{equation}\label{eq:Ener-dis-p-1}
\begin{split}
\rho(\delta_t \bs u_h^n,\bs u_h^n)  + \eta\|\nabla\bs u_h^n\|^2 + \zeta\|\curl\bs u_h^n\|^2 - 2\zeta(\curl\bs u_h^n,\bs \omega_h^n)  = \mu_0 c(\bs u_h^n,\bs m_h^n,\bs H_h^n) \\
+ \frac{\mu_0}{2}(\bs u_h^n\times \bs k_h^n,\bs H_h^n) + \frac{\mu_0}{2}(\bs m_h^n\times\curl\bs u_h^n,\bs H_h^n).
\end{split}
\end{equation}
Taking $\bs s_h = \bs \omega_h^n \in \bs\Sigma_h$ in \eqref{eq:FHD-weak-disf-3}, we have
\begin{equation}
\label{eq:Ener-dis-p-2}
\begin{split}
\rho\kappa(\delta_t\bs\omega_h^n,\bs\omega_h^n) + \eta^\prime\|\nabla\omega_h^n\|^2 + (\eta^\prime + \lambda^\prime)\|\div\bs\omega_h^n\|^2 - 2\zeta(\curl\bs u_h^n,\bs \omega_h^n)
 + 4\zeta\|\bs\omega_h^n\|^2 = \mu_0(\bs m_h^n\times \bs H_h^n,\bs\omega_h^n).
\end{split}
\end{equation}
Taking $\bs F_h = \bs m_h^n\in \bs V_h^n$ in \eqref{eq:FHD-weak-disf-4} and using the fact that $(\bs z_h^n,\bs k_h^n) = (\bs u_h^n\times \bs m_h^n,\bs k_h^n) = (\bs m_h^n,\curl\bs z_h^n)$, we obtain
\begin{equation}
\label{eq:Ener-dis-p-3}
(\delta_t\bs m_h^n,\bs m_h^n) +\sigma\|\div\bs m_h^n\|^2 + (\curl\bs k_h^n,\bs m_h^n) + \frac{1}{\tau}\|\bs m_h^n\|^2 - \frac{\chi_0}{\tau} (\bs H_h^n,\bs m_h^n) = 0.
\end{equation}
Taking $\bs G_h = \bs H_h^n$ and $\bs G_h = \bs m_h^n$ in \eqref{eq:FHD-weak-disf-7},   $r_h = \phi_h^n$ in \eqref{eq:FHD-weak-disf-8} and $\bs G_h = Q_h^v \bs H_e$ in \eqref{eq:FHD-weak-disf-7}, we get
\begin{align*}
& \|\bs H_h^n\|^2 = -(\phi_h^n,\div\bs H_h^n),\qquad\qquad (\bs H_h^n,\bs m_h^n) = -(\phi_h^n,\div\bs m_h^n),\\
& (\div\bs H_h^n\phi_h^n) + (\div\bs m_h^n,\phi_h^n) = -\frac{1}{\mu_0}(\div\bs H_e,\phi_h^n) = \frac{1}{\mu_0}(\bs H_h^n,\bs H_e^n),
\end{align*}
which imply
$$
(\bs H_h^n,\bs m_h^n) = -(\phi_h^n,\div\bs m_h^n) = -\|\bs H_h^n\|^2 - \frac{1}{\mu_0}(\bs H_h^n,\bs H_e^n).
$$
This relation, together with \eqref{eq:Ener-dis-p-3} and   the fact that $(\curl\bs k_h^n,\bs m_h^n) = \|\bs k_h^n\|^2$,  yields
\begin{equation}
\label{eq:Ener-dis-p-4}
(\delta_t\bs m_h^n,\bs m_h^n) +\sigma\|\div\bs m_h^n\|^2 + \sigma\|\bs k_h^n\|^2 + \frac{1}{\tau}\|\bs m_h^n\|^2 + \frac{\chi_0}{\tau}\|\bs H_h^n\|^2 = -\frac{\chi_0}{\tau \mu_0}(\bs H_h^n,\bs H_e^n).
\end{equation}
Adding equations \eqref{eq:Ener-dis-p-1}, \eqref{eq:Ener-dis-p-2} and \eqref{eq:Ener-dis-p-4} together, we get
\begin{equation}
\label{eq:Ener-dis-p-5}
\begin{split}
&\rho(\delta_t  \bs u_h^n,\bs u_h^n) + \rho\kappa(\delta_t\bs\omega_h^n,\bs\omega_h^n) + (\delta_t\bs m_h^n,\bs m_h^n) + \eta\|\nabla\bs u_h^n\|^2 + \eta^\prime\|\nabla\bs\omega_h^n\|^2 + (\eta^\prime+\lambda^\prime)\|\div\bs\omega_h^n\|^2 + \sigma\|\div\bs m_h^n\|^2\\
&\qquad  + \sigma\|\bs k_h^n\|^2 + \frac{1}{\tau}\|\bs m_h^n\|^2  + \frac{\chi_0}{\tau}\|\bs H_h^n\|^2  + \zeta\|\curl\bs u_h^n - 2\bs\omega_h^n\|^2 \\
= & \mu_0 c(\bs u_h^n,\bs m_h^n,\bs H_h^n) + \frac{\mu_0}{2}(\bs u_h^n\times \bs k_h^n,\bs H_h^n) + \frac{\mu_0}{2}(\bs m_h^n\times\curl\bs u_h^n,\bs H_h^n)  + \mu_0(\bs m_h^n\times \bs H_h^n,\bs\omega_h^n) - \frac{\chi_0}{\tau\mu_0}(\bs H_h^n,\bs H_e^n).
\end{split}
\end{equation}
Taking $r_h = \phi_h^n$ in \eqref{eq:FHD-weak-disf-8} at time levels $t_n$ and $t_{n-1}$, we get
$$
(\div\bs H_h^n,\phi_h^n) + (\div\bs m_h^n,\phi_h^n) = -\frac{1}{\mu_0}(\div\bs H_e^n,\phi_h^n)
$$
and
$$
(\div\bs H_h^{n-1},\phi_h^n) + (\div\bs m_h^{n-1},\phi_h^n) = -\frac{1}{\mu_0}(\div\bs H_e^{n-1},\phi_h^n).
$$
These two equations give
\begin{equation}\label{eq:Ener-dis-p-6}
(\div\delta_t\bs H_h^n,\phi_h^n) + (\div\delta_t\bs m_h^n,\phi_h^n) = -\frac{1}{\mu_0}(\div\delta_t\bs H_e^n,\phi_h^n).
\end{equation}
Taking $\bs G_h = \delta_t\bs H_h^n$, $\bs G_h = \delta_t\bs m_h^n$ and $\bs G = \pi_h^v\delta_t\bs H_e^n$ in \eqref{eq:FHD-weak-disf-7}, respectively, we have
\begin{align*}
& (\div\delta_t\bs H_h^n,\phi_h^n) = -(\bs H_h^n,\delta_t\bs H_h^n),\quad (\div\delta_t\bs m_h^n,\phi_h^n) = -(\bs H_h^n,\delta_t\bs m_h^n), \\
& (\div\delta_t\bs H_e^n,\phi_h^n) = -(\bs H_h^n,\pi_h^v\delta_t\bs H_e^n).
\end{align*}
Taking $\bs F_h = \bs H_h^n$ in \eqref{eq:FHD-weak-disf-4}, we get
\begin{equation}
\label{eq:Ener-dis-p-7}
\begin{split}
(\delta_t\bs m_h^n,\bs H_h^n) - c(\bs u_h^n,\bs m_h^n,\bs H_h^n) + \sigma(\div\bs m_h^n,\div\bs H_h^n) - \frac{1}{2}(\bs m_h^n\times\curl\bs u_h^n,\bs H_h^n)- \frac{1}{2}(\bs u_h^n\times\bs k_h^n,\bs H_h^n) & \\
 - \frac{1}{2}(\curl\bs z_h^n,\bs H_h^n) - (\bs \omega_h^n\times \bs m_h^n,\bs H_h^n) 
+ \frac{1}{\tau}(\bs m_h^n,\bs H_h^n) - \frac{\chi_0}{\tau}\|\bs H_h^n\|^2  & = 0.
\end{split}
\end{equation}
Thus, we have
\begin{equation}
\label{eq:Ener-dis-p-8}
\begin{split}
(\delta_t\bs H_h^n,\bs H_h^n) & = (\div\delta_t\bs m_h^n,\phi_h^n) + \frac{1}{\mu_0} (\div\delta_t\bs H_e^n,\phi_h^n)  = -(\delta_t\bs m_h^n,\bs H_h^n) - \frac{1}{\mu_0}(\pi_h^v\delta_t\bs H_e^n,\bs H_h^n) \\
& = -\frac{1}{\mu_0}(\pi_h^v\delta_t\bs H_e^n,\bs H_h^n) - c(\bs u_h^n,\bs m_h^n,\bs H_h^n) + \sigma(\div\bs m_h^n,\div\bs H_h^n) - (\bs \omega_h^n\times\bs m_h^n,\bs H_h^n)\\
& \quad  - \frac{1}{2}(\bs m_h^n\times\curl\bs u_h^n,\bs H_h^n) - \frac{1}{2}(\bs u_h^n\times \bs k_h^n,\bs H_h^n)-\frac{1+\chi_0}{\tau}\|\bs H_h^n\|^2 - \frac{1}{\tau\mu_0}(\bs H_e^n,\bs H_h^n).
\end{split}
\end{equation}
Multiplying \eqref{eq:Ener-dis-p-8} with $\mu_0$, adding the resulting equation with \eqref{eq:Ener-dis-p-5}, and using the fact that $(\div\bs m_h^n,\div\bs H_h^n) = -\|\div\bs H_h^n\|^2 - \frac{1}{\mu_0}(\div\bs H_e^n,\div\bs H_h^n)$, we obtain
\begin{align*}
& \rho(\delta_t\bs u_h^n,\bs u_h^n) + \rho\kappa (\delta_t\bs\omega_h^n,\bs\omega_h^n) + (\delta_t\bs m_h^n,\bs m_h^n) + \mu_0(\delta_t\bs H_h^n,\bs H_h^n) + \mathcal F_h^n \\
= & -\frac{\chi_0+1}{\tau\mu_0}(\bs H_e^n,\bs H_h^n) - (\pi_h^v\delta_t\bs H_e^n,\bs H_h^n)-\sigma(\div\bs H_e^n,\div\bs H_h^n)\\
\leq & \frac{\chi_0+1}{\tau\mu_0}\|\bs H_e^n\|\|\bs H_h^n\| + C \|\delta_t\bs H_e^n\|\|\bs H_h^n\| + \sigma\|\div\bs H_e^n\|\|\div\bs H_h^n\|.
\end{align*}
Using the Cauchy-Schwarz inequality and the Young's inequality, we have
$$
\mathcal E_h^n + 2\dt\mathcal F_h^n \leq \mathcal E_h^{n-1} + C\dt\left( \|\bs H_e^n\|_{H(\div)} +\|\delta_t\bs H_e^n\| \right)\|\bs H_h^n\|_{H(\div)} .
$$
Finally, the desired result follows from Lemma \ref{lem:Gr}.
\end{proof}

\subsection{Error estimates}
To carry out error analysis for the fully discrete scheme \eqref{eq:FHD-weak-disf-1} - \eqref{eq:FHD-weak-disf-8}, we first introduce the following notations:
 $$
 \begin{array}{lll}
  \xi_u^n: = \pi_h^S\bs u(\cdot,t_n) - \bs u_h^n,& \quad \xi_\omega^n := \pi_h^\Sigma \bs\omega(\cdot,t_n) - \bs\omega_h^n, &\quad\xi_m^n := I_h^d\bs m(\cdot,t_n) - \bs m_h^n,\\
 \xi_H^n := I_h^d\bs H(\cdot,t_n) - \bs H_h^n,
& \quad \xi_p^n := Q_h^L \tilde p(\cdot,t_n) - \tilde p_h^n, &\quad \xi_z^n := I_h^c\bs z(\cdot,t_n) - \bs z_h^n,\\
 \xi_k^n := I_h^c\bs k(\cdot,t_n) - \bs k_h^n, &\quad\xi_\phi^n := Q_h^W\phi(\cdot,t_n) - \phi_h^n, &\quad \theta_u^n := \pi_h^S\bs u(\cdot,t_n) - \bs u(\cdot,t_n),\\
  \theta_\omega^n := \pi_h^\Sigma\bs\omega(\cdot,t_n) - \bs\omega(\cdot,t_n), & \quad \theta_m^n := I_h^d\bs m(\cdot,t_n) - \bs m(\cdot,t_n), & \quad \theta_H^n := I_h^d\bs H(\cdot,t_n) - \bs H(\cdot,t_n), \\
 \theta_p^n := Q_h^L\tilde p(\cdot,t_n) - \tilde p(\cdot,t_n), & \quad\theta_z^n := I_h^c\bs z(\cdot,t_n) - \bs z(\cdot,t_n), &\quad \theta_k^n := I_h^c\bs k(\cdot,t_n) - \bs k(\cdot,t_n) \\
 \theta_\phi^n := Q_h^W\phi(\cdot,t_n) - \phi(\cdot,t_n). & & 
 \end{array}
 $$
 From \eqref{eq:FHD-weak-1}, \eqref{eq:FHD-weak-3} , \eqref{eq:FHD-weak-4},  \eqref{eq:FHD-weak-disf-1}, \eqref{eq:FHD-weak-disf-3} and \eqref{eq:FHD-weak-disf-4}  we get the following error equations:
\begin{equation}\label{eq:err-full-1}
\begin{split}
\rho(\delta_t\xi_u^n,\bs v_h) + \eta(\nabla\xi_u^n,\nabla\bs v_h) + \zeta (\curl\xi_u^n,\curl\bs v_h) - (\xi_p^n,\div\bs v_h) - 2\zeta(\xi_\omega^n,\curl \bs v_h)\\
 = -\rho f_{11}^n(\bs v_h) + \mu_0 f_{12}^n(\bs v_h) + \frac{\mu_0}{2} f_{13}^n(\bs v_h) 
+ \frac{\mu_0}{2} f_{14}^n(\bs v_h) +\tilde e_1^n(\bs v_h) 
\end{split}\qquad \forall~\bs v_h \in \bs S_h,
\end{equation}
\begin{equation}\label{eq:err-full-2}
\begin{split}
\rho\kappa(\delta_t\xi_\omega^n,\bs s_h) + \eta^\prime (\nabla\xi_\omega^n,\nabla\bs s_h) + (\eta^\prime + \lambda^\prime) (\div\xi_\omega^n,\div\bs s_h) + 4\zeta(\xi_\omega^n,\bs s_h)- 2\zeta (\curl \xi_u^n,\bs s_h) \\
 = -\rho\kappa f_{21}^n(\bs s_h) + \mu_0 f_{22}^n(\bs s_h) + \tilde e_2^n(s_h)
\end{split}\qquad\forall~\bs s_h \in \bs\Sigma_h,
\end{equation}
\begin{equation}\label{eq:err-full-3}
\begin{split}
(\delta_t\xi_m^n,\bs F_h) +\sigma(\div\xi_m^n,\div\bs F_h) + \sigma(\curl\xi_k^n,\bs F_h)- \frac{1}{2}(\curl\xi_z^n,\bs F_h)+ \frac{1}{\tau}(\xi_m^n,\bs F_h)- \frac{\chi_0}{\tau}(\xi_H^n,\bs F_h)\\   
= f_{31}^n(\bs F_h) + \frac{1}{2} f_{32}^n(\bs F_h) + \frac{1}{2} f_{33}^n(\bs F_h) 
+ f_{34}^n(\bs F_h) + \tilde e_3^n(\bs F_h)
\end{split}\quad\forall\bs F_h \in \bs V_h,
\end{equation}
with
\begin{align*}
&\tilde e_1^n(\bs v_h) := \rho(\delta_t\theta_u^n,\bs v_h) + \rho(\delta_t\bs u^n - \partial_t\bs u^n,\bs v_h) + \eta(\nabla\theta_u^n,\nabla\bs v_h) + \zeta(\curl\theta_u,\curl\bs v_h)  - (\theta_p^n,\div\bs v_h) - 2\zeta(\theta_\omega^n,\curl\bs v_h),\\
&\tilde e_2^n(\bs s_h) := \rho\kappa(\delta_t\theta_\omega^n,\bs s_h)+\rho\kappa(\delta_t\bs\omega^n-\partial_t\bs\omega^n,\bs s_h) + \eta^\prime (\nabla\theta_\omega^n,\nabla\bs v_h)  + (\eta^\prime + \lambda^\prime) (\div\theta_\omega^n,\div\bs s_h) \\
&\qquad\qquad- 2\zeta(\curl\theta_\omega^n,\bs s_h) + 4\zeta(\theta_\omega^n,\bs s_h),\\
&\tilde e_3^n(\bs F_h) := (\delta_t\theta_m^n,\bs F_h) + (\delta_t\bs m^n - \partial_t\bs m^n,\bs F_h)  + \sigma(\curl\theta_k^n,\bs F_h) + \frac{1}{\tau}(\theta_m^n,\bs F_h) - \frac{\chi_0}{\tau}(\theta_H^n,\bs F_h)-\frac{1}{2}(\curl\theta_z^n,\bs F_h).
\end{align*}
We also define
\begin{align*}
&\tilde J_1^n := \tilde e_1^n(\xi_u^n) + \tilde e_2^n(\xi_\omega^n) + \tilde e_3^n(\xi_m^n) - (\div\theta_u^n,\xi_p^n) + \frac{\chi_0}{\tau}(\theta_H^n,\xi_H^n) +\frac{\chi_0}{\tau} (\theta_H^n,\xi_m^n) + \sigma(\theta_k^n,\xi_k^n), \\
&\tilde J_2^n := -\tilde e_3^n(\xi_H^n) + \frac{1}{\tau}(\theta_H^n,\xi_H^n) + \frac{1}{\tau}(\theta_H^n,\xi_m^n).
\end{align*}
Similar to Lemma \ref{lem:err-ID}, the following conclusion holds:
\begin{lemma}\label{lem:err-ID-f}
There holds
\begin{align*}
\rho(\delta_t\xi_u^n,\xi_u^n) + \rho\kappa(\delta_t\xi_\omega^n,\xi_\omega^n) + (\delta_t\xi_m^n,\xi_m^n) + \mu_0(\delta_t\xi_H^n,\xi_H^n) + \eta\|\nabla\xi_u^n\|^2 
+\eta^\prime \|\nabla\xi_\omega^n\|^2 
+ \sigma\|\div\xi_m^n\|^2 + \sigma\|\xi_k^n\|^2\\ 
+ \mu_0\sigma\|\div\xi_H^n\|^2+ \frac{1}{\tau}\|\xi_m^n\|^2
+ (\eta^\prime+\lambda^\prime) \|\div\xi_\omega^n\|^2  
+ \frac{\chi_0+(1+\chi_0)\mu_0}{\tau}\|\xi_H^n\|^2 + \zeta\|\curl\xi_u^n - 2\xi_\omega^n\|^2\\
 = \tilde J_1^n + \mu_0\tilde J_2^n - J_3^n + J_4^n + J_5^n + J_6^n.  
\end{align*}
\end{lemma}
\begin{proof}
Taking $\bs v_h = \xi_u^n$ in \eqref{eq:err-full-1}, $\bs s_h = \xi_\omega^n\in \bs\Sigma_h$ in \eqref{eq:err-full-2}, and $\bs F_h = \xi_m^n \in \bs V_h$ in \eqref{eq:err-full-3}, we obtain
\begin{equation}
\label{eq:ID-f1}
\begin{split}
\rho(\delta_t\xi_u^n,\xi_u^n) + \eta\|\nabla\xi_u^n\|^2 + \zeta\|\curl\xi_u^n\|^2 - 2\zeta(\xi_\omega^n,\curl\xi_u^n) = -\rho f_{11}^n(\xi_u^n) + \mu_0f_{12}^n(\xi_u^n) + \frac{\mu_0}{2}f_{13}^n(\xi_u^n) \\
+ \frac{\mu_0}{2}f_{14}^n(\xi_u^n) + \tilde e_1(\xi_u^n) -(\div\theta_u^n,\xi_p^n),
\end{split}
\end{equation}
\begin{equation}\label{eq:ID-f2}
\begin{split}
\rho\kappa(\delta_t\xi_\omega^n,\xi_\omega^n) + \eta^\prime \|\nabla\xi_\omega^n\|^2 + (\eta^\prime + \lambda^\prime)\|\div\xi_\omega^n\|^2 - 2\zeta(\curl\xi_u^n,\xi_\omega^n)
 + 4\zeta\|\xi_\omega^n\|^2 = -\rho\kappa f_{21}^n(\xi_\omega^n) + \mu_0f_{22}^n(\xi_\omega^n) +\tilde e_2^n(\xi_\omega^n),
\end{split}
\end{equation}
and
\begin{equation}
\label{eq:ID-f3}
\begin{split}
(\delta_t\xi_m^n,\xi_m^n)+\sigma\|\div\xi_m^n\|^2 + \sigma(\curl\xi_k^n,\xi_m^n) -\frac{1}{2}(\curl\xi_z^n,\xi_m^n) + \frac{1}{\tau}\|\xi_m^n\|^2 
- \frac{\chi_0}{\tau}(\xi_H^n,\xi_m^n) = f_{31}^n(\xi_m^n) \\+ \frac{1}{2}f_{32}^n(\xi_m^n) + \frac{1}{2}f_{33}^n(\xi_m^n) 
+ f_{34}^n(\xi_m^n) +\tilde e_3^n(\xi_m^n).
\end{split}
\end{equation}
Equations \eqref{eq:FHD-weak-6} and \eqref{eq:FHD-weak-disf-6} imply
$$
(\xi_k^n,\bs\Theta_h) - (\xi_m^n,\curl\bs\Theta_h) = (\theta_k^n,\bs\Theta_h) - (\theta_m^n,\bs\Theta_h)\qquad\forall~\bs\Theta_h \in \bs U_h.
$$
Taking $\bs\Theta_h = \xi_k^n$ in the above equation, we obtain
\begin{equation}\label{eq:ID-f4}
\|\xi_k^n\|^2 - (\xi_m^n,\curl\xi_k^n) = (\theta_k^n,\xi_k^n) - (\theta_m^n,\xi_k^n).
\end{equation}
Equations \eqref{eq:FHD-weak-7} and \eqref{eq:FHD-weak-disf-7} imply
\begin{equation}\label{eq:ID-f5}
(\xi_H^n,\bs G_h) + (\xi_\phi^n,\div\bs G_h) = (\theta_H^n,\bs G_h)\qquad\forall~\bs G_h \in \bs V_h.
\end{equation}
Equations \eqref{eq:FHD-weak-8} and \eqref{eq:FHD-weak-disf-8} imply
\begin{equation}\label{eq:ID-f6}
(\div\xi_H^n,r_h) + (\div\xi_m^n,r_h) = 0\qquad\forall~r_h \in W_h
\end{equation}
Taking $\bs G_h = \xi_H^n$ and $\bs G_h = \xi_m^n$ in \eqref{eq:ID-f5}, and taking $r_h = \xi_\phi^n$ in \eqref{eq:ID-f6}, we have
\begin{align*}
& \|\xi_H^n\|^2 + (\xi_\phi^n,\div\xi_H^n) = (\theta_H^n,\xi_H^n),\qquad\qquad
 (\xi_H^n,\xi_m^n) + (\xi_\phi^n,\div\xi_m^n) = (\theta_H^n,\xi_m^n),\\
& (\div\xi_H^n,\xi_\phi^n) + (\div\xi_m^n,\xi_\phi^n) = 0.
\end{align*}
These equations lead to 
\begin{equation}\label{eq:ID-f7}
(\xi_m^n,\xi_H^n)= -(\xi_\phi^n,\div\xi_m^n) + (\theta_H^n,\xi_m^n)  = -\|\xi_H^n\|^2 + (\theta_H^n,\xi_H^n) + (\theta_H^n,\xi_m^n)
\end{equation}
Substituting   \eqref{eq:ID-f4} and \eqref{eq:ID-f7} into \eqref{eq:ID-f3}, and 
adding the resulting equation with \eqref{eq:ID-f1} and \eqref{eq:ID-f2}, we obtain
\begin{equation}\label{eq:ID-f8}
\begin{split}
& \rho(\delta_t\xi_u^n,\xi_u^n) + \rho\kappa(\delta_t\xi_\omega^n,\xi_\omega^n) + (\delta_t\xi_m^n,\xi_m^n)  + \eta\|\nabla\xi_u^n\|^2 + \eta^\prime \|\nabla\xi_\omega^n\|^2 + (\eta^\prime + \lambda^\prime) \|\div\xi_\omega^n\|^2 + \sigma\|\div\xi_m^n\|^2 \\
&\qquad \qquad + \sigma\|\xi_k^n\|^2
 + \frac{1}{\tau} \|\xi_m^n\|^2 + \zeta \|\curl\xi_u^n - 2\xi_\omega^n\|^2 + \frac{\chi_0}{\tau} \|\xi_H^n\|^2 \\
 = & \tilde J_1^n + \frac{1}{2}(\curl\xi_z^n,\xi_m^n) - \rho f_{11}^n(\xi_u^n)  + \mu_0f_{12}^n(\xi_u^n) + \frac{\mu_0}{2}f_{13}^n(\xi_u^n) + \frac{\mu_0}{2} f_{14}^n(\xi_u^n)  - \rho\kappa f_{21}^n(\xi_\omega^n)
 \\
&  \qquad\qquad + \mu_0f_{22}^n(\xi_\omega^n) 
+ f_{31}^n(\xi_m^n) + \frac{1}{2}f_{32}^n(\xi_m^n) +\frac{1}{2}f_{33}^n(\xi_m^n) + f_{34}^n(\xi_m^n).
\end{split}
\end{equation}

Taking $\bs F_h = \xi_H^n$ in \eqref{eq:err-full-3} and using the fact that $\curl_h\xi_H^n = \curl_hI_h^d\bs H - \curl_h\bs H_h = Q_h^d\curl\bs H = 0$, we obtain
\begin{equation}\label{eq:ID-f8}
\begin{split}
(\delta_t\xi_m^n,\xi_H^n) + \sigma(\div\xi_m^n,\div\xi_H^n) + \frac{1}{\tau}(\xi_m^n,\xi_H^n) - \frac{\chi_0}{\tau}\|\xi_H^n\|^2 = f_{31}^n(\xi_H^n) 
+ \frac{1}{2}f_{32}^n(\xi_H^n) + \frac{1}{2}f_{33}^n(\xi_H^n) + f_{34}^n(\xi_H^n) + \tilde e_3^n(\xi_H^n)
\end{split}
\end{equation}
Equations \eqref{eq:FHD-weak-8} and \eqref{eq:FHD-weak-disf-8} imply
\begin{equation}\label{eq:ID-f9}
\mu_0(\div\xi_H^n + \div\xi_m^n,r_h) = 0\qquad\forall~r_h \in W_h.
\end{equation}
The fact that $\curl_h\xi_h^n = 0$ implies there exists $\psi_h \in W_h$ such that $\xi_H = \grad_h \psi_h$. Subtracting \eqref{eq:ID-f9} at the time level $t_{n-1}$ from \eqref{eq:ID-f9} at the time level $t_n$, and choosing the resulting equation with $r_h = \psi_h$, we obtain
$$
(\delta_t\xi_m^n,\xi_H^n) = -(\delta_t\xi_H^n,\xi_H^n).
$$
Taking $r_h = \div\xi_H^n$ in \eqref{eq:ID-f9}, we get
$$
\|\div\xi_H^n\|^2  = -(\div\xi_m^n,\div\xi_H^n).
$$
Thus, we have
\begin{equation}\label{eq:ID-f10}
\begin{split}
(\delta_t\xi_H^n,\xi_H^n) + \sigma\|\div\xi_H^n\|^2 + \frac{1+\chi_0}{\tau}\|\xi_H^n\|^2= -f_{31}^n(\xi_H^n) 
- \frac{1}{2}f_{32}^n(\xi_H^n) 
- \frac{1}{2}f_{33}^n(\xi_H^n) - f_{34}^n(\xi_H^n) + \tilde J_2^n.
\end{split}
\end{equation}
Multiplying \eqref{eq:ID-f10} with $\mu_0$, and adding the resulting equation with \eqref{eq:ID-f4}, we  finally get the desired result.
\end{proof}

Similar as Lemma \ref{lem:J1J2}, by using the Taylor's formula 
$$
v(t_{n-1}) = v(t_n) - \dt v_t(t_n) + \int_{t_n}^{t_{n-1}} v_{tt}(s)(s- t_n) \dd s,
$$
 the following estimates for $\tilde J_1$ and $\tilde J2$ can be obtained.
\begin{lemma}\label{lem:J1J2-f}
Under  Assumptions \ref{ass:initial-data} - \ref{ass:omega}, there hold
\begin{align*}
 \tilde J_1& \lesssim \frac{h^2}{\dt} \left(\|\xi_u^n\|\int_{t_{n-1}}^{t_n} \|\bs u_t\|_2\dd s  + \|\xi_\omega^n\|\int_{t_{n-1}}^{t_n}\|\bs\omega_t\|_2\dd s\right) + \frac{h}{\dt}\|\xi_m^n\| \int_{t_{n-1}}^{t_n} \|\bs m_t\|_1\dd s \\
 & + \left(\|\xi_u^n\|\int_{t_{n-1}}^{t_n} \|\bs u_{tt}\| \dd s + \|\xi_\omega^n\|\int_{t_{n-1}}^{t_n}\|\bs\omega_{tt}\|\dd s + \|\xi_m^n\| \int_{t_{n-1}}^{t_n}\|\bs m_{tt}\|\dd s\right) \\
&+ h^2\|\bs \omega^n\|_2(\|\nabla\xi_u^n\|+\|\xi_\omega^n\|) + h(\|\bs u^n\|_2 +\|p^n\|_1)\|\nabla\xi_u^n\| + h\|\bs \omega^n\|_2\|\nabla\xi_\omega^n\| + h\|\bs \omega^n\|_2\|\xi_\omega^n\|\\
& + h(\|\bs m^n\|_1 + \|\bs H^n\|_1)\|\xi_m^n\| + h\|\bs H^n\|_1\|\xi_H^n\| + h\|\bs u^n\|_2\|\xi_p^n\| + h(\|\bs k^n\|_1 + \|\bs z^n\|_1)\|\xi_k^n\|,\\
\tilde J_2 &\lesssim \frac{h}{\dt}\|\xi_H^n\| \int_{t_{n-1}}^{t_n} \|\bs m_t\|_1\dd s + \|\xi_H^n\|\int_{t_{n-1}}^{t_n} \|\bs m_{tt}\| \dd s  + h(\|\bs m^n\|_1 + \|\bs H^n\|_1)\|\xi_H^n\| + h\|\bs H^n\|_1\|\xi_m^n\|.
\end{align*}
\end{lemma}

Now we turn to estimate $\xi_p^n$ and   have the following result.
\begin{lemma}\label{lem:xip-f}
Under   Assumptions \ref{ass:initial-data} - \ref{ass:omega}, for any $1\leq L \leq N$ there holds
\begin{align*}
\dt\sum\limits_{n=1}^L&  \|\xi_p^n\| \lesssim \|\xi_u^L\| + h^2(\int_0^{t_L}\|\bs u_t\|_2\dd s +\|\bs\omega\|_{L^\infty(H^2)})  + h(\|\bs u\|_{L^\infty(H^2)} + \|p\|_{L^\infty(H^1)}) \\
& +h(\|\bs u\|_{L^\infty(H^2)}^2 + \|\bs m\|_{L^\infty(H^2)}\|\bs H\|_{L^\infty(H^2)} + \|\bs k\|_{L^\infty(H^1)}\|\bs H\|_{L^\infty(H^1)} ) \\
& + \dt^{3/2}\int_0^{t_L}\|\bs u_{tt}\|^2\dd s   + \dt\sum\limits_{n=1}^L(\|\nabla\xi_u^n\| + \|\xi_\omega^n\|) \\
& +h\dt\sum\limits_{n = 1}^L (\|\bs u^n\|_2\|\nabla\xi_u^n\| +\|\bs H^n\|_2\|\xi_m^n\|_n + \|\bs H^n\|_2\|\xi_k^n\|)\\
& + \dt\sum\limits_{n=1}^L \left(\|\bs u^n\|_2\|\nabla\xi_u^n\| + \|\xi_u^n\|^{1/2}\|\nabla\xi_u^n\|^{3/2} +\|\bs H^n\|_2\|\xi_m^n\| + \|\bs H^n\|_1\|\div\xi_m^n\|\right) \\
& + \dt\sum\limits_{n=1}^L \left( \|\xi_m^n\|^{1/2}\|\xi_m^n\|_{n}^{1/2}\|\div\xi_H^n\|  +\|\bs m^n\|_2\|\div\xi_H^n\| +\|\bs H^n\|_1\|\xi_k^n\| + \|\bs k^n\|_1\|\xi_H^n\|\|\right).
\end{align*}
\end{lemma}
\begin{proof}
The inf-sup condition \eqref{eq:inf-sup} and the equation \eqref{eq:err-full-1} imply that for any $1 \leq L \leq N$,  
\begin{align*}
\dt\sum\limits_{n = 1}^L \|\xi_p^n\| & \lesssim \dt \sup\limits_{\bs v_h \in \bs S_h} \sum\limits_{n = 1}^L \frac{(\xi_p^n,\div\bs v_h)}{\|\bs v_h\|_1} \\
& \lesssim \|\xi_u^L\| + \dt\sum\limits_{n=1}^L \|\nabla\xi_u^n\| + \dt\sum\limits_{n=1}^L\|\xi_\omega^n\| + h^2(\int_0^{t_L}\|\bs u_t\|_{2}\dd s + \|\bs\omega\|_{L^\infty(H^2)}) \\
& \quad + h(\|\bs u\|_{L^\infty(H^2)} + \|\tilde p\|_{L^\infty(H^1)}) + \dt^{3/2}\int_0^{t_L} \|\bs u_{tt}\|^2\dd s+ \tilde T_1,
\end{align*}
with
$$
\tilde T_1 := \sup\limits_{\bs v_h \in \bs S_h} \dt\sum\limits_{n=1}^L \frac{\left|\rho f_{11}^n(\bs v_h) - \mu_0 f_{12}^n(\bs v_h) - \frac{\mu_0}{2} f_{13}^n(\bs v_h) - \frac{\mu_0}{2} f_{14}^n(\bs v_h)\right|}{\|\bs v_h\|_1}.
$$
Following a similar routine as in the proof of Lemma \ref{lem:xip}, for any $\bs v_h\in \bs S_h$ we have
\begin{align*}
|f_{11}^n(\bs v_h) |
& \lesssim h(\|\bs u^n\|_2^2 +\|\bs u^n\|_2\|\nabla\xi_u^n\|)\|\bs v_h\|_1+ (\|\bs u^n\|_2\|\nabla\xi_u^n\| +\|\xi_u^n\|^{1/2}\|\nabla\xi_u^n\|^{3/2})\|\bs v_h\|_1 ,
\end{align*}
\begin{align*}
|f_{12}^n(\bs v_h)| 
\lesssim & h\|\bs m^n\|_2\|\bs H^n\|_2 \|\nabla\bs v_h\| + \|\bs H^n\|_2\|\xi_m^n\|\|\bs v_h\|_1 +\|\bs H^n\|_1\|\div\xi_m^n\|\|\bs v_h\|_1\\
& + \|\xi_m^n\|^{1/2}\|\xi_m^n\|_n^{1/2}\|\div\xi_H^n\|\|\bs v_h\|_1 + \|\bs m^n\|_2\|\div\xi_H^n\|\|\bs v_h\|_1,
\end{align*}
\begin{align*}
|f_{13}^n(\bs v_h)| \lesssim & h(\|\bs k^n\|_1\|\bs H^n\|_1 + \|\bs H^n\|_2\|\xi_k^n\|) \|\bs v_h\|_1 + (\|\bs H^n\|_1\|\xi_k^n\| + \|\bs k^n\|_1\|\xi_H^n\|)\|\bs v_h\|_1\\
& + \|\xi_k^n\|\|\xi_H^n\|^{1/2}\|\xi_H^n\|_n^{1/2},
\end{align*}
and
\begin{align*}
|f_{14}^n(\bs v_h)| 
\lesssim & h(\|\bs m^n\|_1\|\|\bs H^n\|_2 + \|\bs H^n\|_2\|\xi_m^n\|_n)\|\bs v_h\|_1 + (\|\bs H^n\|_2\|\xi_m^n\| +\|\bs m\|_2\|\xi_H^n\|)\|\bs v_h\|_1 \\
& +\|\xi_m^n\|^{1/2}\|\xi_m^n\|_n^{1/2}\|\div\xi_H^n\|\|\bs v_h\|_1.
\end{align*}
As a result, the desired result follows.
\end{proof}

\begin{lemma}
\label{lem:err-xi-f}
Under   Assumptions \ref{ass:initial-data} - \ref{ass:omega}, there holds
$$
\sup\limits_{1\leq n\leq N}\left(\mathfrak E_h^n + \dt\sum\limits_{l = 1}^n \mathfrak F_h^l  \right) \lesssim h^2 + \dt^2.
$$ 
\end{lemma}
\begin{proof}
By Lemmas \ref{lem:J3} - \ref{lem:J6}, \ref{lem:err-ID-f} and \ref{lem:J1J2-f} we have
\begin{equation}\label{eq:err-pf-1}
\frac{1}{\dt} \left(\mathfrak E_h^n - \mathfrak E_h^{n-1} \right) + \mathfrak F_h^n \leq C h^2\mathfrak T_1^n +C h\mathfrak T_2^n  + C\mathfrak T_3^n + C \mathfrak T_4^n,
\end{equation}
where
\begin{align*}
\mathfrak T_1^n &: = \frac{1}{\dt}\left(\|\xi_u^n\|\int_{t_{n-1}}^{t_n} \|\bs u_t\|_2\dd s + \|\xi_\omega^n\| \int_{t_{n-1}}^{t_n} \|\bs\omega_t\|_2\right) + \|\bs \omega^n\|_2(\|\nabla\xi_u^n\|+\|\xi_\omega^n\|) \\
&\quad +  \|\bs u^n\|_2^2\|\nabla\xi_u^n\| + \|\bs m^n\|_2\|\bs\omega^n\|_2\|\xi_H^n\| + \|\bs u^n\|_2\|\bs\omega^n\|_2\|\nabla\xi_\omega^n\|, 
\end{align*}
\begin{align*}
\mathfrak T_2^n & :=\frac{1}{\dt} \left(\|\xi_m^n\|\int_{t_{n-1}}^{t_n}\|\bs m_t\|_1\dd t   + \|\xi_H^n\|\int_{t_{n-1}}^{t_n} \|\bs m_t\|_1\dd s \right)\\
&\quad  +(\|\bs u^n\|_2^2 +\|\bs m^n\|_2\|\bs H^n\|_2 +\|\bs k^n\|_1\|\bs H^n\|_2) \|\xi_u^n\| + \|\bs u^n\|_2\|\xi_p^n\| +\|\bs u^n\|_2\|\bs m^n\|_1 \|\xi_m^n\|_n\\
&\quad + (\|\bs u^n\|_2+\|p^n\|_1+\|\bs H^n\|_2(\|\bs m^n\|_2+\|\bs k^n\|_1))\|\nabla\xi_u^n\| \\
&\quad + (\|\bs\omega^n\|_2+\|\bs u^n\|_2\|\bs\omega^n\|_2+\|\bs H^n\|_2\|\bs m^n\|_2)\|\xi_\omega^n\| +\|\bs\omega^n\|_2 \|\nabla\xi_\omega^n\| \\
&\quad +(\|\bs m^n\|_1 +\|\bs H^n\|_1 +\|\bs u^n\|_2(\|\bs m^n\|_2 + \|\bs k^n\|_1)  + \|\bs m^n\|_2(\|\bs\omega^n\|_2+\|\bs m^n\|_1)) \|\xi_m^n\| \\
&\quad  + \left( \|\bs m^n\|_1 + \|\bs H^n\|_1 +(\|\bs u^n\|_2+\|\bs\omega^n\|_2)\|\bs m^n\|_2 +\|\bs k^n\|_1\|\bs u^n\|_2\right) \|\xi_H^n\|\\
& \quad   + (\|\bs z^n\|_1 +\|\bs u^n\|_2\|\bs m^n\|_2) \|\xi_k^n\| + \|\bs k^n\|_1\|\xi_z^n\|+ \|\bs u^n\|_2\|\bs m^n\|_1 \|\div\xi_H^n\|,
\end{align*}
\begin{align*}
\mathfrak T_3^n & := \|\bs u^n\|_2 \|\xi_u^n\|\|\nabla\xi_u^n\| + \|\bs\omega^n\|_2 \|\xi_u^n\|\|\nabla\xi_\omega^n\| + \|\bs\omega^n\|_2\|\nabla\xi_u^n\|\|\xi_\omega^n\| \\
&\quad +(\|\bs H^n\|_1+\|\bs m^n\|_1)\|\xi_m^n\|\|\nabla\xi_\omega^n\| + \|\bs H^n\|_2\|\xi_u^n\|\|\xi_k^n\| + \|\bs u^n\|_2\|\xi_k^n\|\|\xi_H^n\|\\
& \quad + (\|\bs H^n\|_2+\|\bs k^n\|_1 +\|\bs m^n\|_2)\|\nabla\xi_u^n\|\|\xi_m^n\| + (\|\bs H^n\|_2+\|\bs m^n\|_2)\|\xi_u^n\|\|\xi_m^n\|_n\\
&\quad + (\|\bs u^n\|_2 +\|\bs\omega^n\|_2)\|\xi_m^n\|\|\div\xi_H^n\| +\|\bs u^n\|_2 \|\xi_H^n\|\|\div\xi_m^n\|,\\
\mathfrak T_4^n &: = \|\xi_u^n\|\int_{t_{n-1}}^{t_n}\|\bs u_{tt}\|\dd s + \|\xi_\omega^n\|\int_{t_{n-1}}^{t_n}\|\bs\omega_{tt}\|\dd s  + (\|\xi_m^n\| + \|\xi_H^n\|) \|\int_{t_{n-1}}^{t_n}\|\bs m_{tt}\|\dd s.
\end{align*}
For any $1\leq L \leq N$, summing up \eqref{eq:err-pf-1} with $n$ from $1$ to $L$, we get
\begin{equation}\label{eq:err-pr-f-2}
\mathfrak E_h^L + \dt\sum\limits_{n = 1}^L \mathfrak F_h^n  \leq Ch^2\dt \sum\limits_{n = 1}^L \mathfrak T_1^n + Ch\dt \sum\limits_{n = 1}^L \mathfrak T_2^n+ C\dt \sum\limits_{n = 1}^L \mathfrak T_3^n + C\dt \sum\limits_{n = 1}^L \mathfrak T_4^n. 
\end{equation}
Using Lemmas \ref{lem:xiz} and \ref{lem:xip-f}, the inverse inequality, the Young's inequality, and the regularity results in Theorems \ref{the:reg-2} - \ref{the:reg-aux}, we have
\begin{align*}
Ch^2\dt\sum\limits_{n = 1}^L \mathfrak T_1^n & \leq Ch^2 + C\dt\sum\limits_{n = 1}^L \left(\|\xi_u^n\|^2 + \|\xi_\omega^n\|^2 + \|\xi_H^n\|^2 \right),\\
Ch\dt\sum\limits_{n = 1}^L \mathfrak T_2^n & \leq Ch^2 + C\dt^3+ \frac{\rho}{2}\|\xi_u^L\|^2 + C\dt\sum\limits_{n = 1}^L (\|\xi_u^n\|^2 + \|\xi_\omega^n\|^2 + \|\xi_m^n\|^2 + \|\xi_H^n\|^2 )  + \frac{\eta}{4}\dt\sum\limits_{n = 1}^L\|\nabla\xi_u^n\|^2 \\
&\quad+ \frac{\eta^\prime}{4}\dt\sum\limits_{n=1}^L\|\nabla\xi_\omega^n\|^2 + \frac{\sigma}{8}\dt\sum\limits_{n = 1}^L (\|\xi_m^n\|_n^2 +\|\xi_k^n\|^2)
 + \frac{\mu_0\sigma}{4}\dt\sum\limits_{n = 1}^L \|\div\xi_H^n\|^2,\\
C\dt\sum\limits_{n = 1}^L\mathfrak T_3^n & \leq \frac{\eta}{4}\dt\sum\limits_{n = 1}^L\|\nabla\xi_u^n\|^2 + \frac{\eta^\prime}{4}\dt\sum\limits_{n=1}^L\|\nabla\xi_\omega^n\|^2 + \frac{\sigma}{8}\dt\sum\limits_{n = 1}^L (\|\xi_m^n\|_n^2 +\|\xi_k^n\|^2) \\
&
\quad + \frac{\mu_0\sigma}{4}\dt\sum\limits_{n = 1}^L \|\div\xi_H^n\|^2 + C\dt\sum\limits_{n = 1}^L (\|\xi_u^n\|^2 + \|\xi_\omega^n\|^2 + \|\xi_m^n\|^2 + \|\xi_H^n\|^2 ) ,\\
C\dt\sum\limits_{n = 1}^L\mathfrak T_4^n & \leq C\dt^2 + C\dt\sum\limits_{n = 1}^L (\|\xi_u^n\|^2 + \|\xi_\omega^n\|^2 + \|\xi_m^n\|^2 + \|\xi_H^n\|^2 ) ,
\end{align*}
which, together with Lemma \eqref{eq:err-pr-f-2}, indicate 
\begin{equation*}
\mathfrak E_h^L + \dt\sum\limits_{n = 1}^L \mathfrak F_h^n  \lesssim h^2+ \dt^2 + \dt\sum\limits_{n = 1}^L \left(\mathfrak E_h^n  + \dt\sum\limits_{l = 1}^n \mathfrak F_h^l \right).
\end{equation*}
This inequality plus the discrete Gr\"onwall's inequality (cf. \citet{Clark1987}) implies the desired result.
\end{proof}

Using Lemma \ref{lem:err-xi-f} and the triangular inequality, we have the following error estimates of the fully discrete scheme \eqref{eq:FHD-weak-disf-1} - \eqref{eq:FHD-weak-disf-8}.
\begin{theorem}\label{estimate-full}
Under   Assumptions \ref{ass:initial-data} - \ref{ass:omega}, there hold
$$
\sup\limits_{1\leq n \leq N}\left(\hbox{err}_1^n + \dt\sum\limits_{l = 1}^n \hbox{err}_2^n \right) \lesssim h^2 + \dt^2
$$
and
$$
\dt\sum\limits_{n = 1}^L\|\tilde p(\cdot, t_n) - \tilde p_h^n\| \lesssim h + \dt.
$$
\end{theorem}

\section{Numerical experiments}\label{sec:num}
In this section, we provide numerical examples to verify the performance of the fully discrete scheme \eqref{eq:FHD-weak-disf-1} - \eqref{eq:FHD-weak-disf-8}. The numerical experiments are performed by using the iFEM package by \citet{Chen2009}, and the nonlinear system \eqref{eq:FHD-weak-disf-1} - \eqref{eq:FHD-weak-disf-8} is solve by the following quasi-Newton iteration with $M = 2$:
\begin{alg}
Given $\bs u_h^{n-1}$, $\bs\omega_h^{n-1}$ and $\bs m_h^{n-1}$, to find $\bs u_h^n$, $p_h^n$, $\bs\omega_h^n$, $\bs m_h^n$, $\bs z_h^n$, $\bs k_h^n$, $\bs H_h^n$ and $\phi_h^n$ through three steps:
\begin{itemize}
  \item[Step 1.] Set $\bs u_h^- = \bs u_h^{n-1}$, $\bs\omega_h^- = \bs \omega_h^{n-1}$ and $\bs m_h^- = \bs m_h^{n-1}$.
  \item[Step 2.] For $\hbox{iter} = 1,2,\dots,M$ do
  \begin{itemize}
  \item[(1)] Solving the saddle point system: Find $\bs H_h \in \bs V_h$ and $\phi_h\in W_h$ such that
  \begin{equation*}
  \left\{ 
  \begin{array}{ll}
  (\bs H_h,\bs G_h) + (\phi_h,\div\bs G_h) = 0 & \forall~\bs G_h \in \bs V_h,\\
  \mu_0(\div\bs H_h,r_h) = -(\div \bs H_e^n,r_h) - \mu_0(\div\bs m_h^-,r_h) &\forall~r_h \in W_h.
  \end{array}
  \right.
  \end{equation*}
  \item[(2)] Solving the angular momentum equation: Find $\bs\omega_h \in \bs \Sigma_h$ such that for all $\bs s_h \in \bs\Sigma_h$, it holds
  \begin{equation*}
  \begin{split}
  (\rho\kappa + 4\zeta)(\bs\omega_h,\bs s_h) + \eta^\prime\dt(\nabla\bs\omega_h,\nabla\bs s_h) + (\eta^\prime + \lambda^\prime)\dt (\div\bs\omega_h,\div\bs s_h) 
= \rho\kappa(\bs\omega_h^{n-1},\bs s_h) \\
 - \rho\kappa\dt b(\bs u_h^-,\bs\omega_h^-,\bs s_h) 
    + \mu_0\dt (\bs m_h^-\times\bs H_h^-,\bs s_h) 
    + 2\zeta\dt (\curl\bs u_h^-,\bs s_h).
  \end{split}  \end{equation*}
  \item[(3)] Solving the magnetization equation: Find $\bs m_h \in \bs V_h$, $\bs z_h \in \bs U_h$ and $\bs k_h\in \bs U_h$ such that for all $\bs F_h \in \bs V_h$, $\bs\Lambda_h \in \bs U_h$ and $\bs\Theta_h \in \bs U_h$, there hold
  \begin{align*}
  & \left( 1+\frac{\dt}{\tau}\right)(\bs m_h,\bs F_h) + \sigma\dt(\div\bs m_h,\div\bs F_h)   + \sigma\dt(\curl\bs k_h,\bs F_h)- \frac{\dt}{2}(\bs u_h^-\times\bs k_h,\bs F_h)  - \frac{\dt}{2}(\curl\bs z_h,\bs F_h) \\
&\qquad  = (\bs m_h^{n-1},\bs F_h) + \frac{\chi_0\dt}{\tau}(\bs H_h,\bs F_h) + \dt c(\bs u_h^-,\bs m_h^-,\bs F_h)+ \frac{\dt}{2}(\bs m_h^-\times \curl\bs u_h^-,\bs F_h)+ \dt(\bs\omega_h\times \bs m_h^-,\bs F_h),\\
  &(\bs z_h,\bs\Lambda_h)- (\bs u_h^-\times \bs m_h,\bs\Lambda_h) = 0,\\
  & (\bs k_h,\bs\theta_h) - (\bs m_h,\curl\bs\Theta_h) = 0.
  \end{align*}
  \item[(4)] Solving the Navier-Stokes equation: Find $\bs u_h \in \bs S_h$ and $\tilde p_h \in L_h$ such that for all $\bs v_h \in \bs S_h$ and $q_h \in L_h$, there hold
  \begin{align*}
 & \rho(\bs u_h,\bs v_h) + \dt\eta(\nabla\bs u_h,\nabla\bs v_h) + \dt\zeta(\curl\bs u_h,\curl\bs v_h) - \dt(\tilde p_h,\div\bs v_h)  =(\bs u_h^{n-1},\bs v_h) - \dt\rho b(\bs u_h^-,\bs u_h^-,\bs v_h)\\
 &\qquad\quad + \dt\mu_0 c(\bs v_h,\bs m_h,\bs H_h)  + \frac{\mu_0\dt}{2}(\bs v_h\times \bs k_h,\bs H_h) + \frac{\mu_0\dt}{2}(\bs m_h\times\curl\bs v_h,\bs H_h) + 2\dt\zeta(\bs\omega_h,\curl\bs v_h),\\
&  (\div\bs u_h,q_h)  = 0.
  \end{align*}
  \item[(5)] Set $\bs u_h^- = \bs u_h$, $\bs\omega_h^- = \bs\omega_h$ and $\bs m_h^- = \bs m_h$.
  \end{itemize}

  \item[Step 3.] Let $\bs u_h^n = \bs u_h$, $p_h^n = p_h$, $\bs\omega_h^n = \bs\omega_h$, $\bs H_h^n = \bs H_h$, $\bs z_h^n = \bs z_h$, $\bs k_h^n = \bs k_h$ and $\phi_h^n = \phi_h$.
\end{itemize}
\end{alg}

\begin{figure}[!h]
\centering
\includegraphics[height=4.5cm,width=9cm]{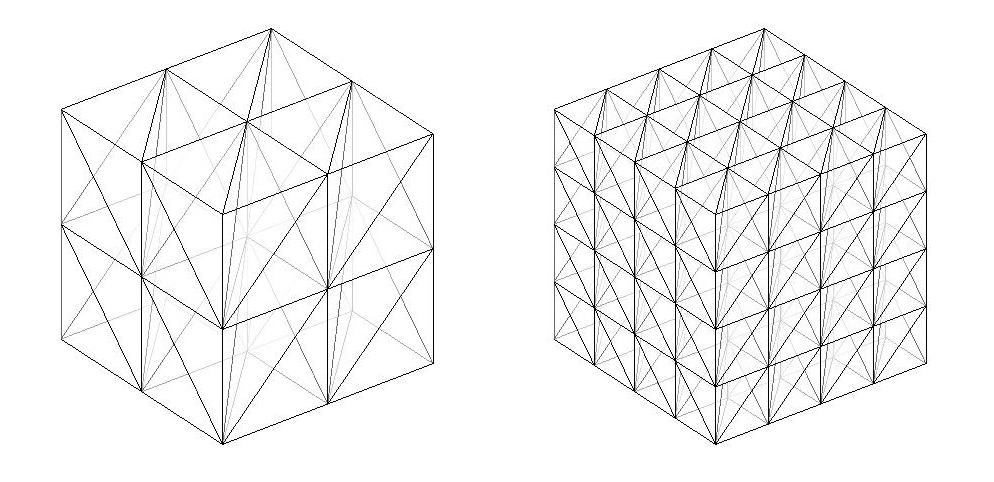}
\caption{ The domain $\Omega = (0,1)^3$: $2\times2\times 2$   (left)     and $4\times 4\times 4$ (right) meshes.}\label{Fig:mesh}
\end{figure}

In the following numerical examples, we take $\Omega = (0,1)^3$ and use $K \times K \times K$ uniform tetrahedral meshes (cf. Fig. \ref{Fig:mesh}) with $K~=~4,~8,~16,~32$. In the first example, we take the temporal step size as $\dt = h/\sqrt{3} = 1/K$. From Theorem \ref{estimate-full}, we easily see that the theoretical accuracy of the scheme is $\mathcal O(h+\dt)$.

\begin{example}
\label{ex:ex-1}
The exact solution, $(\bs u,~\tilde p,~\bs\omega,~\bs m,~\bs H,~\phi)$, of the FHD model \eqref{eq:FHD-1} - \eqref{eq:FHD-6} with $T = 1$ is given by
\begin{align*}
& \bs u(x,y,z,t) = \sin(t)[\sin(\pi y),\sin(\pi z),\sin(\pi x)]^\intercal, \qquad \qquad \tilde p(x,y,z,t) = 120x^2yz - 40 y^3z - 40 y z^3,\\
& \bs\omega(x,y,z,t) = \sin(t) [(x^2-x)(y^2-y)(z^2-z),0,0]^\intercal,\qquad \bs m(x,y,z,t) = \sin(t)[\sin(\pi x)\sin(\pi y)\sin(\pi z),0,0]^\intercal,\\
& \bs H(x,y,z,t) = 2000\sin(t)\begin{pmatrix}
(2x-1)(x^2-x)(y^2-y)^2(z^2-z)^2 \\
(x^2-x)^2(2y-1)(y^2-y)(z^2-z)^2\\
(x^2-x)^2(y^2-y)^2(2z-1)(z^2-z)
\end{pmatrix},\\
& \phi(x,y,z,t) = 1000\sin(t)(x^2-x)^2(y^2-y)^2(z^2-z)^2. 
\end{align*}
The parameters $\rho,~\kappa,~\eta,~\zeta,~\mu_0,~\sigma,~\eta^\prime,~\lambda^\prime,~\tau$ and $\chi_0$ are all chosen as $1$. Numerical results of the relative errors of the discrete solutions at the ending time $T=1$ are listed in Table \ref{tab:ex-1}.
\end{example}
\begin{table}[h!]
\footnotesize
\begin{center}
\caption{Errors and convergence orders at   $T = 1$ for Example \ref{ex:ex-1}: $\dt=1/K$.}
\label{tab:ex-1}
\begin{tabular}{c|c|c|c|c|c|c}
\hline  $K$ & $\frac{\|\bs u(T) - \bs u_{h}^N\|}{\|\bs u(T)\|}$  &$ \frac{|\bs u(T) - \bs u_{h}^N|_1}{|\bs u(T)|_1} $ & $\frac{\|\tilde p(T) - \tilde p_h^N\| }{\|\tilde p(T)\|}$&$ \frac{\|\bs m(T) - \bs m_{h}^N\|}{\|\bs m(T)\|}$ & $\frac{\|\nabla\cdot(\bs m(T) - \bs m_{h}^N)\|}{\|\nabla\cdot\bs m(T)\|}$ & $\frac{\|\bs H(T) - \bs H_h^N\|}{\|\bs H(T)\|}$ \\
\hline
$4$ & $0.0554$ & $0.2270$ & $0.0678$ & $0.3204$ & $0.2848$ &    $0.4461$   \\
\hline
$8$ & $0.0140$  & $0.1117$  & $0.0162$  & $0.1638$ & $0.1446$ &   $0.2276$  \\
\hline
$ 16$ & $0.0035$  & $0.0555$ & $0.0041$ & $0.0824$ & $0.0740$ &   $0.1145$ \\
\hline
$ 32$& $0.0009$  &  $0.0277$  & $0.0011$ &  $0.0413$ & $0.0388$  &  $0.0574$ \\
\hline
order & $1.9942$ & $1.0111$ &  $1.9701$ & $0.9858$ &  $0.9585$ &   $0.9865 $ \\ 
\hline
\hline
$K$ &  $\frac{\|\nabla\cdot(\bs H(T) - \bs H_h^N)\|}{\|\nabla\cdot\bs H(T)\|}$ & $\frac{\|\bs z(T) - \bs z_h^N\|}{\|\bs z(T)\|}$ & $\frac{\|\bs k(T) - \bs k_h^N\|}{\|\bs k(T)\|}$ & $\frac{\|\bs\omega(T) - \bs\omega_h^N\|}{\|\bs\omega(T)\|}$ & $\frac{|\bs\omega(T) - \bs\omega_h^N|_1}{|\bs\omega(T)|_1}$ & $\frac{\|\phi(T) - \phi_h^T\|}{\|\phi(T)\|}$ \\
\hline
$ 4$ &  $0.4627$ & $0.3249$ & $0.3089$ & $0.4362$ & $0.7637$ &    $0.4026$\\
\hline
$ 8$ & $0.2379$ & $0.1624$ & $0.1599$ & $0.1734$ & $0.3389$ &    $0.1924$\\
\hline
$ 16$  & $0.1200$ & $0.0818$ & $0.0817$ & $0.0608$ & $0.1503$ &   $0.0950$\\
\hline
$ 32$ & $0.0603$ & $0.0411$ & $0.0420$ & $0.0219$ & $0.0704$ &   $0.0473$\\
\hline
order & $0.9803$ & $0.9940$ & $0.9598$ & $1.4395$ & $1.1462$ &   $1.0295$\\
\hline
\end{tabular}
\end{center}
\end{table}

\begin{example}\label{ex:ex-2}
The exact solution, $(\bs u,~\tilde p,~\bs\omega,~\bs m,~\bs H,~\phi)$, of the FHD model \eqref{eq:FHD-1} - \eqref{eq:FHD-6} with $T = 4$ is given by
\begin{align*}
& \bs u(x,y,z,t) = e^{-t}[\sin(\pi y),\sin(\pi z),\sin(\pi x)]^\intercal, \qquad\qquad \tilde p(x,y,z,t) = 120x^2yz - 40 y^3z - 40 y z^3,\\
& \bs\omega(x,y,z,t) = e^{-t} [(x^2-x)(y^2-y)(z^2-z),0,0]^\intercal,\qquad \bs m(x,y,z,t) = e^{-t}[\sin(\pi x)\sin(\pi y)\sin(\pi z),0,0]^\intercal,\\
& \bs H(x,y,z,t) = 2000e^{-t}\begin{pmatrix}
(2x-1)(x^2-x)(y^2-y)^2(z^2-z)^2 \\
(x^2-x)^2(2y-1)(y^2-y)(z^2-z)^2\\
(x^2-x)^2(y^2-y)^2(2z-1)(z^2-z)
\end{pmatrix},\\
& \phi(x,y,z,t) = 1000 e^{-t}(x^2-x)^2(y^2-y)^2(z^2-z)^2. 
\end{align*} 
The parameters $\rho,~\kappa,~\eta,~\zeta,~\mu_0,~\sigma,~\eta^\prime,~\lambda^\prime,~\tau$ and $\chi_0$ are all chosen as $1$. Numerical results of the relative errors of the discrete solutions at the ending time $T=2$ are listed in Table \ref{tab:ex-2}.
\end{example}

\begin{table}[h!]
\footnotesize
\begin{center}
\caption{Errors and convergence orders at   $T = 2$ for Example \ref{ex:ex-2}: $\dt=1/K$.}
\label{tab:ex-2}
\begin{tabular}{c|c|c|c|c|c|c}
\hline  $K$ & $\frac{\|\bs u(T) - \bs u_{h}^N\|}{\|\bs u(T)\|}$  &$ \frac{|\bs u(T) - \bs u_{h}^N|_1}{|\bs u(T)|_1} $ & $\frac{\|\tilde p(T) - \tilde p_h^N\| }{\|\tilde p(T)\|}$&$ \frac{\|\bs m(T) - \bs m_{h}^N\|}{\|\bs m(T)\|}$ & $\frac{\|\nabla\cdot(\bs m(T) - \bs m_{h}^N)\|}{\|\nabla\cdot\bs m(T)\|}$ & $\frac{\|\bs H(T) - \bs H_h^N\|}{\|\bs H(T)\|}$ \\
\hline
$4$ & $0.0617$ &   $0.4169$ &  $0.0620$  &  $0.3201$  &  $0.2840$ &   $0.4452$   \\
\hline
$8$ & $0.0144$ &   $0.1457$ &   $0.0140$  &  $0.1637$  &  $0.1441$  &  $0.2274$   \\
\hline
$ 16$ & $0.0035$ &   $0.0604$ &   $0.0034$ &  $ 0.0824$ &  $ 0.0738$ &  $ 0.1145$ \\
\hline
$ 32$& $0.0009$ &  $0.0284$  &  $0.0009$ & $0.0412$ &   $0.0387$ &  $0.0573$  \\
\hline
order & $2.0528$ &  $1.2925$ &  $2.0584$  &  $0.9854$ &  $0.9584$  &  $0.9857$ \\ 
\hline
\hline
$K$ &  $\frac{\|\nabla\cdot(\bs H(T) - \bs H_h^N)\|}{\|\nabla\cdot\bs H(T)\|}$ & $\frac{\|\bs z(T) - \bs z_h^N\|}{\|\bs z(T)\|}$ & $\frac{\|\bs k(T) - \bs k_h^N\|}{\|\bs k(T)\|}$ & $\frac{\|\bs\omega(T) - \bs\omega_h^N\|}{\|\bs\omega(T)\|}$ & $\frac{|\bs\omega(T) - \bs\omega_h^N|_1}{|\bs\omega(T)|_1}$ & $\frac{\|\phi(T) - \phi_h^T\|}{\|\phi(T)\|}$ \\
\hline
$ 4$ &  $0.4627$ &   $0.3223$ & $0.3085$ & $ 0.4058$ & $0.7764$ &  $0.4006$\\
\hline
$ 8$ & $0.2379$ &  $0.1617$ & $0.1598$ & $0.1382$ & $0.3219$ &  $0.1917$\\
\hline
$ 16$  & $0.1200$ & $0.0816$ & $0.0817$ &  $0.0357$ & $0.1393$ & $ 0.0948$\\
\hline
$ 32$ & $0.0602$ & $0.0411$ & $0.0420$ & $0.0088$
& $0.0663$ &  $0.0473$\\
\hline
order & $0.9804$ & $0.9907$ & $0.9593$ &  $1.8447$ & $1.1832$ & $1.0278$\\
\hline
\end{tabular}
\end{center}
\end{table}

\begin{figure}[!h]
\centering
\includegraphics[height=6cm]{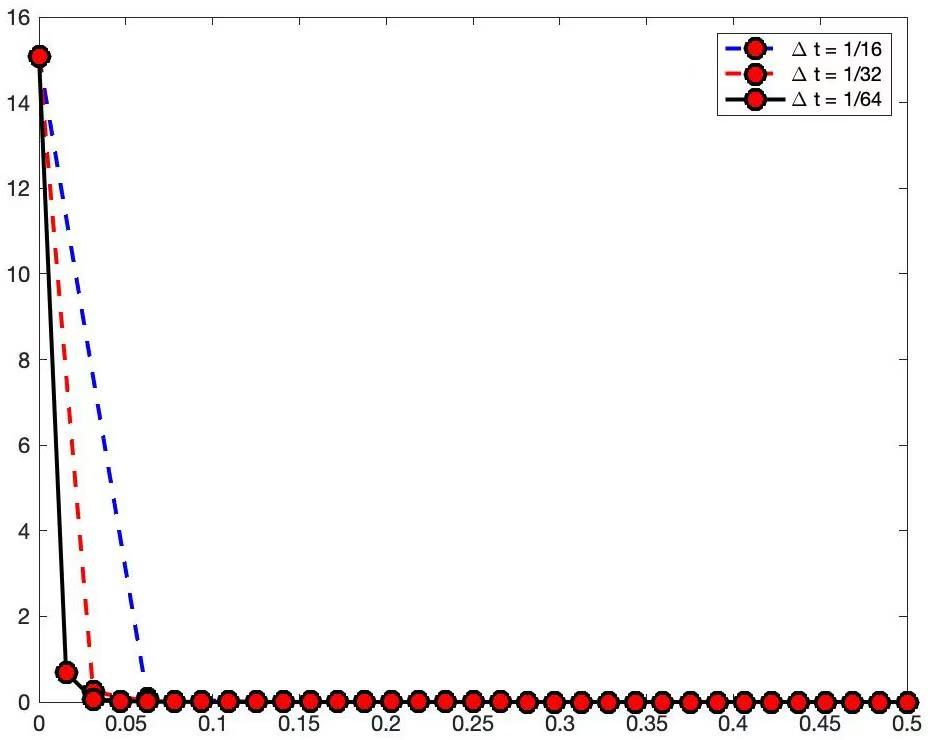}
\caption{ The energy curve with $K = 1/16$ and $\dt = 1/16,~1/32,~1/64$.}\label{Fig:ener}
\end{figure}

\begin{example}[Energy test] \label{ex:3}This example is to investigate the energy decaying phenomenon of the scheme (cf Theorem \ref{the:ener-dis-f}). We consider the FHD model \eqref{eq:FHD-1} - \eqref{eq:FHD-6} with the initial value functions
as
\begin{align*}
& \bs u_0(x,y,z) = [\sin(\pi y),\sin(\pi z),\sin(\pi x)]^\intercal, \qquad \qquad\bs\omega_0(x,y,z) = [(x^2-x)(y^2-y)(z^2-z),0,0]^\intercal,\\
& \bs m_0(x,y,z) = [\sin(\pi x)\sin(\pi y)\sin(\pi z),0,0]^\intercal,
\end{align*}
and the external magnetic field $\bs H_e = 0$. We show in Figure \ref{Fig:ener} the discrete energy curve at the spatial temporal meshes with $K = 1/16$ and $\dt = 1/16,~1/32,~1/64$.
\end{example}

From Tables \ref{tab:ex-1} and \ref{tab:ex-2} and Figure \ref{Fig:ener}, we have the following observations:
\begin{itemize}
  \item The errors in $L^2$ norm for $\bs u$ and $\tilde p$ have the second order of convergence rates, which is better than the theoretical prediction, maybe because  the finite element spaces for   $\bs u$ and $\tilde p$ contain the piecewise polynomials of degree up to $1$. 
  \item The errors in $H^1$ semi-norm for $\bs u$ and $\bs\omega$, in  $H(\div)$ norms for $\bs H$ and $\bs m$, and in $L^2$ norm for $z$, $k$ and $\phi$, all have the first (optimal) order rates.
  \item When there is no source term, i.e. the external magnetic field $\bs H_e = \bs 0$,   the   discrete energy $\tilde{\mathcal E}_h$ decays with time, which is conformable to  the theoretical prediction  in Theorem  \ref{the:ener-dis-f}. 
\end{itemize}




%
%
%
%




%
\bibliographystyle{cas-model2-names}

\bibliography{cas-refs}



\end{document}